\documentclass[a4paper,12pt,twoside]{article}

\usepackage[top=2.2cm, bottom=2.8cm, left=2.2cm, right=2.2cm]{geometry}
\usepackage[pdftex]{graphicx}
\usepackage[pstarrows]{pict2e}
\usepackage{amsfonts,amstext,amscd,bezier,amsthm,amssymb}
\usepackage[centertags]{amsmath}
\usepackage{mathtools}
\usepackage[integrals]{wasysym}
\usepackage[nice]{nicefrac}
\usepackage{stmaryrd}
\usepackage[mathcal]{euscript}
\usepackage[latin1]{inputenc}
\usepackage[T1]{fontenc}
\usepackage[english]{babel}
\usepackage{cite}
\usepackage{mathptmx}
\usepackage{color}
  \definecolor{gray}{rgb}{0.65,0.65,0.65}
\usepackage{enumitem}
\usepackage{anyfontsize}

\newcommand{\theoname}{Theorem}
\newcommand{\lemmname}{Lemma}
\newcommand{\coroname}{Corollary}
\newcommand{\propname}{Proposition}
\newcommand{\definame}{Definition}
\newcommand{\hyponame}{Hypothesis}
\newcommand{\remkname}{Remark}
\newcommand{\explname}{Example}
\newcommand{\prfname}{Proof}

\theoremstyle{plain}
\newtheorem{theo}{\theoname}[section]
\newtheorem{lemm}[theo]{\lemmname}
\newtheorem{coro}[theo]{\coroname}
\newtheorem{prop}[theo]{\propname}
\theoremstyle{definition}
\newtheorem{defi}[theo]{\definame}
\newtheorem{hypo}[theo]{\hyponame}
\newtheorem{remk}[theo]{\remkname}
\newenvironment{expl}%
  {\refstepcounter{theo}%
    \begin{list}{}{%
    \setlength{\topsep}{0pt}%
    \setlength{\leftmargin}{0pt}%
    \setlength{\rightmargin}{0pt}%
    \setlength{\listparindent}{\parindent}%
    \setlength{\itemindent}{0pt}%
    \setlength{\parsep}{\parskip}}%
    \item[]{\bf \explname\ \thetheo. }}%
  {\hspace*{\fill} $\square$ \end{list} \medskip}
\newenvironment{prf}[1][\prfname]%
  {\begin{list}{}{%
    \setlength{\topsep}{0pt}%
    \setlength{\leftmargin}{0pt}%
    \setlength{\rightmargin}{0pt}%
    \setlength{\listparindent}{\parindent}%
    \setlength{\itemindent}{0pt}%
    \setlength{\parsep}{\parskip}}%
    \item[]{\bf #1. }}%
  {\hspace*{\fill} $\blacksquare$ \end{list} \medskip}

\DeclareMathOperator{\id}{Id}

\newcommand{\suchthat}{\;|\:}
\newcommand{\middlesuchthat}{\;\middle|\:}
\newcommand{\transp}{{\mathrm{T}}}
\newcommand{\R}{\mathbb{R}}

\newcommand{\norm}[1]{\left\lVert #1\right\lVert}
\newcommand{\abs}[1]{\left\lvert #1\right \rvert}

\newcommand{\Bigabs}[1]{\Bigl\lvert #1\Bigr \rvert}

\newcommand{\floor}[1]{\left\lfloor #1 \right\rfloor}
\newcommand{\ceil}[1]{\left\lceil #1 \right\rceil}

\newcommand{\scalprod}[2]{\left\langle #1, #2 \right\rangle}

\numberwithin{figure}{section}
\numberwithin{equation}{section}

\begin{document}

\makeatletter
\renewcommand{\@makecaption}[2]{%
  \vspace*{\abovecaptionskip}%
  \centering {\sc #1:} {\small\it #2}%
  \vspace*{\belowcaptionskip}%
}
\makeatother

\theoremstyle{plain}
\newtheorem{step}{Step}


\title{Persistently damped transport on a network of circles\thanks{This research was partially supported by the iCODE Institute, research project of the IDEX Paris-Saclay, and by the Hadamard Mathematics LabEx (LMH) through the grant number ANR-11-LABX-0056-LMH in the ``Programme des Investissements d'Avenir''. \newline\hspace*{\parindent} \textbf{Keywords:} transport equation, persistent excitation, exponential stability, multistructures \newline\hspace*{\parindent} \textbf{Mathematics Subject Classification:} 35R02, 35B35, 35C05, 35L40, 93C20}}
\author{
Yacine Chitour\thanks{Laboratoire des Signaux et Syst\`emes, Sup\'elec, and Universit\'e Paris Sud, Orsay, France.} ,
Guilherme Mazanti\thanks{CMAP, \'Ecole Polytechnique, Palaiseau, France.} {} \thanks{Team GECO, Inria Saclay, Palaiseau, France.} ,
Mario Sigalotti\footnotemark[\value{footnote}] \addtocounter{footnote}{-1} \footnotemark[\value{footnote}]
}
\maketitle

\begin{abstract}
In this paper we address the exponential stability of a system of transport equations with intermittent damping on a network of $N \geq 2$ circles intersecting at a single point $O$. The $N$ equations are coupled through a linear mixing of their values at $O$, described by a matrix $M$. The activity of the intermittent damping is determined by persistently exciting signals, all belonging to a fixed class. The main result is that, under suitable hypotheses on $M$ and on the rationality of the ratios between the lengths of the circles, such a system is exponentially stable, uniformly with respect to the persistently exciting signals. The proof relies on an explicit formula for the solutions of this system, which allows one to track down the effects of the intermittent damping.
\end{abstract}

\tableofcontents


\section{Introduction}

Consider the following system of $N \geq 2$ coupled transport equations,
\begin{equation}
\left\{
\begin{aligned}
 & \partial_t u_i(t, x) + \partial_x u_i(t, x) + \alpha_i(t) \chi_i(x) u_i(t, x) = 0, & & t \ge 0,\; x \in [0, L_i],\; 1\le i \le N_d, \\
 & \partial_t u_i(t, x) + \partial_x u_i(t, x) = 0, & \quad & t \ge 0,\; x \in [0, L_i],\; N_d+1\le i \le N, \\
 & u_i(t, 0) = \sum_{j = 1}^N m_{ij} u_j(t, L_j), & & t \ge 0,\; 1\le i \le N, \\
 & u_i(0, x) = u_{i, 0}(x), & & x \in [0, L_i],\; 1\le i \le N. 
\end{aligned}
\right.
\label{IntroTransport}
\end{equation}
For $1 \leq i \leq N$, the corresponding transport equation is defined in the space domain $[0, L_i]$ with $L_i > 0$. The integer $N_d$ denotes the number of equations with a damping term. For $1 \leq i \leq N_d$, the activity of the damping of the $i$-th equation in space is determined by the function $\chi_i$, which is assumed to be the characteristic function of an interval $[a_i, b_i] \subset [0, L_i]$ with $a_i < b_i$, whereas its activity in time is determined by the function $\alpha_i$, which is assumed to be a signal in $L^\infty(\mathbb R, [0, 1])$. The coupling between the $N$ transport equations is determined by the coefficients $m_{ij} \in \mathbb R$ for $1\le i,j \le N$. The goal of the paper consists in studying the stability properties of \eqref{IntroTransport} when the signals $\alpha_i$ are persistently exciting, as defined below.

\begin{defi}
Let $T$, $\mu$ be two positive constants with $T \geq \mu > 0$. A function $\alpha \in L^\infty(\mathbb R, [0, 1])$ is said to be a \emph{$(T, \mu)$-persistently exciting signal} if, for every $t \in \mathbb R$, one has
\begin{equation}
\int_t^{t + T} \alpha(s) ds \geq \mu.
\label{IntroCondPE}
\end{equation}
The set of all $(T, \mu)$-persistently exciting signals is denoted by $\mathcal G(T, \mu)$. 
\end{defi}

\begin{figure}[ht]

\centering

\setlength{\unitlength}{1cm}

\begin{picture}(6, 6)(-3, -3)

\put(0, 0){\circle*{0.2}}

\linethickness{0.02\unitlength}
\moveto(0, 0)
\curveto(-1.5,  0.375)(-3   ,  1.5  )(-3   , -0.75 )
\curveto(-3  , -0.75 )(-3   , -2.025)(-2.25, -2.025)
\curveto(-1.5, -2.025)(-0.75, -0.525)( 0   , 0     )
\strokepath
\moveto(0, 0)
\curveto(-2.25, -3)(-0.75, -3)(0, -3)
\curveto( 0.75, -3)( 2.25, -3)(0,  0)
\strokepath
\moveto(0, 0)
\curveto(1.5,  0.09375)(3   ,  0.75 )(3   , -0.75)
\curveto(3  , -0.75   )(3   , -2.25 )(2.25, -2.25)
\curveto(1.5, -2.25   )(0.75, -0.525)(0   ,  0   )
\strokepath
\moveto(0, 0)
\curveto(1.5, 0.375)(3, 0.75)(3, 1.875)
\circlearc{1.875}{1.875}{1.125}{0}{90}
\lineto(1.125, 3)
\curveto(0, 3)(0, 1.125)(0, 0)
\strokepath
\moveto(0, 0)
\curveto(-0.1875, 0.75)( 0   , 3    )(-1.5, 3   )
\curveto(-2.25  , 3   )(-3   , 3    )(-3  , 2.25)
\curveto(-3     , 1.5 )(-0.75, 0.375)( 0  , 0   )
\strokepath

\linethickness{0.07\unitlength}
\put( 2.5575 ,  2.77875){\vector(-1,  1  ){0}}
\put(-2.8125 ,  2.76   ){\vector(-2, -1.1){0}}
\put(-3      , -0.825  ){\vector( 0, -1  ){0}}
\put( 0.15   , -3      ){\vector( 1,  0  ){0}}
\put( 2.90625, -1.5    ){\vector( 1,  3.1){0}}

\linethickness{0.07\unitlength}%
\textcolor{red}{\put(-2.616, -1.8924){\Line(-0.06476484201, -0.076193931776)(0.06476484201, 0.076193931776)}}%
\textcolor{red}{\cbezier(-2.616, -1.8924)(-2.52, -1.974)(-2.4, -2.025)(-2.25, -2.025)}%
\textcolor{red}{\put(-1.8, -1.8648){\Line(0.055546753394, -0.0831538224461)(-0.055546753394, 0.0831538224461)}}%
\textcolor{red}{\cbezier(-2.25, -2.025)(-2.1, -2.025)(-1.95, -1.965)(-1.8, -1.8648)}%
\textcolor{red}{\put(0.756, -2.919){\Line(0.0321582426358, -0.0946881588721)(-0.0321582426358, 0.0946881588721)}}%
\textcolor{red}{\put(1.0546875, -1.734375){\Line(0.0932004671541, 0.0362446261155)(-0.0932004671541, -0.0362446261155)}}%
\textcolor{red}{\cbezier(0.756, -2.919)(1.11375, -2.7975)(1.35, -2.49375)(1.0546875, -1.734375)}%
\textcolor{red}{\put(1.704, 0.2085){\Line(0.0085989834578, -0.099629601442)(-0.0085989834578, 0.099629601442)}}%
\textcolor{red}{\put(2.0625, 0.2265625){\Line(0.0020828813682, 0.0999783056728)(-0.0020828813682, -0.0999783056728)}}%
\textcolor{red}{\cbezier(1.704, 0.2085)(1.83, 0.219375)(1.95, 0.225)(2.0625, 0.22265625)}%

\put(-1.8, -0.75 ){\makebox(0, 0){$C_1$}}
\put( 0  , -2    ){\makebox(0, 0){$C_2$}}
\put( 2  , -1    ){\makebox(0, 0){$C_3$}}
\put( 1.5,  1.5  ){\makebox(0, 0){$C_4$}}
\put(-1.5,  1.875){\makebox(0, 0){$C_5$}}
\put( 0.1,  0.2  ){\makebox(0, 0)[bl]{$O$}}
\end{picture}
\caption{Network corresponding to $N=5$ and $N_d=3$.
}
\label{FigIntroTransport}
\end{figure}
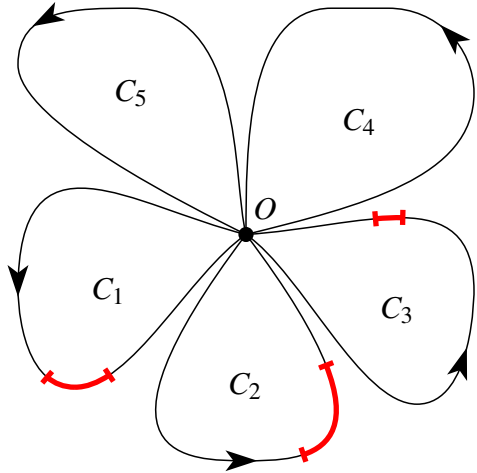

System~\eqref{IntroTransport} is a system of $N$ transport equations defined on intervals $[0, L_i]$, $1\le i \le N$, which may be identified with circles $C_1, C_2, \dotsc, C_N$ of respective lengths $L_1, L_2, \dotsc, L_N$. Moreover, there exists a point $O$ such that any two distinct circles only intersect at $O$ (see Figure~\ref{FigIntroTransport}). The transmission condition at $O$ can be written as
\begin{equation}
\begin{pmatrix}
u_1(t, 0) \\
u_2(t, 0) \\
\vdots \\
u_N(t, 0) \\
\end{pmatrix} = M
\begin{pmatrix}
u_1(t, L_1) \\
u_2(t, L_2) \\
\vdots \\
u_N(t, L_N) \\
\end{pmatrix},
\label{TransmissionMatrix}
\end{equation}
where $M = (m_{ij})_{1 \leq i, j \leq N}$ is called the \emph{transmission matrix} of the system. The topology of the network considered in this paper is star-shaped with respect to the point $O$. Note that any other network configuration falls into the present framework by a suitable choice of transition matrix $M$, namely, the fact that two circles $C_i$ and $C_j$ are not inward-outward adjacent translates to $m_{ij} = 0$. For $1 \leq i \leq N_d$, the transport equation of $u_i$ is damped on the support $[a_i, b_i]$ of $\chi_i$, represented in Figure~\ref{FigIntroTransport}. The damping is subject to the signal $\alpha_i$, which can be zero on certain time intervals. When all the $\alpha_i$ take their values in $\{0, 1\}$, \eqref{IntroTransport} can be seen as a switched system, where the switching signal $\alpha_i$ controls the damping action on the interval $[a_i, b_i]$ of the circle $C_i$.

Switching occurs in several control applications, motivating the study of systems with switched or intermittent actuators (see \cite{Gugat2008Optimal, Gugat2010Stars, Hante2009Modeling, Liberzon2003Switching, Margaliot2006Stability, Shorten2007Stability} and references therein). In this context, the activity of the actuator is guaranteed by appropriate conditions, for instance existence of a positive dwell-time or average dwell-time \cite{Liberzon2003Switching}. In this paper we rely instead on the integral condition \eqref{IntroCondPE} to guarantee the damping activity. This condition finds its origin in problems of identification and adaptive control \cite{Anderson1986Stability, Andersson2002Degenerate, Anderson1977Exponential}, where it is used in a more general form as a necessary and sufficient condition for the global exponential stability of some linear time-dependent systems. Linear control systems of the kind
\[\dot x(t) = A x(t) + \alpha(t) B u(t),\]
where a persistently exciting signal $\alpha$ determines the activity of the control $u$, are usually called \emph{persistently excited systems}. They have been considered in the finite-dimensional setting in \cite{Chitour2010Stabilization, Mazanti2013Stabilization, Mazanti2014Stabilization, Chaillet2008Uniform, Chaillet2007Towards, Chitour2014Growth}, dealing mostly with problems concerning stabilizability by a linear feedback law (see \cite{Chitour2013Stabilization} for a thorough presentation of persistently excited systems). In such systems, the persistently exciting signal $\alpha$ is a convenient tool to model several phenomena, such as failures in links between systems, resource allocation, or other internal or external processes that affect control efficiency. Further discussion on the practical interest of persistently excited systems can be found, e.g., in \cite{Loria2005PE}. Its infinite-dimensional counterpart is much less present in the literature, due to the fact that finite-dimensional results cannot be straightforwardly generalized (see \cite{Hante2012Conditions}, where illustrative examples are provided together with sufficient conditions ensuring stability).

System~\eqref{IntroTransport} is a ``toy model'' to study infinite-dimensional systems under persistent excitation. It is a simple case of a \emph{multi-body structure}. These type of problems model strings, membranes, plates, by partial differential equations defined on several coupled domains. They are an active research subject attracting much interest due to both its applications and the complexity of its analysis (see \cite{Ali2001Partial, Bastin2007Lyapunov, Lagnese1994Modeling, Lumer1980Connecting, Alabau2013Indirect, Bressan2014Flows} and references therein).

Particular multi-body systems that attracted much interest in recent years are  networks of strings \cite{Dager2006Wave, Gugat2010Stars, Zuazua2013Control, Valein2009Stabilization}. Their constituent elements are one-dimensional vibrating strings distributed along a planar graph. By decomposing each one-dimensional wave equation intro traveling waves according to D'Alembert decomposition, one can replace an edge of the graph by a pair of oriented edges and consider the transport equation in each edge. Hence, when \eqref{IntroTransport} is undamped (i.e., when $N_d=0$), it actually represents the D'Alembert decomposition of a star-shaped network of strings. The damping term in \eqref{IntroTransport} does not come from the above decomposition of the wave equation and thus the results of this paper cannot be directly applied to networks of strings. 

This paper addresses the issue of exponential stability of \eqref{IntroTransport}, uniformly with respect to the signals $\alpha_i$ in a class $\mathcal G(T, \mu)$: given $T \geq \mu > 0$, is system \eqref{IntroTransport}  uniformly exponentially stable with respect to $\alpha_i \in \mathcal G(T, \mu)$, $1\le i \le N_d$? The answer clearly depends on the transmission matrix $M$, since this matrix can amplify or reduce the solutions when they pass through $O$, as well as on the rationality of the ratios $L_i/L_j$, since periodic solutions may exist when they are rational (see Sections~\ref{SecExamples} and \ref{SecHypo} below). The main result of this paper is the following.

\begin{theo}
\label{MainTheoIntro}
Suppose that $N\ge 2$, $N_d \geq 1$, $\abs{M}_{\ell^1} \leq 1$, $m_{ij} \not = 0$ for every $1\le i,j \le N$, and that there exist $1\le i_*,j_* \le N$ such that ${L_{i_*}}/{L_{j_*}} \notin \mathbb Q$. Then, for every $T \geq \mu > 0$, there exist $C, \gamma > 0$ such that, for every $p \in [1, +\infty]$, every initial condition $u_{i, 0} \in L^p(0, L_i)$, $1\le i \le N$, and every choice of signals $\alpha_i \in \mathcal G(T, \mu)$, $1\le i \le N_d$, the corresponding solution of \eqref{IntroTransport} satisfies
\[\sum_{i=1}^N \norm{u_i(t)}_{L^p(0, L_i)} \leq C e^{-\gamma t} \sum_{i=1}^N \norm{u_{i, 0}}_{L^p(0, L_i)}, \qquad \forall t \geq 0.\]
\end{theo}

Our argument is based on explicit formulas for the solutions of system \eqref{IntroTransport}, which allow one to efficiently track down the effects of the persistency of the damping. This approach can be worked out since system \eqref{IntroTransport} consists of constant-speed transport equations with local damping. On the other hand, the usual techniques from PDE control, such as Carleman estimates, spectral criteria, Ingham estimates or microlocal analysis, do not seem well-adapted here, since they do not allow to handle the effects due to the time-dependency induced by the persistently exciting  signals $\alpha_i$. Extensions of our result to the case of state-dependent speed of transport and non-local damping would probably require more refined techniques.

The idea of relying on explicit representations for solutions of \eqref{IntroTransport} to address control and identification issues has already been used in \cite{Suzuki2013Analysis, Gugat2012Contamination}. Note, however, that, in these two references, rational dependence assumptions were necessary to derive tractable explicit formulas, which is not the case in this paper.

The paper is organized as follows. In Section~\ref{SecNotations}, we give the main notations and definitions used through this paper, discuss the well-posedness of \eqref{IntroTransport} and explain the role of the hypotheses of Theorem~\ref{MainTheoIntro}. Section~\ref{SecExplicit} provides the explicit formula for the solutions of \eqref{IntroTransport}, first in the undamped case, where the notations are simpler and the formulas easier to write, and then in the general case. Our main result is proved in Section~\ref{SecConvergence}, where we study the asymptotic behavior of coefficients appearing in the explicit solution obtained in Section~\ref{SecExplicit}. We finally collect in a series of appendices various technical results used in the paper.


\section{Notations, definitions and preliminary facts}
\label{SecNotations}

In this paper, we denote by $\mathbb Z$ the set of all integers, $\mathbb N = \{0, 1, 2, 3, \dotsc\}$ the set of nonnegative integers, $\mathbb N^\ast = \{1, 2, 3, \dotsc\}$ the set of positive integers, $\mathbb Q$ the set of rational numbers, $\mathbb R$ the set of real numbers, $\mathbb R_+ = \left[0, +\infty\right)$ the set of nonnegative real numbers and $\mathbb R_+^\ast = (0, +\infty)$ the set of positive real numbers. For $a, b \in \mathbb R$, let $\llbracket a, b \rrbracket = [a, b] \cap \mathbb Z$, with the convention that $[a, b] = \emptyset$ if $a > b$.

The set of $d \times m$ matrices with real coefficients is denoted by $\mathcal M_{d, m}(\mathbb R)$, or simply by $\mathcal M_d(\mathbb R)$ when $d = m$. As usual, we identify column matrices in $\mathcal M_{d, 1}(\mathbb R)$ with vectors in $\mathbb R^d$. The identity matrix in $\mathcal M_d(\mathbb R)$ is denoted by $\id_d$. For $p \in [1, +\infty]$, $\abs{\cdot}_{\ell^p}$ indicates both the $\ell^p$ norm of a vector of $\mathbb R^d$ and the induced matrix norm of a linear map on $\mathbb R^d$.

All Banach and Hilbert spaces considered are supposed to be real. The scalar product between two elements $u, v$ on a Hilbert space $\mathsf X$ is denoted by $\scalprod{u}{v}_{\mathsf X}$, and the norm of an element $u$ on a Banach space $\mathsf X$ is denoted by $\norm{u}_{\mathsf X}$. The indices $\mathsf X$ will be dropped from the previous notations when the context is clear. 

We shall refer to linear operators in a Banach space $\mathsf X$ simply as {\em operators}. The domain of an operator $T$ on $\mathsf X$ is denoted by $D(T)$. The set of all bounded operators from a Banach space $\mathsf X$ to a Banach space $\mathsf Y$ is denoted by $\mathcal L(\mathsf X, \mathsf Y)$ and is endowed with its usual norm $\norm{\cdot}_{\mathcal L(\mathsf X, \mathsf Y)}$. The notation $\mathcal L(\mathsf X)$ is used for $\mathcal L(\mathsf X, \mathsf X)$.

For two topological spaces $X$ and $Y$, $\mathcal C^0(X, Y)$ denotes the set of continuous functions from $X$ to $Y$. If $I$ is an interval of $\mathbb R$, $X$ is a Banach space and $k \in \mathbb N$, $\mathcal C^k(I,X)$ denotes the set of $k$-times differentiable $X$-valued functions defined on $I$, and $\mathcal C^k_c(I,X)$ is the subset of $\mathcal C^k(I,X)$ of  compactly supported functions. When $X=\mathbb R$ we omit it from the notation. 

For $x \in \mathbb R$, $\floor{x} \in \mathbb Z$ (resp. $\ceil{x}$) denotes the greatest (resp. smallest) integer $k \in \mathbb Z$ such that $k \leq x$ (resp. $k \geq x$). For $y > 0$, we denote by $\{x\}_y$ the number $\{x\}_y = x - \floor{\nicefrac{x}{y}}y$. For every $n,m \in \mathbb N$ with $m\leq n$, we use the binomial coefficient notation $\binom{n}{m}= \frac{n!}{m! (n-m)!}.$ We use $\# F$ and $\delta_{ij}$ to denote, respectively, the cardinality of a finite set $F$ and the Kronecker symbol of $i,j\in \mathbb Z$.

We refer to System~\eqref{IntroTransport} as being {\em undamped} by setting  $\alpha_i \equiv 0$ for every $i \in \llbracket 1, N_d\rrbracket$, in which case it is written as
\begin{equation}
\left\{
\begin{aligned}
 & \partial_t u_i(t, x) + \partial_x u_i(t, x) = 0, & \quad & t \in \mathbb R_+,\; x \in [0, L_i],\; i \in \llbracket 1, N\rrbracket, \\
 & u_i(t, 0) = \sum_{j = 1}^N m_{ij} u_j(t, L_j), & & t \in \mathbb R_+,\; i \in \llbracket 1, N\rrbracket, \\
 & u_i(0, x) = u_{i, 0}(x), & & x \in [0, L_i],\; i \in \llbracket 1, N\rrbracket.
\end{aligned}
\right.
\label{Undamped}
\end{equation}
We say that System~\eqref{IntroTransport} has an {\em always active damping} if $\alpha_i \equiv 1$ for every $i \in \llbracket 1, N_d \rrbracket$, in which case it becomes
\begin{equation}
\left\{
\begin{aligned}
 & \partial_t u_i(t, x) + \partial_x u_i(t, x) + \chi_i(x) u_i(t, x) = 0, & & t \in \mathbb R_+,\; x \in [0, L_i],\; i \in \llbracket 1, N_d \rrbracket, \\
 & \partial_t u_i(t, x) + \partial_x u_i(t, x) = 0, & \quad & t \in \mathbb R_+,\; x \in [0, L_i],\; i \in \llbracket N_d + 1, N\rrbracket, \\
 & u_i(t, 0) = \sum_{j = 1}^N m_{ij} u_j(t, L_j), & & t \in \mathbb R_+,\; i \in \llbracket 1, N\rrbracket, \\
 & u_i(0, x) = u_{i, 0}(x), & & x \in [0, L_i],\; i \in \llbracket 1, N\rrbracket.
\end{aligned}
\right.
\label{AlwaysActive}
\end{equation}
The general case of \eqref{IntroTransport} can be written as 
\begin{equation}
\left\{
\begin{aligned}
 & \partial_t u_i(t, x) + \partial_x u_i(t, x)+ \alpha_i(t)\chi_i(t)u_i(t,x)= 0, & \quad & t \in \mathbb R_+,\; x \in [0, L_i],\; i \in \llbracket 1, N\rrbracket, \\
 & u_i(t, 0) = \sum_{j = 1}^N m_{ij} u_j(t, L_j), & & t \in \mathbb R_+,\; i \in \llbracket 1, N\rrbracket, \\
 & u_i(0, x) = u_{i, 0}(x), & & x \in [0, L_i],\; i \in \llbracket 1, N\rrbracket,
\end{aligned}
\right.
\label{damped0}
\end{equation}
with the convention that $\alpha_i\equiv 1$ and $a_i=b_i=L_i$ for $i \in \llbracket N_d+1, N\rrbracket$, implying that $\chi_i= 0$ almost everywhere in $[0,L_i]$ for $i \in \llbracket N_d+1, N\rrbracket$. In the case where $\alpha_1,\dots,\alpha_{N_d}$ belong to a class $\mathcal G(T, \mu)$ for the same fixed $T \geq \mu > 0$, System~\eqref{damped0} is referred to as a {\em persistently damped system}.

\begin{remk}
Assuming that 
all the persistently exciting signals $\alpha_1,\dotsc, \alpha_{N_d}$,  in \eqref{damped0} belong to the class $\mathcal G(T, \mu)$, with the same constants $T \geq \mu > 0$,  is not  actually a restriction. Indeed, if $\alpha_i \in \mathcal G(T_i, \mu_i)$ with $T_i \geq \mu_i > 0$ for $i \in \llbracket 1, N_d \rrbracket$, then we clearly have, for every $i \in \llbracket 1, N_d \rrbracket$, $\alpha_i \in \mathcal G(T, \mu)$ with $T = \max_{i \in \llbracket 1, N_d \rrbracket} T_i$ and $\mu = \min_{i \in \llbracket 1, N_d \rrbracket} \mu_i$.
\end{remk}


\subsection{Formulation and well-posedness of the Cauchy problem}
\label{SecWellPosed}

The goal of this section consists in providing a rigorous definition for a solution of \eqref{damped0} and in guaranteeing that, given any initial data, the required solution exists, is unique and depends continuously on the initial data. 

\begin{defi}
\label{DefABi}
Let $p \in \left[1, +\infty\right)$. We set $\mathsf X_p = \prod_{i=1}^N L^p (0, L_i)$, endowed with the norm $\norm{z}_{\mathsf X_p} = \left(\sum_{i=1}^N \norm{u_i}_{L^p(0, L_i)}^p\right)^{\nicefrac{1}{p}}$ for $z = (u_1, \dotsc, u_N) \in \mathsf X_p$.

We define the operator $A : D(A) \subset \mathsf X_p \to \mathsf X_p$ on its domain $D(A)$ by
\begin{equation}
\begin{gathered}
D(A) = \left\{(u_1, \dotsc, u_N) \in \prod_{i=1}^N W^{1, p}(0, L_i) \middlesuchthat \forall i \in \llbracket 1, N\rrbracket,\; u_i(0) = \sum_{j=1}^N m_{ij} u_j(L_j) \right\}, \\
A(u_1, \dotsc, u_N) = \left(-\frac{d u_1}{dx}, \dotsc, -\frac{d u_N}{dx}\right).
\end{gathered}
\label{DefA}
\end{equation}

For $i \in \llbracket 1, N_d\rrbracket$, we define the operator $B_i \in \mathcal L(\mathsf X_p)$ by
\[B_i(u_1, \dotsc, u_N) = (0, \dotsc, 0, -\chi_i u_i, 0, \dotsc, 0),\]
where the term $-\chi_i u_i$ is in the $i$-th position.
\end{defi}

\begin{remk}
\label{RemkP1Infty}
Even though Theorem~\ref{MainTheoIntro} is stated for every $p \in [1, +\infty]$, we restrict ourselves  in the sequel to the case $p \in \left[1, +\infty\right)$. The main reason for this is that, when $p = +\infty$, the domain $D(A)$ of the operator $A$ defined by \eqref{DefA} is not dense in $\prod_{i=1}^N L^\infty(0, L_i)$, and thus some of our arguments given for $p$ finite do not apply. However, once we prove Theorem~\ref{MainTheoIntro} for $p \in \left[1, +\infty\right)$, we  obtain the case $p = +\infty$ by suitable continuity arguments, as  detailed in Remark~\ref{RemkInfty}.
\end{remk}

With the operators $A$ and $B_i$ defined above, System~\eqref{damped0} can be written as
\begin{equation}
\left\{
\begin{aligned}
\dot z(t) & = A z(t) + \sum_{i=1}^{N_d} \alpha_i(t) B_i z(t), \\
z(0) & = z_0,
\end{aligned}
\right.
\label{TranspOperatorPE}
\end{equation}
with $z_0 = (u_{1, 0}, \dotsc, u_{N, 0})$ and $\alpha_1,\dots,\alpha_{N_d}\in L^\infty(\R,[0,1])$. The case of the undamped system \eqref{Undamped} can be written as
\begin{equation}
\left\{
\begin{aligned}
\dot z(t) & = A z(t), \\
z(0) & = z_0,
\end{aligned}
\right.
\label{TranspOperatorUndamped}
\end{equation}
and the system with an always active damping \eqref{AlwaysActive} becomes
\begin{equation}
\left\{
\begin{aligned}
\dot z(t) & = A z(t) + \sum_{i=1}^{N_d} B_i z(t), \\
z(0) & = z_0.
\end{aligned}
\right.
\label{TranspOperatorAlwaysActive}
\end{equation}

The well-posedness of \eqref{TranspOperatorPE} is established in the sense of the following theorem, whose proof is deferred in Appendix~\ref{AppendWellPosed}.

\begin{theo}
\label{TheoWellPosed}
Let $p \in \left[1, +\infty\right)$ and $\alpha_i \in L^\infty(\mathbb R, [0, 1])$ for $i \in \llbracket 1, N_d\rrbracket$. There exists a unique evolution family $\{T(t, s)\}_{t \geq s \geq 0}$ of bounded operators in $\mathsf X_p$ such that, for every $s \geq 0$ and $z_0 \in D(A)$, $t \mapsto z(t) = T(t, s) z_0$ is the unique continuous function such that $z(s) = z_0$, $z(t) \in D(A)$ for every $t \geq s$, $z$ is differentiable for almost every $t \geq s$, $\dot z \in L^\infty_{\text{\rm loc}}(\left[s, +\infty\right), \mathsf X_p)$, and $\dot z(t) = A z(t) + \sum_{i=1}^{N_d} \alpha_i(t) B_i z(t)$ for almost every $t \geq s$.
\end{theo}

The definition of an evolution family is recalled in Appendix~\ref{AppendWellPosed}. The function $z$ in Theorem~\ref{TheoWellPosed} is said to be a \emph{regular solution} of \eqref{TranspOperatorPE} with initial condition $z_0 \in D(A)$. When $z_0 \in \mathsf X_p \setminus D(A)$, the function $t \mapsto z(t) = T(t, s) z_0$ is still well-defined and continuous, and is said to be a \emph{mild solution} of \eqref{TranspOperatorPE}. We use the word \emph{solution} to refer to both regular and mild solutions, according to the context.

Theorem~\ref{TheoWellPosed} also provides solutions to \eqref{TranspOperatorUndamped} and \eqref{TranspOperatorAlwaysActive} as particular cases. Since these equations are time-independent, we can actually obtain more regular solutions, thanks to the fact that $A$ and $A + \sum_{i=1}^{N_d} B_i$ generate strongly continuous semigroups, as we detail in Appendix~\ref{AppendWellPosed}.

\begin{theo}
\label{TheoWellPosedUndamped}
Let $p \in \left[1, +\infty\right)$. The operators $A$ and $A + \sum_{i=1}^{N_d} B_i$ generate strongly continuous semigroups $\{e^{tA}\}_{t \geq 0}$ and $\{e^{t(A + \sum_{i=1}^{N_d} B_i)}\}_{t \geq 0}$. In particular, for every $z_0 \in D(A)$, the function $t \mapsto e^{tA} z_0$ is the unique function in $\mathcal C^0(\mathbb R_+, D(A)) \cap \mathcal C^1(\mathbb R_+, \mathsf X)$ satisfying \eqref{TranspOperatorUndamped} and the function $t \mapsto e^{t(A + \sum_{i=1}^{N_d} B_i)} z_0$ is the unique function in $\mathcal C^0(\mathbb R_+, D(A)) \cap \mathcal C^1(\mathbb R_+, \mathsf X)$ satisfying \eqref{TranspOperatorAlwaysActive}.
\end{theo}


\subsection{Some examples of asymptotic behavior}
\label{SecExamples}

It is useful to have in mind some illustrative examples of the asymptotic behaviors of \eqref{IntroTransport} under no damping, an always active damping and a persistent damping, respectively. 

\begin{expl}
\label{ExplOneCircle}
Consider the case of a single transport equation on a circle of length $L$,
\begin{equation}
\label{OneCircle}
\left\{
\begin{aligned}
 & \partial_t u(t, x) + \partial_x u(t, x) + \alpha(t) \chi(x) u(t, x) = 0, & & t \in \mathbb R_+,\; x \in [0, L], \\
 & u(t, 0) = u(t, L), & & t \in \mathbb R_+, \\
 & u(0, x) = u_0(x), & & x \in [0, L], \\
 & \alpha \in \mathcal G(T, \mu),
\end{aligned}
\right.
\end{equation}
where $\chi$ is the characteristic function of the interval $[a, b] \subset [0, L]$. This corresponds to \eqref{damped0} with persistent damping, $N_d = N = 1$, and $m_{11} = 1$. Due to the condition $u(t, 0) = u(t, L)$, it can be seen as a transport equation on a circle of length $L$.

When \eqref{OneCircle} is undamped, all its solutions are $L$-periodic. Indeed, for $u_0 \in \mathsf X_p = L^p(0, L)$, the corresponding solution of \eqref{OneCircle} is $u(t, x) = u_0(\{x - t\}_{L})$, where we recall that $\{x\}_y = x - \floor{\nicefrac{x}{y}} y$, and this function is clearly $L$-periodic.

When \eqref{OneCircle} has an always active damping, all its solutions converge exponentially to zero. Indeed,  every solution of \eqref{OneCircle} satisfies $u(t, x) = e^{-(b - a)} u(t - L, x)$ for every $x \in [0, L]$ and $t \geq L$, so that $\norm{u(t)}_{L^p(0, L)} = e^{-(b-a)} \norm{u(t - L)}_{L^p(0, L)}$. It is also clear that $\norm{u(t)}_{L^p(0, L)} \leq \norm{u_0}_{L^p(0, L)}$ for every $t \geq 0$, and so $\norm{u(t)}_{L^p(0, L)} \leq C e^{-\gamma t} \norm{u_0}_{L^p(0, L)}$ for every $t \geq 0$, with $\gamma = \frac{b-a}{L}$ and $C = e^{\gamma L}$.

When \eqref{OneCircle} has a persistent damping and the damping interval $[a, b]$ is a proper subset of $[0, L]$, there exist $T > \mu > 0$, a persistently exciting signal $\alpha \in \mathcal G(T, \mu)$ and a nontrivial initial condition $u_0 \in L^p(0, L)$ such that the corresponding solution of \eqref{OneCircle} is $L$-periodic. Indeed, suppose that $a = 0$ and $b < L$. Take $u_0 \in \mathcal C^\infty([0, L]) \setminus \{0\}$ such that the support of $u_0$ is contained in $[b, \frac{b+L}{2}]$. Take $\alpha \in L^\infty(\mathbb R, [0, 1])$ defined by
\[
\alpha(t) = 
\begin{dcases*}
1, & if $0 \leq \{t\}_L \leq \frac{L - b}{2}$, \\
0, & if $\{t\}_L > \frac{L - b}{2}$.
\end{dcases*}
\]
Then $\alpha \in \mathcal G\left(L, \frac{L - b}{2}\right)$ and one can easily verify that the corresponding solution $u(t, x)$ of \eqref{OneCircle} is equal to $u_0(\{x - t\}_L)$. Hence \eqref{OneCircle} admits a $L$-periodic solution. 
\end{expl}

Example~\ref{ExplOneCircle} shows that the asymptotic behavior of \eqref{damped0} can be different if the damping is always active or if it is submitted to a persistently exciting signal, and this is due to the fact that the support of the solution may not be in the damping interval $[a, b]$ when the damping is active. We now consider a second example showing that, when we have more than one circle, the rationality of the ratios ${L_i}/{L_j}$ for $i \not = j$ plays an important role in the asymptotic behavior.

\begin{expl}
\label{ExplTwoCircles}
Consider the case of System~\eqref{damped0} with persistent damping, $N = 2$, $N_d = 1$, and $m_{ij} = \nicefrac{1}{2}$ for $i, j \in \{1, 2\}$, i.e.,
\begin{equation}
\left\{
\begin{aligned}
 & \partial_t u_1(t, x) + \partial_x u_1(t, x) + \alpha(t) \chi(x) u_1(t, x) = 0, & & t \in \mathbb R_+,\; x \in [0, L_1], \\
 & \partial_t u_2(t, x) + \partial_x u_2(t, x) = 0, & \quad & t \in \mathbb R_+,\; x \in [0, L_2], \\
 & u_1(t, 0) = u_2(t, 0) = \frac{u_1(t, L_1) + u_2(t, L_2)}{2}, & & t \in \mathbb R_+, \\
 & u_i(0, x) = u_{i, 0}(x), & & x \in [0, L_i],\; i \in \{1, 2\}, \\
 & \alpha \in \mathcal G(T, \mu), & &
\end{aligned}
\right.
\label{TwoCircles}
\end{equation}
where $\chi$ is the characteristic function of the interval $[a, b] \subset [0, L_1]$. In order to simplify the discussion, let us fix $p=2$ and set $\mathsf X_2 = L^2(0, L_1) \times L^2(0, L_2)$.

When \eqref{TwoCircles} is undamped its asymptotic behavior depends on the rationality of the ratio ${L_1}/{L_2}$, as stated in the next theorem, which is proved in Appendix~\ref{AppendAsymptotic}.

\begin{theo}
\label{TheoAsympTwoCircles}
Consider \eqref{TwoCircles} with $\chi \equiv 0$.
\begin{enumerate}[label={\bf\roman*.}, ref={\roman*}]
\item If ${L_1}/{L_2} \notin \mathbb Q$, each solution converges to a constant function $(\lambda, \lambda) \in \mathsf X_2$ with $\lambda \in \mathbb R$.
\label{CaseIrrational}
\item If ${L_1}/{L_2} \in \mathbb Q$, there exists a non-constant periodic solution.
\label{CaseRational}
\end{enumerate}
\end{theo}

When \eqref{TwoCircles} has an always active damping all solutions converge exponentially to zero, independently of  the rationality of the ratio ${L_1}/{L_2}$, as it follows, for instance, from Remark~\ref{TheoConvNonPE}.

When \eqref{TwoCircles} has a persistent damping, the rationality of the ratio ${L_1}/{L_2}$ plays once again a role in the asymptotic behavior of the system: if ${L_1}/{L_2} \notin \mathbb Q$, all its solutions converge exponentially to zero, as it follows from our main result, Theorem~\ref{MainTheoIntro}. However, if ${L_1}/{L_2} \in \mathbb Q$ and the damping interval $[a, b]$ is small enough, there exist $T > \mu > 0$, a persistently exciting signal $\alpha \in \mathcal G(T, \mu)$ and a nontrivial initial condition $u_0 \in L^p(0, L)$ such that the corresponding solution of \eqref{TwoCircles} is periodic, as we show in Appendix~\ref{SecPEPeriodic}.
\end{expl}

Both in Example~\ref{ExplOneCircle} and in Example~\ref{ExplTwoCircles} in the case ${L_1}/{L_2} \in \mathbb Q$, the lack of exponential stability is illustrated by the existence of a periodic solution for the persistently damped system which is actually a solution to the undamped one for which a persistently exciting signal $\alpha$ inactivates the damping whenever the support of the solution passes through the damping interval. The heuristic of the proof of  Theorem~\ref{MainTheoIntro} is that the support of every initial condition  spreads with time and eventually covers the entire network. Hence every solution of the persistently excited system  eventually passes through a damping interval at a time where the damping is active.


\subsection{Discussion on the hypotheses of Theorem~\ref{MainTheoIntro}}
\label{SecHypo}

Recall that the two main assumptions of Theorem~\ref{MainTheoIntro} are the following. 

\begin{hypo}
\label{HypoIrrational}
There exist $i_*, j_* \in \llbracket 1, N\rrbracket$ such that ${L_{i_*}}/{L_{j_*}} \notin \mathbb Q$.
\end{hypo}

\begin{hypo}
\label{HypoM}
The matrix $M$ satisfies $\abs{M}_{\ell^1} \leq 1$ and $m_{ij} \not = 0$  for every $i, j \in \llbracket 1, N \rrbracket$.
\end{hypo}

At the light of Example~\ref{ExplTwoCircles}, one cannot expect exponential stability of \eqref{damped0} with persistent damping in general if ${L_i}/{L_j} \in \mathbb Q$ for every $i, j \in \llbracket 1, N\rrbracket$. This is why it is reasonable to make Hypothesis~\ref{HypoIrrational}.

Even though the well-posedness of \eqref{damped0} discussed in Section~\ref{SecWellPosed} and the explicit formula for its solutions given later in Section~\ref{SecExplicit} are obtained for every $M \in \mathcal M_N(\mathbb R)$, the asymptotic behavior of System~\eqref{damped0} clearly depends on the  choice of the matrix $M$, since this matrix determines the coupling among the $N$ transport equations. 

The hypothesis $\abs{M}_{\ell^1} \leq 1$ can be written as
\begin{equation}
\sum_{i=1}^{N} \abs{m_{ij}} \leq 1, \qquad \forall j \in \llbracket 1, N\rrbracket.
\label{CondNormL1}
\end{equation}
The coefficient $m_{ij}$ can be interpreted as the proportion of mass in the circle $C_j$ that goes to the circle $C_i$ as it passes the contact point $O$. Hence, \eqref{CondNormL1} states that, for every $j \in \llbracket 1, N\rrbracket$, the total mass arriving at the circles $C_i$ from the circle $C_j$ is less than or equal the total mass leaving the circle $C_j$, which means that the mass never increases while passing through the junction.

The hypothesis $m_{ij}\ne 0$ for all $i,j$ can be seen as a \emph{strong mixing} of the solutions at the junction. It is designed to avoid reducibility phenomena which may be an obstruction to uniform exponential stability. Consider for instance the case $M = \id_N$ with $N \geq 2$. Then \eqref{damped0} is reduced to $N$ uncoupled transport equations on circles, each of them of the form \eqref{OneCircle}. In that case, if there exists at least one index $i \in \llbracket 1, N\rrbracket$ such that $b_i - a_i < L_i$, then there exist solutions not converging to $0$ as $t \to +\infty$, even if there is damping, cf. Example~\ref{ExplOneCircle}.

\begin{remk}
Equation~\eqref{CondNormL1} is satisfied when $M$ is  \emph{left stochastic}, i.e., $m_{ij} \geq 0$ for every $i, j \in \llbracket 1, N \rrbracket$ and $\sum_{i=1}^N m_{ij} = 1$ for every $ j \in \llbracket 1, N \rrbracket$. Note that left stochasticity of $M$ is equivalent for the undamped system \eqref{Undamped} to the preservation of $\sum_{i=1}^N \int_0^{L_i} u_i(t, x) dx$ and monotonicity of the solutions with respect to the initial conditions. 
\end{remk}


\section{Explicit solution}
\label{SecExplicit}

This section provides a general formula for the explicit solution of \eqref{damped0}. We first prove our formula in Section~\ref{SecSolutionUndamped} in the simpler case of the undamped system \eqref{Undamped},  before turning to the general case in Section~\ref{SecGeneralFormula}. The coefficients appearing in the formula will be characterized in Section~\ref{SecCoefficients}.


\subsection{The undamped system}
\label{SecSolutionUndamped}

Remark that, in order to obtain an explicit formula for $u_i(t, x)$ for $i \in \llbracket 1, N \rrbracket$, $t \geq 0$ and $x \in [0, L_i]$, it suffices to obtain a formula for $u_i(t, 0)$ for $i \in \llbracket 1, N \rrbracket$ and $t \geq 0$.
Indeed, it is immediate to derive the following. 

\begin{lemm}
\label{PropOnlyXZeroMatters}
Let $(u_{1, 0}, \dotsc, u_{N, 0}) \in D(A)$ and let $(u_1, \dotsc, u_N) \in \mathcal C^0(\mathbb R_+, D(A)) \cap \mathcal C^1(\mathbb R_+, \mathsf X_p)$ be the corresponding solution of \eqref{Undamped}. Then, for every $i \in \llbracket 1, N \rrbracket$, $t \geq 0$, and $x \in [0, L_i]$, we have
\begin{equation}
u_i(t, x) = 
\begin{dcases*}
u_{i, 0}(x - t), & if $0 \leq t \leq x$, \\
u_i(t - x, 0), & if $t \geq x$.
\end{dcases*}
\label{OnlyXZeroMatters}
\end{equation}
\end{lemm}

In order to express $u_i(t, 0)$ in terms of the initial condition $(u_{1, 0}, \dotsc, u_{N, 0}) \in D(A)$, we need to introduce some notation.

\begin{defi}
\label{DefMathfrak}
\begin{enumerate}[label={\bf \roman*.}, ref={\roman*}]
\item We define $\mathfrak N = \mathbb N^N$ and, for $i \in \llbracket 1, N \rrbracket$,  $\mathfrak N_i = \mathbb N^{i-1} \times \{0\} \times \mathbb N^{N - i}$.
\item We write $\mathbf 0 = (0, 0, \dotsc, 0) \in \mathfrak N$ and, for every $j \in \llbracket 1, N\rrbracket$ and $\mathfrak n = (n_1, \dotsc, n_N) \in \mathfrak N$,  $\mathbf 1_j = (\delta_{ij})_{i = 1, \dotsc, N} \in \mathfrak N$ and $\hat{\mathfrak n}_j = (n_1, n_2, \dotsc, n_{j-1}, 0, n_{j+1}, \dotsc, n_N) = \mathfrak n - n_j \mathbf 1_j \in \mathfrak N_j$.
\item We define the function $L: \mathfrak N \to \mathbb R_+$ by
\[L(n_1, \dotsc, n_N) = \sum_{i=1}^N n_i L_i.\]
\end{enumerate}
\end{defi}

With these notations, the general formula for the solutions of \eqref{Undamped} can be written as follows.

\begin{theo}
\label{TheoSolExplicite}
Let $(u_{1, 0}, \dotsc, u_{N, 0}) \in D(A)$. The corresponding solution $(u_1, \dotsc, u_N)$ of \eqref{Undamped} is given by
\begin{equation}
u_i(t, x) = 
\begin{dcases*}
u_{i, 0}(x - t), & if $0 \leq t \leq x$, \\
u_i(t - x, 0), & if $t \geq x$,
\end{dcases*}
\label{ExplicitSolution}
\end{equation}
with
\begin{equation}
u_i(t, 0) = \sum_{j=1}^N \sum_{\substack{\mathfrak n \in \mathfrak N_j \\ L(\mathfrak n) \leq t}} \beta^{(i)}_{j, \mathfrak n + \floor{\frac{t - L(\mathfrak n)}{L_j}} \mathbf 1_j} u_{j, 0}\left(L_j - \left\{t - L(\mathfrak n)\right\}_{L_j}\right),
\label{ExplicitPhi}
\end{equation}
and where the coefficients $\beta^{(i)}_{j, \mathfrak n}$ are defined by the relations
\begin{subequations}
\label{Beta}
\begin{equation}
\beta^{(i)}_{j, \mathbf 0} = m_{ij}, \qquad i, j \in \llbracket 1, N\rrbracket,
\label{BetaInitial}
\end{equation}
and
\begin{equation}
\beta^{(i)}_{j, \mathfrak n} = \sum_{\substack{k=1 \\ n_k \geq 1}}^N m_{kj} \beta^{(i)}_{k, \mathfrak n - \mathbf 1_k}, \qquad i, j \in \llbracket 1, N\rrbracket,\; \mathfrak n \in \mathfrak N \setminus \{\mathbf 0\}.
\label{BetaRecurrence}
\end{equation}
\end{subequations}
\end{theo}
The above result follows by iterating Equations~\eqref{OnlyXZeroMatters} together with Equation~\eqref{TransmissionMatrix}. Indeed, using the notations of the theorem, one has for $i\in \llbracket 1, N\rrbracket$, 
\begin{equation}\label{toto1}
u_i(t,0)=\sum_{j=1}^Nm_{ij}u_j(t,L_j).
\end{equation}
According to Equation~\eqref{OnlyXZeroMatters}, each $u_j(t,L_j)$ is either equal to $u_{j,0}(L_j-t)$ or $u_j(t-L_j, 0)$ whether $t\leq L_j$ or not. In the latter case, we express $u_j(t-L_j, 0)$
by using Equation~\eqref{toto1} and we repeat the procedure a finite number of times until obtaining $u_i(t, 0)$ as a linear combination involving only evaluations of the initial condition at finitely many points on the circles. This yields Equation~\eqref{ExplicitPhi} with an explicit expression of both the coefficients of this linear combination and the points on the circles. 

The complete proof of Theorem~\ref{TheoSolExplicite} is provided in Appendix~\ref{AppendExplicit} and consists in verifying that the explicit formula given in the above statement is indeed the solution of \eqref{Undamped}. 

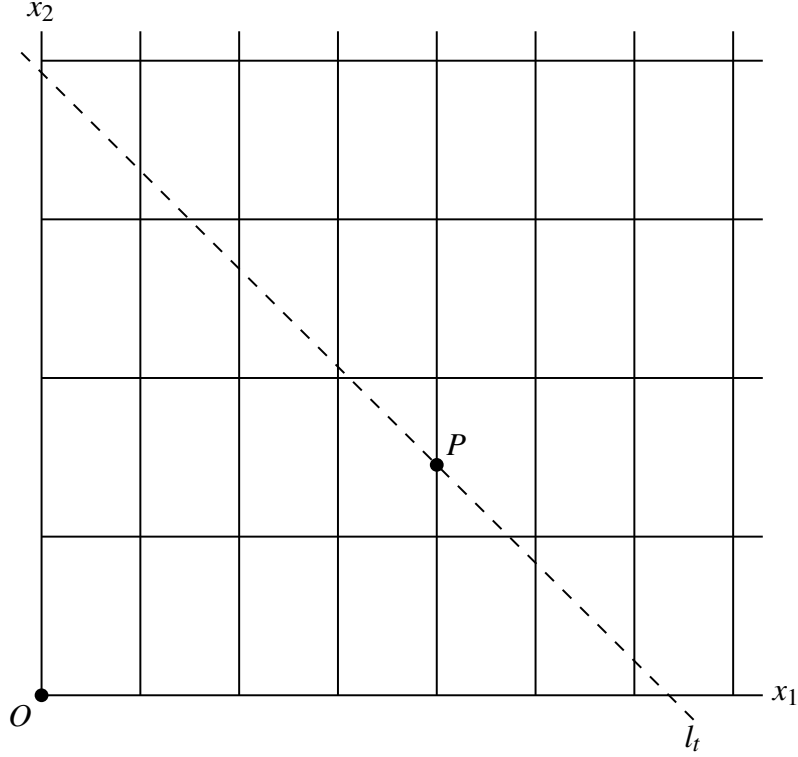
\begin{figure}[ht]
\centering

\setlength{\unitlength}{1.3cm}

\begin{picture}(8.6, 8.2)(-0.5, -0.7)

\linethickness{0.02\unitlength}

\multiput(0, 0)(0, 1.618034){5}{\Line(0, 0)(7.3, 0       )}
\multiput(0, 0)(1, 0       ){8}{\Line(0, 0)(0  , 6.772136)}

\put(0, 0       ){\circle*{0.14}}
\put(4, 2.35    ){\circle*{0.14}}

\linethickness{0.02\unitlength}
\multiput(6.6, -0.25)(-0.176777, 0.176777){39}{\line(-1, 1){0.0883883}}

\put(-0.1, -0.1    ){\makebox(0, 0)[tr]{$O$}}
\put( 6.6, -0.3    ){\makebox(0, 0)[t]{$l_t$}}
\put( 4.1, 2.45    ){\makebox(0, 0)[bl]{$P$}}
\put( 7.4, 0       ){\makebox(0, 0)[l]{$x_1$}}
\put( 0  , 6.872136){\makebox(0, 0)[b]{$x_2$}}
\end{picture}
\caption{Geometric construction for the explicit formula for the solution of \eqref{Undamped} in the case $N = 2$. }
\label{FigExplicite3}
\end{figure}

We next provide with Figure~\ref{FigExplicite3} a geometric interpretation of \eqref{ExplicitPhi} in the case $N=2$. The point $O$ is identified with the origin of the plane $(x_1,\ x_2)$ and the horizontal (resp. vertical) segments in the grid represented in  Figure~\ref{FigExplicite3} correspond to identical copies of the circle $C_1$ (resp. the circle $C_2$). 
The intersection of the dashed line $l_t:\ x_1+x_2=t$ and the grid exactly represents the set of points of the circles where the initial condition $(u_{1,0},u_{2,0})$ is evaluated in Equation~\eqref{ExplicitPhi}. Note that the coefficients in Equation~\eqref{ExplicitPhi} appearing in front of the evaluation of the   initial condition at  $P$ can be expressed as a sum of products of the $m_{ij}$'s, each product corresponding to a path on the grid between $P$ and $O$.


\subsection{Formula for the explicit solution in the general case}
\label{SecGeneralFormula}

We first notice that, as in Lemma~\ref{PropOnlyXZeroMatters}, it suffices to study $u_i(t, 0)$ for every $t \geq 0$ and $i \in \llbracket 1, N\rrbracket$ in order to obtain the whole solution $(u_1(t), \dotsc, u_N(t))$. Recall that by convention we have set $\alpha_i\equiv 1$ and 
$a_i=b_i$ (and thus $\chi_i= 0$ almost everywhere) for $i\in \llbracket  N_d+1,N \rrbracket$. 

\begin{prop}
\label{PropOnlyXZeroMattersPE}
Let $(u_{1, 0}, \dotsc, u_{N, 0}) \in D(A)$. Then the corresponding solution $(u_1, \dotsc, u_N)$ of \eqref{damped0} satisfies, for $i \in \llbracket 1, N \rrbracket$,
\begin{equation}
\label{FlowU1UiPE}
u_i(t, x)=\begin{dcases*}
u_{i, 0}(x - t) \exp\left({-\int_{[0,t]\cap [t-x+a_i,t-x+b_i]} \alpha_i(s) ds}\right), & if $0 \leq t \leq x$, \\
u_i(t - x, 0) \exp\left({-\int_{[0,t]\cap [t-x+a_i,t-x+b_i]} \alpha_i(s) ds}\right), & if $t \geq x$. \\
\end{dcases*}
\end{equation}
\end{prop}

\begin{prf}
Let $i \in \llbracket 1, N \rrbracket$. Equation~\eqref{FlowU1UiPE} is obtained by integrating the differential equation 
$$
\frac{d}{ds}u_i(t+s,x+s)=-\alpha_i(t+s)\chi_i(x+s)u_i(t+s,x+s),
$$
on the interval $[-t,0]$ if $t\leq x$ and on $[-x,0]$ if $t\geq x$. 
\end{prf}

Thanks to the fact that all the exponential decays appearing in \eqref{FlowU1UiPE} are upper bounded by $1$, one obtains trivially the following corollary.

\begin{coro}
\label{CoroBoundPE}
If $(u_1, \dotsc, u_N)$ is the solution of \eqref{damped0} with an initial condition $(u_{1, 0}, \dotsc, u_{N, 0})$, then, for $i \in \llbracket 1, N_d\rrbracket$, $u_i$ satisfies the estimate
\[
\abs{u_i(t, x)} \leq 
\begin{dcases*}
\abs{u_{i, 0} (x - t)}, & if $0 \leq t \leq x$, \\
\abs{u_i(t - x, 0)}, & if $t \geq x$.
\end{dcases*}
\]

For every $p \in [1, +\infty]$, $i \in \llbracket 1, N\rrbracket$, and $t \geq L_i$, we have
\[\norm{u_i(t, \cdot)}_{L^p(0, L_i)} \leq \norm{u_i(\cdot, 0)}_{L^p(t - L_i, t)},\]
with equality if $i \in \llbracket N_d + 1, N\rrbracket$.
\end{coro}

This corollary allows us to replace the spatial $L^p$-norm of $u_i$ at a given time $t$ by its $L^p$-norm in a time interval of length $L_i$ at the fixed position $x = 0$.

We can now write the explicit formula for the solutions of \eqref{damped0} using the notations from Definition~\ref{DefMathfrak}. The proof follows the same steps as that of Theorem~\ref{TheoSolExplicite}.

\begin{theo}
\label{TheoPEExplicit}
Let $(u_{1, 0}, \dotsc, u_{N, 0}) \in D(A)$. The corresponding solution $(u_1, \dotsc, u_N)$ of \eqref{damped0} is given by \eqref{FlowU1UiPE}, where $u_i(t, 0)$ is given for $t \geq 0$ by
\begin{equation}
u_i(t, 0) = \sum_{j=1}^N \sum_{\substack{\mathfrak n \in \mathfrak N_j \\ L(\mathfrak n) \leq t}} \vartheta^{(i)}_{j, \mathfrak n + \floor{\frac{t - L(\mathfrak n)}{L_j}} \mathbf 1_j, L_j - \left\{t - L(\mathfrak n)\right\}_{L_j}, t} u_{j, 0}\left(L_j - \left\{t - L(\mathfrak n)\right\}_{L_j}\right)
\label{PhiPEFinal}
\end{equation}
and the coefficients $\vartheta^{(i)}_{j, \mathfrak n, x, t}$ are defined for $i, j \in \llbracket 1, N\rrbracket$, $\mathfrak n \in \mathfrak N$, $x \in [0, L_j]$ and $t \in \R$ by
\begin{subequations}
\label{Vartheta}
\begin{equation}
\vartheta^{(i)}_{j, \mathfrak n, x, t} = \varepsilon_{j, \mathfrak n, x, t}\vartheta^{(i)}_{j, \mathfrak n, L_j, t} , 
\label{ThetaIX-PE}
\end{equation}
with
\begin{equation}
\varepsilon_{j, \mathfrak n, x, t} = 
\exp\left(-\int_{I_{j, \mathfrak n, x, t}}\alpha_j(s) ds
\right),
\label{EpsilonIX-PE}
\end{equation}
where $I_{j, \mathfrak n, x, t}=[t - L(\mathfrak n) - L_j +\max(x,a_j),t - L(\mathfrak n) - L_j+b_j]$, 
and 
\begin{align}
\vartheta^{(i)}_{j, \mathbf 0, L_j, t} &= m_{ij},
\label{VarthetaInitial}\\
\vartheta^{(i)}_{j, \mathfrak n, L_j, t} 
&= \sum_{\substack{k=1 \\ n_k \geq 1}}^N m_{kj} \vartheta^{(i)}_{k, \mathfrak n - \mathbf 1_k, 0, t}.
\label{VarthetaRecurrence}
\end{align}
\end{subequations}
\end{theo}

\begin{figure}[ht]
\centering

\setlength{\unitlength}{1.3cm}

\begin{picture}(8.6, 8.2)(-0.5, -0.7)

\linethickness{0.02\unitlength}

\multiput(0, 0)(0, 1.618034){5}{\Line(0, 0)(7.3, 0       )}
\multiput(0, 0)(1, 0       ){8}{\Line(0, 0)(0  , 6.772136)}

\linethickness{0.04\unitlength}
\multiput(0, 0)(0, 1.618034){5}{%
\multiput(0, 0)(1, 0       ){7}{%
\textcolor{red}{%
\Line(0.2,  0  )(0.6, 0  )%
\Line(0.2, -0.1)(0.2, 0.1)%
\Line(0.6, -0.1)(0.6, 0.1)%
}}}%

\put(0, 0       ){\circle*{0.14}}
\put(4, 2.35    ){\circle*{0.14}}

\linethickness{0.02\unitlength}
\multiput(6.6, -0.25)(-0.176777, 0.176777){39}{\line(-1, 1){0.0883883}}

\put(-0.1, -0.1    ){\makebox(0, 0)[tr]{$O$}}
\put( 6.6, -0.3    ){\makebox(0, 0)[t]{$l_t$}}
\put( 4.1, 2.45    ){\makebox(0, 0)[bl]{$P$}}
\put( 7.4, 0       ){\makebox(0, 0)[l]{$x_1$}}
\put( 0  , 6.872136){\makebox(0, 0)[b]{$x_2$}}
\end{picture}
\caption{Geometric interpretation of the explicit formula \eqref{PhiPEFinal} for the solution of \eqref{damped0} in the case $N = 2$.}
\label{FigTranspDamped}
\end{figure}
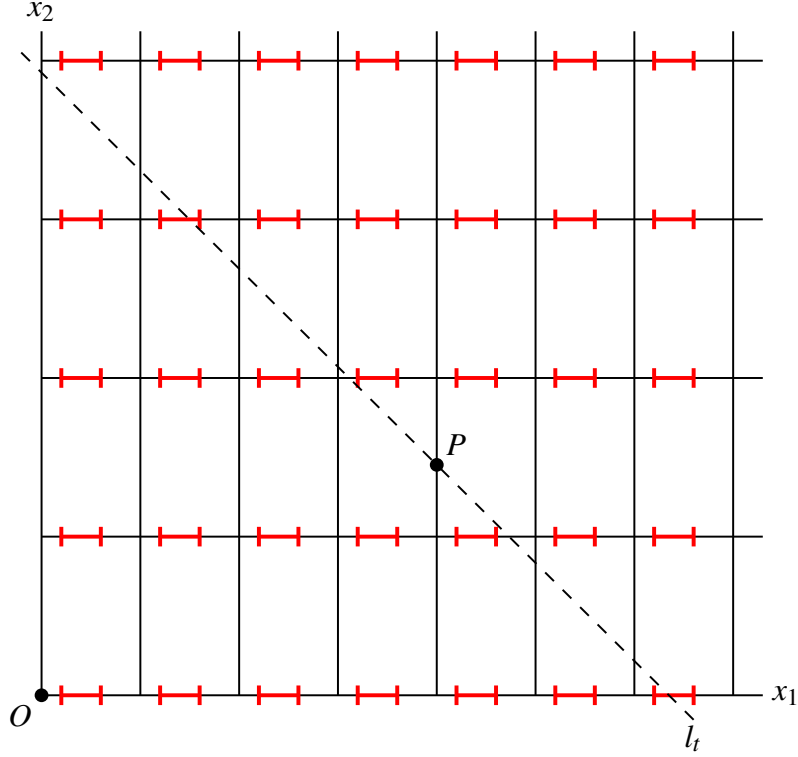

Let us provide a geometrical interpretation of the above theorem in the case $N=2$, $N_d=1$ with damping on the circle $C_1$. With respect to Figure~\ref{FigExplicite3}, Figure~\ref{FigTranspDamped} now includes segments corresponding to the intervals $[a_1, b_1]$ in $C_1$, on which the solution is damped. Similarly to \eqref{ExplicitPhi}, the new explicit formula~\eqref{PhiPEFinal} expresses $u_i(t, 0)$ as a linear combination involving only evaluations of the initial condition at finitely many points on the circles. The coefficients, which in \eqref{ExplicitPhi} were sum of products of the $m_{ij}$'s, each product corresponding to a path between $P$ and $O$, have now an analogous  expression, with the following modification: each factor of the original product is multiplied by an  additional term of the type $\varepsilon_{j, \mathfrak n, x, t}$, which takes into account the effect of the damping along the path. 

\begin{remk}
The explicit formula for the solution of the undamped equation \eqref{Undamped} given in Theorem~\ref{TheoSolExplicite} can be obtained as a particular case of Theorem~\ref{TheoPEExplicit} by setting $\alpha_j \equiv 0$ for every $j \in \llbracket 1, N_d\rrbracket$. 
 Similarly, we can obtain the explicit formula for the solution of \eqref{AlwaysActive} as a particular case of Theorem~\ref{TheoPEExplicit} by setting $\alpha_j \equiv 1$ for every $j \in \llbracket 1, N_d\rrbracket$, yielding $\varepsilon_{j,\mathfrak n,x,t}= e^{-\mathrm{meas}([a_j,b_j]\cap [x,b_j])}$.
\end{remk}

\begin{remk}\label{densityDAinXp}
In the general case $(u_{1, 0}, \dotsc, u_{N, 0}) \in \mathsf X_p$, the mild solution $(u_1, \dotsc, u_N)$ of \eqref{damped0} can still be characterized by \eqref{FlowU1UiPE} and \eqref{PhiPEFinal} (yielding an equality in $\mathsf X_p$ for every $t\geq 0$). This follows by a simple density argument of $D(A)$ in $\mathsf X_p$.
\end{remk}


\subsection{Recursive  formula for the coefficients}
\label{SecCoefficients}

We now wish to determine a recursive  formula with $K$ steps for the coefficients $\vartheta^{(i)}_{j, \mathfrak n, x, t}$ appearing in the expression of the explicit solution \eqref{PhiPEFinal}. For $v \in \llbracket 1, N\rrbracket^K$ and $k \in \llbracket 1, N\rrbracket$, we denote
\[\varphi_{k, K}(v) = \sum_{s=1}^K \delta_{k v_s} = \# \{s \in \llbracket 1, K\rrbracket \suchthat v_s = k\},\]
and, for $\mathfrak n \in \mathfrak N$ with $\abs{\mathfrak n}_{\ell^1} \geq K$, we introduce the set
\[\Phi_K(\mathfrak n) = \{v \in \llbracket 1, N\rrbracket^K \suchthat n_j \geq \varphi_{j, K}(v) \text{ for all } j \in \llbracket 1, N\rrbracket\}.\]
Then we have the following result.

\begin{prop}
\label{PropAverage}
Let $K \in \mathbb N^\ast$ and suppose that $\mathfrak n \in \mathfrak N$ is such that $\abs{\mathfrak n}_{\ell^1} \geq K$. Then, for every $i, j \in \llbracket 1, N\rrbracket$ and $t \in\R$, we have
\begin{equation}
\vartheta^{(i)}_{j, \mathfrak n, L_j, t} = \sum_{v \in \Phi_K(\mathfrak n)} \left[\left(m_{v_1 j} \prod_{s=2}^K m_{v_s v_{s-1}}\right) \left(\prod_{s=1}^K \varepsilon_{v_s, \mathfrak n - \sum_{r=1}^s \mathbf 1_{v_r}, 0, t}\right) \vartheta^{(i)}_{v_K, \mathfrak n - \sum_{s=1}^K \mathbf 1_{v_s}, L_{v_K}, t}\right].
\label{VarthetaAverage}
\end{equation}
\end{prop}

\begin{prf}
The proof is done by induction on $K$. If $K = 1$, we have
\[\Phi_1(\mathfrak n) = \{v \in \llbracket 1, N\rrbracket \suchthat n_j \geq \delta_{j v} \text{ for all } j \in \llbracket 1, N\rrbracket\} = \{v \in \llbracket 1, N\rrbracket \suchthat n_v \geq 1\},\]
and so, by \eqref{ThetaIX-PE} and \eqref{VarthetaRecurrence},
\[\sum_{v \in \Phi_1(\mathfrak n)} \left[m_{v j}\varepsilon_{v, \mathfrak n - \mathbf 1_v, 0, t}\vartheta^{(i)}_{v, \mathfrak n - \mathbf 1_v, L_v, t} \right] = \sum_{\substack{v=1 \\ n_v \geq 1}}^N m_{v j} \varepsilon_{v, \mathfrak n - \mathbf 1_v, 0, t} \vartheta^{(i)}_{v, \mathfrak n - \mathbf 1_v, L_v, t} = \vartheta^{(i)}_{j, \mathfrak n, L_j, t}.\]

Suppose now that $K \in \mathbb N^\ast$ with $K \leq \abs{\mathfrak n}_{\ell^1}$ and \eqref{VarthetaAverage} holds true for $K-1$. Then we have, by \eqref{VarthetaRecurrence},
\[
\begin{aligned}
\vartheta^{(i)}_{j, \mathfrak n, L_j, t} & = \sum_{v^\prime \in \Phi_{K-1}(\mathfrak n)} \left[\left(m_{v^\prime_1 j} \prod_{s=2}^{K-1} m_{v^\prime_s v^\prime_{s-1}}\right)  \left(\prod_{s=1}^{K-1} \varepsilon_{v^\prime_s, \mathfrak n - \sum_{r=1}^s \mathbf 1_{v^\prime_r}, 0, t}\right)\vartheta^{(i)}_{v^\prime_{K-1}, \mathfrak n - \sum_{s=1}^{K-1} \mathbf 1_{v^\prime_s}, L_{v^\prime_{K-1}}, t}\right] \\
 & = \sum_{v^\prime \in \Phi_{K-1}(\mathfrak n)} \left[\left(m_{v^\prime_1 j} \prod_{s=2}^{K-1} m_{v^\prime_s v^\prime_{s-1}}\right)\left(\prod_{s=1}^{K-1} \varepsilon_{v^\prime_s, \mathfrak n - \sum_{r=1}^s \mathbf 1_{v^\prime_r}, 0, t}\right)\right. \\
 & \hphantom{=\sum_{v^\prime \in \Phi_{K-1}(\mathfrak n)}[} \left. \left(\sum_{\substack{k = 1 \\ n_k > \varphi_{k, K - 1}(v^\prime)}}^N m_{k v^\prime_{K-1}}  \varepsilon_{k, \mathfrak n - \sum_{s=1}^{K-1} \mathbf 1_{v^\prime_s} - \mathbf 1_k, 0, t}\vartheta^{(i)}_{k, \mathfrak n - \sum_{s=1}^{K-1} \mathbf 1_{v^\prime_s} - \mathbf 1_k, L_k, t} \right) \right] \\
 & = \sum_{v \in \Phi_K(\mathfrak n)} \left[\left(m_{v_1 j} \prod_{s=2}^{K} m_{v_s v_{s-1}}\right)  \left(\prod_{s=1}^{K} \varepsilon_{v_s, \mathfrak n - \sum_{r=1}^s \mathbf 1_{v_r}, 0, t}\right)\vartheta^{(i)}_{v_{K}, \mathfrak n - \sum_{s=1}^{K} \mathbf 1_{v_s}, L_{v_{K}}, t}\right],
\end{aligned}
\]
where we take $v = (v^\prime, v_N = k) = (v^\prime_1, \dotsc, v^\prime_{N-1}, k)$ and notice that $\varphi_{j, K}(v^\prime, k) = \varphi_{j, K-1}(v^\prime) + \delta_{jk}$ for every $j \in \llbracket 1, N\rrbracket$, so that
\begin{multline*}
\{(v^\prime, k) \in \Phi_{K-1}(\mathfrak n) \times \llbracket 1, N\rrbracket \suchthat n_k > \varphi_{k, K-1}(v^\prime)\} \\ = \{(v^\prime, k) \in \llbracket 1, N\rrbracket^{K-1} \times \llbracket 1, N\rrbracket \suchthat n_j \geq \varphi_{j, K-1}(v^\prime) \text { for all } j \in \llbracket 1, N\rrbracket \text{ and } n_k > \varphi_{k, K-1}(v^\prime)\} \\ = \{v \in \llbracket 1, N\rrbracket^K \suchthat n_j \geq \varphi_{j, K}(v) \text{ for all } j \in \llbracket 1, N\rrbracket\} = \Phi_K(\mathfrak n).
\end{multline*}
This proves \eqref{VarthetaAverage} by induction.
\end{prf}

\begin{remk}
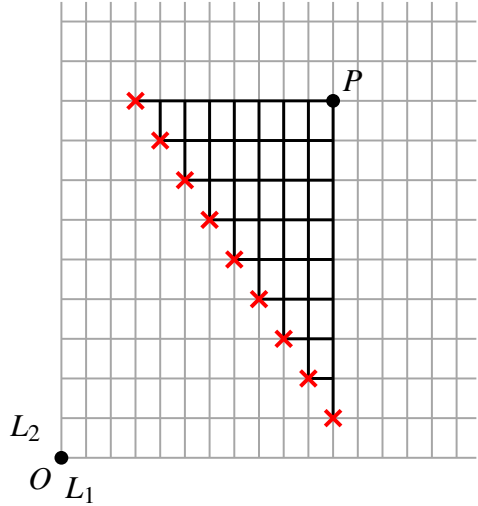
\begin{figure}[ht]
\centering

\setlength{\unitlength}{1.3cm}

\begin{picture}(5, 5.5)(-0.5, -0.5)

\linethickness{0.02\unitlength}%
\textcolor{gray}{\multiput(0, 0)(0   , 0.404508){12}{\Line(0, 0)(4.2, 0       )}}%
\textcolor{gray}{\multiput(0, 0)(0.25, 0       ){17}{\Line(0, 0)(0  , 4.649593)}}%

\linethickness{0.03\unitlength}
\put(2.75    , 3.640576){\Line(0, 0)(-2   , 0)}
\put(2.75    , 3.236068){\Line(0, 0)(-1.75, 0)}
\put(2.75    , 2.831559){\Line(0, 0)(-1.5 , 0)}
\put(2.75    , 2.427051){\Line(0, 0)(-1.25, 0)}
\put(2.75    , 2.022542){\Line(0, 0)(-1   , 0)}
\put(2.75    , 1.618034){\Line(0, 0)(-0.75, 0)}
\put(2.75    , 1.213525){\Line(0, 0)(-0.5 , 0)}
\put(2.75    , 0.809017){\Line(0, 0)(-0.25, 0)}
\put(2.75    , 3.640576){\Line(0, 0)(0, -3.236068)}
\put(2.5     , 3.640576){\Line(0, 0)(0, -2.831559)}
\put(2.25    , 3.640576){\Line(0, 0)(0, -2.427051)}
\put(2       , 3.640576){\Line(0, 0)(0, -2.022542)}
\put(1.75    , 3.640576){\Line(0, 0)(0, -1.618034)}
\put(1.5     , 3.640576){\Line(0, 0)(0, -1.213525)}
\put(1.25    , 3.640576){\Line(0, 0)(0, -0.809017)}
\put(1       , 3.640576){\Line(0, 0)(0, -0.404508)}
%
%
\linethickness{0.04\unitlength}
\multiput(0.75, 3.640576)(0.25, -0.40450850){9}{\textcolor{red}{\Line(-0.08, -0.08)(0.08, 0.08)\Line(-0.08, 0.08)(0.08, -0.08)}}%

\put(0       , 0       ){\circle*{0.14}}
\put(2.75    , 3.640576){\circle*{0.14}}

\put(-0.1 , -0.1 ){\makebox(0, 0)[tr]{$O$}}
\put( 2.85,  3.74){\makebox(0, 0)[bl]{$P$}}
\put( 0.2 , -0.2 ){\makebox(0, 0)[t]{$L_1$}}
\put(-0.2 ,  0.3 ){\makebox(0, 0)[r]{$L_2$}}
\end{picture}

\caption{Here $N=2$, $\mathfrak n=(11,9)$ and $K=8$. }
\label{FigAverage}
\end{figure}

In terms of the construction of Figure~\ref{FigTranspDamped}, each $v \in \Phi_K(\mathfrak n)$
corresponds to a path between $P$ and a point represented by a cross in Figure~\ref{FigAverage}. The term
\[\prod_{s=1}^K \varepsilon_{v_s, \mathfrak n - \sum_{r=1}^s \mathbf 1_{v_r}, 0, t}\]
in \eqref{VarthetaAverage} represents the total decay along this path.
\end{remk}

\begin{remk}
Under Hypothesis~\ref{HypoM}, the terms $m_{v_1 j} \prod_{s=2}^K m_{v_s v_{s-1}}$ in \eqref{VarthetaAverage} satisfy
\begin{equation}
\sum_{v \in \Phi_K(\mathfrak n)} \abs{m_{v_1 j} \prod_{s=2}^K m_{v_s v_{s-1}}} \leq 1.
\label{SumLeq1}
\end{equation}
Indeed, we have
\[\sum_{v \in \Phi_K(\mathfrak n)} \abs{m_{v_1 j} \prod_{s=2}^K m_{v_s v_{s-1}}} \leq \sum_{v \in \llbracket 1, N\rrbracket^K} \abs{m_{v_1 j} \prod_{s=2}^K m_{v_s v_{s-1}}} = \sum_{v_1=1}^N \dotsi \sum_{v_K=1}^N \abs{m_{v_K v_{K-1}}} \dotsm \abs{m_{v_1 j}} \leq 1,\]
by applying iteratively \eqref{CondNormL1}. We also remark that Proposition~\ref{PropAverage} implies that the coefficient $\vartheta^{(i)}_{j, \mathfrak n, L_j, t}$ belongs to the convex hull of the points $\pm \vartheta^{(i)}_{v_K, \mathfrak n - \sum_{s=1}^K \mathbf 1_{v_s}, L_{v_K}, t}$, for $v$ in $\Phi_K(\mathfrak n)$.
\end{remk}


\section{Proof of the main result}
\label{SecConvergence}

We now study the asymptotic behavior of the solutions of \eqref{damped0} with persistent damping, through their explicit formula obtained in Section~\ref{SecExplicit}. We first show that, in order to obtain the exponential stability of the solutions of \eqref{damped0}, it suffices to obtain the exponential convergence as $\abs{\mathfrak n}_{\ell^1} \to +\infty$ of the coefficients $\vartheta^{(i)}_{j, \mathfrak n, x, t}$ of the explicit formula \eqref{PhiPEFinal}. 


\subsection{Convergence of the coefficients implies convergence of the solution}

The section is devoted to the proof of the following result.

\begin{prop}
\label{LemmConvCoefConvSol}
Let $\mathcal F\subset L^\infty(\R,[0,1])$. Suppose that there exist constants $C_0, \gamma_0 > 0$ such that, for every $\alpha_k \in \mathcal F$, $k \in \llbracket 1, N_d\rrbracket$, we have
\[
\abs{\vartheta^{(i)}_{j, \mathfrak n, x, t}} \leq C_0 e^{-\gamma_0 \abs{\mathfrak n}_{\ell^1}}, \qquad \forall i, j \in \llbracket 1, N\rrbracket,\; \forall \mathfrak n \in \mathfrak N,\; \forall x \in [0, L_j],\; \forall t \in \R.
\]
Then there exist constants $C, \gamma > 0$ such that, for every $p \in \left[1, +\infty\right)$ and every initial condition $z_0 \in \mathsf X_p$, the corresponding solution $z(t)$ of \eqref{damped0} satisfies
\begin{equation}
\norm{z(t)}_{\mathsf X_p} \leq C e^{-\gamma t} \norm{z_0}_{\mathsf X_p}, \qquad \forall t \geq 0.
\label{SolConvExpo}
\end{equation}
\end{prop}

\begin{remk}
The conclusion \eqref{SolConvExpo} of Proposition~\ref{LemmConvCoefConvSol} can be written, in terms of the evolution family $\{T(t, s)\}_{t \geq s \geq 0}$ associated with \eqref{damped0}, as
\[\norm{T(t, 0)}_{\mathcal L(\mathsf X_p)} \leq C e^{-\gamma t}, \qquad t \geq 0.\]
When the class $\mathcal F$  is invariant by time-translation (e.g., for $\mathcal F=\mathcal G(T, \mu)$), this is actually equivalent to
\[\norm{T(t, s)}_{\mathcal L(\mathsf X_p)} \leq C e^{-\gamma (t - s)}, \qquad t \geq s \geq 0.\]
\end{remk}

\begin{prf}
Set $L_{\max} = \max_{i \in \llbracket 1, N\rrbracket}L_i$ and $L_{\min} = \min_{i \in \llbracket 1, N\rrbracket} L_i$. Take $z_0 = (u_{1, 0}, \dotsc, u_{N, 0}) \in \mathsf X_p$ and denote by $z(t) = (u_1(t), \dotsc, u_N(t))$ the corresponding solution of \eqref{damped0}. By Corollary~\ref{CoroBoundPE}, Theorem~\ref{TheoPEExplicit}, and Remark~\ref{densityDAinXp}, we have, for $t \geq L_{\max}$,
\begin{equation}
\norm{z(t)}_{\mathsf X_p}^p = \sum_{i=1}^N \norm{u_i(t)}_{L^p(0, L_i)}^p \leq \sum_{i=1}^N \norm{u_i(\cdot, 0)}_{L^p(t - L_i, t)}^p,
\label{EstimSolution}
\end{equation}
with $u_i(t, 0)$ given by \eqref{PhiPEFinal}. Denoting
\[Y_j(t) = \#\{\mathfrak n \in \mathfrak N_j \suchthat L(\mathfrak n) \leq t\},\]
 we have
\begin{multline}
\norm{u_i(\cdot, 0)}_{L^p(t - L_i, t)}^p = \int_{t - L_i}^t \abs{u_i(s, 0)}^p ds \displaybreak[0]\\%
 \leq N^{p-1} \sum_{j=1}^N \int_{t - L_i}^t \Bigabs{\sum_{\substack{\mathfrak n \in \mathfrak N_j \\ L(\mathfrak n) \leq s}} \vartheta^{(i)}_{j, \mathfrak n + \floor{\frac{s - L(\mathfrak n)}{L_j}} \mathbf 1_j, L_j - \{s - L(\mathfrak n)\}_{L_j},s} u_{j, 0}(L_j - \{s - L(\mathfrak n)\}_{L_j})}^p ds \displaybreak[0] \\
 \leq N^{p-1} \sum_{j=1}^N \int_{t - L_i}^t Y_j(s)^{p-1} \sum_{\substack{\mathfrak n \in \mathfrak N_j \\ L(\mathfrak n) \leq s}} \Bigabs{\vartheta^{(i)}_{j, \mathfrak n + \floor{\frac{s - L(\mathfrak n)}{L_j}} \mathbf 1_j, L_j - \{s - L(\mathfrak n)\}_{L_j},s} u_{j, 0}(L_j - \{s - L(\mathfrak n)\}_{L_j})}^p ds \displaybreak[0] \\
 \leq N^{p-1} C_0^p \sum_{j=1}^N Y_j(t)^{p-1} \sum_{\substack{\mathfrak n \in \mathfrak N_j \\ L(\mathfrak n) \leq t}} \int_{t - L_i}^t e^{-p \gamma_0 \left(\abs{\mathfrak n}_{\ell^1} + \floor{\frac{s - L(\mathfrak n)}{L_j}}\right)} \abs{u_{j, 0}(L_j - \{s - L(\mathfrak n)\}_{L_j})}^p ds \displaybreak[0] \\
 \leq N^{p-1} C_0^p e^{2p \gamma_0} e^{- \frac{p \gamma_0}{L_{\max}} t} \sum_{j=1}^N Y_j(t)^{p-1} \sum_{\substack{\mathfrak n \in \mathfrak N_j \\ L(\mathfrak n) \leq t}} \int_{t - L_{\max}}^t \abs{u_{j, 0}(L_j - \{s - L(\mathfrak n)\}_{L_j})}^p ds, \\
\label{EstimNormPhi}
\end{multline}
where we use that
\[\abs{\mathfrak n}_{\ell^1} + \floor{\frac{s - L(\mathfrak n)}{L_j}} = \sum_{k=1}^N n_k + \floor{\frac{s - L(\mathfrak n)}{L_j}} \geq \frac{L(\mathfrak n)}{L_{\max}} + \floor{\frac{s - L(\mathfrak n)}{L_{\max}}} \geq \frac{s}{L_{\max}} - 1 \geq \frac{t}{L_{\max}} - 2,\]
for $\mathfrak n \in \mathfrak N_j$ with $L(\mathfrak n) \leq t$ and $s \in [t - L_i, t]$.

According to its definition, $Y_j(t)$ can be upper bounded as follows
\begin{align}
Y_j(t) & \leq \#\{\mathfrak n \in \mathfrak N_j \suchthat n_i L_i \leq t \text{ for all } i \in \llbracket 1, N \rrbracket \setminus\{j\}\}\nonumber \\
 & = \# \left(\left\llbracket 0, \frac{t}{L_1} \right\rrbracket \times \dotsb \times \left\llbracket 0, \frac{t}{L_{j-1}} \right\rrbracket \times \{0\} \times \left\llbracket 0, \frac{t}{L_{j+1}} \right\rrbracket \times \dotsb \times \left\llbracket 0, \frac{t}{L_N} \right\rrbracket \right) \nonumber\\
 & \leq \left(\frac{t}{L_{\min}} + 1\right)^{N-1}.\label{EstimYjt}
\end{align}

We next estimate $\int_{t - L_{\max}}^t \abs{u_{j, 0}(L_j - \{s - L(\mathfrak n)\}_{L_j})}^p ds$ with $j \in \llbracket 1, N\rrbracket$. Notice that $\left[t-L_{\max}, t\right)\allowbreak \subset \cup_{k=k_{\min}}^{k_{\max}} \left[L(\mathfrak n) + k L_j, L(\mathfrak n) + (k+1) L_j\right)$ with 
\[
k_{\min} = \max\{k \in \mathbb Z \suchthat L(\mathfrak n) + k L_j \leq t - L_{\max}\},\quad 
k_{\max} = \min\{k \in \mathbb Z \suchthat L(\mathfrak n) + (k+1) L_j \geq t\}.
\]
We deduce that
\begin{align}
\int_{t - L_{\max}}^t \abs{u_{j, 0}(L_j - \{s - L(\mathfrak n)\}_{L_j})}^p ds \displaybreak[0]&\leq \sum_{k=k_{\min}}^{k_{\max}} \int_0^{L_j} \abs{u_{j, 0}(\sigma)}^p d\sigma\nonumber\\
 & = (k_{\max} - k_{\min} + 1) \norm{u_{j, 0}}_{L^p(0, L_j)}^p\nonumber\\
 &\leq \left(\frac{L_{\max}}{L_j} + 2\right) \norm{u_{j, 0}}_{L^p(0, L_j)}^p.\label{EstimNormPhi2}
\end{align}

Inserting \eqref{EstimYjt} and \eqref{EstimNormPhi2} into \eqref{EstimNormPhi} finally gives \eqref{SolConvExpo} thanks to \eqref{EstimSolution}. Notice that the coefficients $\gamma$ and $C$ can be chosen  to be independent of $p$.
\end{prf}

\begin{remk}
\label{RemkInfty}
Even though the well-posedness of \eqref{TranspOperatorPE} was considered in Section~\ref{SecWellPosed} only for $p \in \left[1, +\infty\right)$, we extend it here below to the case $p = +\infty$ and we verify that Proposition~\ref{LemmConvCoefConvSol} still holds true in this case.

First  set $\mathsf X_{\infty} = \prod_{i=1}^N L^\infty(0, L_i)$ with its usual norm $\norm{z}_{\mathsf X_\infty} = \max_{i \in \llbracket 1, N\rrbracket} \norm{u_i}_{L^\infty(0, L_i)}$ for $z = (u_1, \dotsc,\allowbreak u_N) \in \mathsf X_\infty$. Fix $z_0\in \mathsf X_\infty$. Since $z_0\in \mathsf X_p$ for every $p \in \left[1, +\infty\right)$, then \eqref{TranspOperatorPE} admits a unique mild solution $z(t) = T(t, 0) z_0$ in $\bigcap_{p \in \left[1, +\infty\right)} \mathcal C^0(\mathbb R_+, \mathsf X_p)$ with initial condition $z(0) = z_0$. As noticed in Remark~\ref{densityDAinXp}, $z(t)$ is  characterized as an element of $\mathsf X_p$ by equations \eqref{FlowU1UiPE} and  \eqref{PhiPEFinal}. Hence $z(t) \in \mathsf X_\infty$ for every $t \geq 0$.
We can thus refer to $z(\cdot)$ as the solution of the Cauchy problem \eqref{TranspOperatorPE} in $\mathsf X_\infty$.

Suppose now that the hypotheses of Proposition~\ref{LemmConvCoefConvSol} are satisfied and let $C, \gamma > 0$ be as in its statement. By \eqref{SolConvExpo}, we have
\[\norm{z(t)}_{\mathsf X_p} \leq C e^{-\gamma t} \norm{z_0}_{\mathsf X_p} \leq C \left(\sum_{i=1}^N L_i\right)^{\nicefrac{1}{p}} e^{-\gamma t} \norm{z_0}_{\mathsf X_\infty}.\]
Since $\norm{y}_{\mathsf X_\infty} = \lim_{p \to +\infty} \norm{y}_{\mathsf X_p}$ for every $y \in \mathsf X_\infty$, we conclude that
\[\norm{z(t)}_{\mathsf X_\infty} \leq C e^{-\gamma t} \norm{z_0}_{\mathsf X_\infty}.\]
\end{remk}


\subsection{Preliminary estimates of $\vartheta^{(i)}_{j, \mathfrak n, x, t}$}

In this section we establish estimates on the growth of $\vartheta^{(i)}_{j, \mathfrak n, x, t}$ based on combinatorial arguments.

\begin{prop}
\label{PropRough}
For every $i, j \in \llbracket 1, N\rrbracket$, $\mathfrak n \in \mathfrak N$, $x \in [0, L_j]$, $t \in \R$, and $\alpha_k \in L^\infty(\mathbb R, [0, 1])$, $k \in \llbracket 1, N_d\rrbracket$, we have
\begin{equation}
\abs{\vartheta^{(i)}_{j, \mathfrak n, x, t}} \leq \abs{M}_{\ell^1}^{\abs{\mathfrak n}_{\ell^1} + 1}.
\label{EstimRough}
\end{equation}
\end{prop}

\begin{prf}
We show \eqref{EstimRough} by induction on $\abs{\mathfrak n}_{\ell^1}$. For every $i, j \in \llbracket 1, N\rrbracket$, $x \in [0, L_j]$, $t \in\R$, and $\alpha_k \in L^\infty(\mathbb R, [0, 1])$, $k \in \llbracket 1, N_d\rrbracket$, we have, by \eqref{Vartheta},
\[\abs{\vartheta^{(i)}_{j, \mathbf 0, x, t}} \leq \abs{\vartheta^{(i)}_{j, \mathbf 0, L_j, t}} = \abs{m_{ij}} \leq \abs{M}_{\ell^1}.\]
If $R \in \mathbb N$ is such that \eqref{EstimRough} holds for every $i, j \in \llbracket 1, N\rrbracket$, $\mathfrak n \in \mathfrak N$ with $\abs{\mathfrak n}_{\ell^1} = R$, $x \in [0, L_j]$, $t \in \R$, and $\alpha_k \in L^\infty(\mathbb R, [0, 1])$, $k \in \llbracket 1, N_d\rrbracket$, then, for $\mathfrak n \in \mathfrak N$ with $\abs{\mathfrak n}_{\ell^1} = R+1$, we have, by \eqref{Vartheta},
\[
\abs{\vartheta^{(i)}_{j, \mathfrak n, x, t}} \leq \abs{\vartheta^{(i)}_{j, \mathfrak n, L_j, t}} \leq \sum_{\substack{r=1 \\ n_r \geq 1}}^N \abs{m_{rj}} \abs{\vartheta^{(i)}_{r, \mathfrak n - \mathbf 1_r, L_r, t}} \leq \sum_{r=1}^N \abs{m_{rj}} \abs{M}_{\ell^1}^{R + 1} \leq \abs{M}_{\ell^1}^{R + 2},
\]
since $\abs{M}_{\ell^1} = \max_{j \in \llbracket 1, N\rrbracket} \sum_{r=1}^N \abs{m_{rj}}$. The result thus follows by induction.
\end{prf}

As a  consequence of Propositions~\ref{LemmConvCoefConvSol} and \ref{PropRough} and Remark~\ref{RemkInfty}, we deduce at once the following corollary.  

\begin{coro}
\label{CoroUndamped}
Suppose that $\abs{M}_{\ell^1} < 1$. Then there exist $C, \gamma > 0$ such that, for every $p \in [1, +\infty]$ and every initial condition $z_0 \in \mathsf X_p$, the corresponding solution $z(t)$ of the undamped equation \eqref{Undamped} satisfies
\begin{equation*}
\norm{z(t)}_{\mathsf X_p} \leq C e^{-\gamma t} \norm{z_0}_{\mathsf X_p}, \qquad \forall t \geq 0.
\end{equation*}
\end{coro}

Another trivial but important consequence of Proposition~\ref{PropRough} is that, if $\abs{M}_{\ell^1} \leq 1$, the coefficients $\vartheta^{(i)}_{j, \mathfrak n, x, t}$ are all bounded in absolute value by $1$.

\begin{coro}
\label{CoroLeq1}
Suppose that $\abs{M}_{\ell^1} \leq 1$. Then, for every $i, j \in \llbracket 1, N\rrbracket$, $\mathfrak n \in \mathfrak N$, $x \in [0, L_j]$, $t \in \R$,
and $\alpha_k \in L^\infty(\mathbb R, [0, 1])$, $k \in \llbracket 1, N_d\rrbracket$, we have
\begin{equation}
\abs{\vartheta^{(i)}_{j, \mathfrak n, x, t}} \leq 1.
\label{VarthetaLeq1}
\end{equation}
\end{coro}

Our second estimate on the coefficients $\vartheta^{(i)}_{j, \mathfrak n, x, t}$ is the following.

\begin{lemm}
\label{CoroBinBound}
Suppose that Hypothesis~\ref{HypoM} is satisfied. Then there exists $\nu \in (0, 1)$ such that, for every $i, j, k \in \llbracket 1, N\rrbracket$, $\mathfrak n \in \mathfrak N$, $x \in [0, L_j]$, $t \in\R$, and $\alpha_r \in L^\infty(\mathbb R, [0, 1])$, $r \in \llbracket 1, N_d\rrbracket$, we have
\begin{equation}
\abs{\vartheta^{(i)}_{j, \mathfrak n, x, t}} \leq \binom{\abs{\mathfrak n}_{\ell^1}}{n_k} \nu^{\abs{\mathfrak n}_{\ell^1}}.
\label{BinomialBound}
\end{equation}
\end{lemm}

\begin{prf}
Up to a permutation in the set of indices, we can suppose, without loss of generality, that $k = N$. Let
\[
\mu_N = \max_{j \in \llbracket 1, N\rrbracket} \sum_{i=1}^{N-1} \abs{m_{ij}}, \qquad \nu_N = \max_{j \in \llbracket 1, N\rrbracket} \abs{m_{Nj}}.
\]
By Hypothesis~\ref{HypoM}, we have both $\sum_{i=1}^{N-1} \abs{m_{ij}} < 1$ and $\abs{m_{Nj}} < 1$. Hence $\mu_N, \nu_N \in (0, 1)$.

We prove by induction on $\abs{\mathfrak n}_{\ell^1}$ that
\begin{equation}
\abs{\vartheta^{(i)}_{j, \mathfrak n, x, t}} \leq \binom{\abs{\mathfrak n}_{\ell^1}}{n_k} \mu_N^{\abs{\mathfrak n}_{\ell^1}-n_N}\nu_{N}^{n_N}.
\label{BinomialKBound}
\end{equation}
For every $i, j \in \llbracket 1, N\rrbracket$, $x \in [0, L_j]$, $t \in\R$, and $\alpha_r \in L^\infty(\mathbb R, [0, 1])$, $r \in \llbracket 1, N_d\rrbracket$, we have, by \eqref{Vartheta},
\[\abs{\vartheta^{(i)}_{j, \mathbf 0, x, t}} \leq \abs{\vartheta^{(i)}_{j, \mathbf 0, L_j, t}} = \abs{m_{ij}} \leq 1,\]
so that \eqref{BinomialKBound} is satisfied for $\mathfrak n = \mathbf 0$.

Suppose now that $R \in \mathbb N$ is such that \eqref{BinomialKBound} is satisfied for every $i, j \in \llbracket 1, N\rrbracket$, $\mathfrak n \in \mathfrak N$ with $\abs{\mathfrak n}_{\ell^1} = R$, $x \in [0, L_j]$, $t \in \mathbb R$, 
and $\alpha_r \in L^\infty(\mathbb R, [0, 1])$, $r \in \llbracket 1, N_d\rrbracket$. If $\mathfrak n \in \mathfrak N$ is such that $\abs{\mathfrak n}_{\ell^1} = R + 1$, we have, by \eqref{Vartheta},
\begin{align*}
\abs{\vartheta^{(i)}_{j, \mathfrak n, x, t}} & \leq \abs{\vartheta^{(i)}_{j, \mathfrak n, L_j, t}} \leq \sum_{\substack{r=1 \\ n_r \geq 1}}^{N-1} \abs{m_{rj}} \abs{\vartheta^{(i)}_{r, \mathfrak n - \mathbf 1_r, L_r, t}} + \abs{m_{Nj}} \abs{\vartheta^{(i)}_{N, \mathfrak n - \mathbf 1_N, L_N, t}} \displaybreak[0] \\
 & \leq \sum_{\substack{r=1 \\ n_r \geq 1}}^{N-1} \abs{m_{rj}} \binom{\abs{\mathfrak n}_{\ell^1} - 1}{n_N} \mu_N^{\abs{\mathfrak n}_{\ell^1} - n_N - 1} \nu_N^{n_N} + \abs{m_{Nj}} \binom{\abs{\mathfrak n}_{\ell^1} - 1}{n_N - 1} \mu_N^{\abs{\mathfrak n}_{\ell^1} - n_N} \nu_N^{n_N - 1} \displaybreak[0] \\
 & \leq \binom{\abs{\mathfrak n}_{\ell^1} - 1}{n_N} \mu_N^{\abs{\mathfrak n}_{\ell^1} - n_N} \nu_N^{n_N} + \binom{\abs{\mathfrak n}_{\ell^1} - 1}{n_N - 1} \mu_N^{\abs{\mathfrak n}_{\ell^1} - n_N} \nu_N^{n_N} \displaybreak[0] \\
 & = \binom{\abs{\mathfrak n}_{\ell^1}}{n_N} \mu_N^{\abs{\mathfrak n}_{\ell^1} - n_N} \nu_N^{n_N},
\end{align*}
with the convention that $\vartheta^{(i)}_{N, \mathfrak n - \mathbf 1_N, L_N, t} = 0$ if $n_N = 0$. Hence \eqref{BinomialKBound} holds for $\mathfrak n$, which proves the result by induction. We conclude by taking 
$\nu = \max \{\mu_k, \nu_k \suchthat k \in \llbracket 1, N\rrbracket\}$.
\end{prf}


\subsection{Exponential decay of $\vartheta^{(i)}_{j, \mathfrak n, x, t}$ in $\mathfrak N_b(\rho)$}

The proof of  the exponential decay of the coefficients $\vartheta^{(i)}_{j, \mathfrak n, x, t}$ as $\abs{\mathfrak n}_{\ell^1} \to +\infty$, uniformly with respect to $\alpha_k \in \mathcal G(T, \mu)$, $k \in \llbracket 1, N_d\rrbracket$, 
is split into two cases. We first estimate $\vartheta^{(i)}_{j, \mathfrak n, x, t}$ for  $\mathfrak n$ in a subset $\mathfrak N_b(\rho)$ of $\mathfrak N$, namely when one of the components of $\mathfrak n$ is much smaller than the others. The parameter $\rho\in (0,1)$ is a measure of such a smallness and will be fixed later. For $\mathfrak n\in \mathfrak N_b(\rho)$, the exponential decay of $\vartheta^{(i)}_{j, \mathfrak n, x, t}$ does not result from the presence of the persistent damping but solely from combinatorial considerations. We then proceed in Section~\ref{SecNc} to estimate $\vartheta^{(i)}_{j, \mathfrak n, x, t}$ in the set $\mathfrak N_c(\rho) = \mathfrak N \setminus \mathfrak N_b(\rho)$, where the decay comes from the persistent damping in \eqref{damped0}.

\begin{defi}\label{def_Nb}
For $k \in \llbracket 1, N\rrbracket$ and $\rho \in (0, 1)$, we define
\begin{gather*}
\mathfrak N_b(\rho, k) = \{\mathfrak n = (n_1, \dotsc, n_N) \in \mathfrak N \suchthat n_k \leq \rho \abs{\mathfrak n}_{\ell^1}\}, \\
\mathfrak N_b(\rho) = \bigcup_{k=1}^N \mathfrak N_b(\rho, k), \qquad \mathfrak N_c(\rho) = \mathfrak N \setminus \mathfrak N_b(\rho).
\end{gather*}
\end{defi}

We now deduce from Lemma~\ref{CoroBinBound} the exponential decay of $\vartheta^{(i)}_{j, \mathfrak n, x, t}$ in the set $\mathfrak N_b(\rho)$.

\begin{theo}
\label{TheoVarthetaNb}
Suppose that Hypothesis~\ref{HypoM} is satisfied. There exist $\rho \in \left(0, \nicefrac{1}{2}\right)$ and constants $C, \gamma > 0$ such that, for every $i, j \in \llbracket 1, N\rrbracket$, $\mathfrak n \in \mathfrak N_b(\rho)$, $x \in [0, L_j]$, $t \in \R$, and $\alpha_r \in L^\infty(\mathbb R, [0, 1])$, $r \in \llbracket 1, N_d\rrbracket$, we have
\begin{equation}
\abs{\vartheta^{(i)}_{j, \mathfrak n, x, t}} \leq C e^{-\gamma \abs{\mathfrak n}_{\ell^1}}.
\label{EstimVarthetaNb}
\end{equation}
\end{theo}

\begin{prf}
Let $\nu \in (0, 1)$ be as in Lemma~\ref{CoroBinBound}. 
According to Lemma~\ref{LemmCoeffBinom} in appendix, there exist $\rho \in \left(0, \nicefrac{1}{2}\right)$, $C, \gamma > 0$ such that for every $n \in \mathbb N$ and $k \in \llbracket 0, \rho n \rrbracket$, we have $\binom{n}{k} \nu^{n} \leq C e^{-\gamma n}$. Take $i, j \in \llbracket 1, N\rrbracket$, $\mathfrak n \in \mathfrak N_b(\rho)$, $x \in [0, L_j]$, and $\alpha_r \in L^\infty(\mathbb R, [0, 1])$ for $r \in \llbracket 1, N_d\rrbracket$. Since $\mathfrak n \in \mathfrak N_b(\rho)$, there exists $k \in \llbracket 1, N\rrbracket$ such that $\mathfrak n \in \mathfrak N_b(\rho, k)$, i.e., $n_k \leq \rho \abs{\mathfrak n}_{\ell^1}$. Then, by Lemmas~\ref{CoroBinBound} and \ref{LemmCoeffBinom}, we have
\[
\abs{\vartheta^{(i)}_{j, \mathfrak n, x, t}} \leq \binom{\abs{\mathfrak n}_{\ell^1}}{n_k} \nu^{\abs{\mathfrak n}_{\ell^1}} \leq C e^{-\gamma \abs{\mathfrak n}_{\ell^1}}.
\]
\end{prf}

\begin{remk}\label{TheoConvNonPE}
The above estimate is actually sufficient to derive the conclusion of Theorem~\ref{MainTheoIntro} when the damping is always active; indeed, in this case, one can easily deduce by an inductive argument using \eqref{Vartheta} that
\[\vartheta^{(i)}_{j, \mathfrak n, L_j, t} = \beta^{(i)}_{j, \mathfrak n} e^{-n_1 (b_1 - a_1)} e^{-n_2 (b_2 - a_2)} \dotsm e^{- n_{N_d} (b_{N_d} - a_{N_d})}\]
and the exponential decay in $\mathfrak N_c(\rho)$ follows straightforwardly. Notice that in this case Hypothesis~\ref{HypoIrrational} is not necessary. 
\end{remk}


\subsection{Exponential decay of $\vartheta^{(i)}_{j, \mathfrak n, x, t}$ in $\mathfrak N_c(\rho)$}
\label{SecNc}

In this section, we establish the exponential decay of $\vartheta^{(i)}_{j, \mathfrak n, x, t}$ in the set $\mathfrak N_c(\rho)$. The main difficulty in proving it lies in the fact that $\alpha_i(t)$ can be equal to zero for certain time intervals, so that the term $\varepsilon_{j, \mathfrak n, x, t}$ defined by \eqref{EpsilonIX-PE} can be equal to $1$. Recall that
\begin{equation}
\varepsilon_{j, \mathfrak n, 0, t} = e^{-\int_{t - L(\mathfrak n) - L_j + a_j}^{t - L(\mathfrak n) - L_j + b_j} \alpha_j(s) ds}.
\label{Varepsilon}
\end{equation}
Our goal consists in showing in Lemma~\ref{LemmPEEta} that $\varepsilon_{j, \mathfrak n, 0, t}$ is smaller than a certain value ``often enough''. The first step in this direction is the following lemma.

\begin{lemm}
\label{LemmIPE}
Let $T \geq \mu > 0$ and $j \in \llbracket 1, N_d\rrbracket$. For $\rho > 0$ and $\alpha \in \mathcal G(T, \mu)$, define
\begin{equation}
\mathcal I_{j, \rho, \alpha} = \left\{\tau \in \mathbb R \middlesuchthat \int_{\tau + a_j}^{\tau + b_j} \alpha(s) ds \geq \rho\right\}.
\label{DefIRhoAlpha}
\end{equation}
There exist $\rho_j > 0$ and $\ell_j > 0$, depending only on $\mu$, $T$ and $b_j - a_j$, such that, for every $t \in \mathbb R$ and $\alpha \in \mathcal G(T, \mu)$, $\mathcal I_{j, \rho_j, \alpha} \cap [t, t + T]$ contains an interval of length $\ell_j$.
\end{lemm}

\begin{prf}
We set $\rho_j = \frac{\mu (b_j - a_j)}{2T}$, $\ell_j = \min\{\rho_j, T\}$. Take $\alpha \in \mathcal G(T, \mu)$ and define the function $A: \mathbb R \to \mathbb R$ by
\[A(\tau) = \int_{\tau + a_j}^{\tau + b_j} \alpha(s) ds.\]
Since $\alpha \in L^\infty(\mathbb R, [0, 1])$, $A$ is $1$-Lipschitz continuous. We also have, for every $t \in \mathbb R$,
\begin{equation}
\int_{t}^{t + T} A(\tau) d\tau = \int_{t}^{t + T} \int_{a_j}^{b_j} \alpha(s + \tau) ds d\tau = \int_{a_j}^{b_j} \int_{s + t}^{s + t + T} \alpha(\tau) d\tau ds \geq \mu (b_j - a_j).
\label{PE-A}
\end{equation}
Take $t \in \mathbb R$. There exists $t_\star \in [t, t + T]$ such that $A(t_\star) \geq \frac{\mu (b_j - a_j)}{T} = 2 \rho_j$, for otherwise \eqref{PE-A} would not be satisfied. Since $A$ is $1$-Lipschitz continous, we have $A(\tau) \geq \rho_j$ for $\tau \in [t_\star - \rho_j, t_\star + \rho_j] $, and thus
\[\left[t_\star - \rho_j, t_\star + \rho_j\right]  \cap [t, t + T] \subset \mathcal I_{j, \rho_j, \alpha} \cap [t, t + T].\]
But, since $t_\star \in [t, t + T]$, $[t_\star - \rho_j, t_\star + \rho_j]  \cap [t, t + T]$ is an interval of length at least $\ell_j$, which concludes the proof.
\end{prf}

Lemma~\ref{LemmIPE} translates the persistence of excitation of $\alpha$ into a property on the integrals appearing in \eqref{Varepsilon}.

As remarked in Section~\ref{SecHypo}, one cannot expect to obtain a general result concerning the exponential stability of \eqref{damped0} without taking into account the rationality of the ratios ${L_i}/{L_j}$. The following lemma uses the irrationality of ${L_i}/{L_j}$ for certain $i, j \in \llbracket 1, N\rrbracket$ to give a further step into the understanding of $\varepsilon_{j, \mathfrak n, 0, t}$.

\begin{lemm}
\label{LemmTimeIRhoAlpha}
Let $T \geq \mu > 0$ and let $\rho_j > 0$, $j \in \llbracket 1, N_d\rrbracket$, be as in Lemma~\ref{LemmIPE}. There exists $K \in \mathbb N$ such that, for every $k_1 \in \llbracket 1, N\rrbracket$, $k_2 \in \llbracket 1, N_d\rrbracket $ with ${L_{k_1}}/{L_{k_2}} \notin \mathbb Q$, $t \in \mathbb R$, $\mathfrak n \in \mathfrak N$, and $\alpha \in \mathcal G(T, \mu)$, there exists $\mathfrak r \in \mathfrak N$ with $n_j \leq r_j \leq K + n_j$, $j \in \{k_1, k_2\}$, and $r_j = n_j$ for $j \in \llbracket 1, N\rrbracket \setminus \{k_1, k_2\}$, such that
\[t - L(\mathfrak r) \in \mathcal I_{k_2, \rho_{k_2}, \alpha}.\]
\end{lemm}

\begin{prf}
We shall prove the following simpler statement: for every $k_1 \in \llbracket 1, N\rrbracket$ and $k_2 \in \llbracket 1, N_d\rrbracket $ with ${L_{k_1}}/{L_{k_2}} \notin \mathbb Q$, there exist $N_1 = N_1(k_1, k_2) \in \mathbb N$ and $N_2 = N_2(k_1, k_2) \in \mathbb N$ such that, for every $t \in \mathbb R$, $\mathfrak n \in \mathfrak N$, and $\alpha \in \mathcal G(T, \mu)$, there exist $\mathfrak r \in \mathfrak N$ with $n_{k_j} \leq r_{k_j} \leq N_j + n_{k_j}$, $j \in \{1, 2\}$, and $r_j = n_j$ for $j \in \llbracket 1, N\rrbracket \setminus \{k_1, k_2\}$, such that
\[t - L(\mathfrak r) \in \mathcal I_{k_2, \rho_{k_2}, \alpha}.\]
From this result, one can easily obtain the statement of the lemma by taking
\[K = \max\left\{N_1(k_1, k_2),\; N_2(k_1, k_2) \middlesuchthat k_1 \in \llbracket 1, N\rrbracket,\; k_2 \in \llbracket 1, N_d\rrbracket  \text{ such that } \frac{L_{k_1}}{L_{k_2}} \notin \mathbb Q\right\}.\]

We decompose the argument into two steps.

\setcounter{step}{0}
\begin{step}
\label{Step1}
Definition of the points $x_j$ and $y_j$.
\end{step}

Let $\rho_{k_2} > 0$ and $\ell_{k_2} > 0$ be obtained from $\mu$, $T$ and $b_{k_2} - a_{k_2}$ as in Lemma~\ref{LemmIPE}. Let $\kappa = 3\ceil{{T}/{\ell_{k_2}}}$ and set
\[x_j = \frac{j}{\kappa}T, \qquad j \in \llbracket 0, \kappa \rrbracket,\]
which satisfy $x_j - x_{j-1} = \frac{\kappa}{T} \leq \frac{\ell_{k_2}}{3}$ for $j \in \llbracket 1, \kappa \rrbracket$. Hence, for every interval $J$ of length $\ell_{k_2}$ contained in $[0, T]$, there exists $j \in \llbracket 1, \kappa\rrbracket$ such that $x_{j-1}, x_j \in J$.

We now construct intermediate points between the $x_j$, $j \in \llbracket 0, \kappa \rrbracket$. Since ${L_{k_1}}/{L_{k_2}} \notin \mathbb Q$, the set
\begin{equation}
\{n_1 L_{k_1} + n_2 L_{k_2} \suchthat n_1, n_2 \in \mathbb Z\}
\label{Setn1L1n2L2}
\end{equation}
is dense in $\mathbb R$. Hence we can find $n_{1, j}, n_{2, j} \in \mathbb Z$, $j \in \llbracket 1, \kappa \rrbracket$, such that the numbers $y_j = n_{1, j} L_{k_1} +n_{2, j} L_{k_2}$ satisfy
\begin{equation}
0 = x_0 < y_1 < x_1 < y_2 < x_2 < \dotsb < y_\kappa < x_\kappa = T.
\label{YBetweenX}
\end{equation}
As a consequence, for any interval $J$ of length $\ell_{k_2}$ contained in $[0, T]$, there exists $j \in \llbracket 1, \kappa\rrbracket$ such that $y_j \in J$.

\begin{step}
Characterization of $\mathfrak r \in \mathfrak N$ and conclusion.
\end{step}

Let $N_1^\star = \max\{|n_{1, 1}|, \dotsc, |n_{1, \kappa}|\}$, $N_2^\star = \max\{ |n_{2, 1}|,\dotsc, |n_{2, \kappa}|\}$ and $N_1=2N_1^\star$, $N_2=2N_2^\star$. Take $t \in \mathbb R$ and $\mathfrak n \in \mathfrak N$. For $j \in \llbracket 1, \kappa \rrbracket$, define $\mathfrak r_j = (r_{1, j}, \dotsc, r_{N, j}) \in \mathfrak N$ by
\[r_{k_1, j} = n_{k_1} + n_{1, j} + N_1^\star, \qquad r_{k_2, j} = n_{k_2} + n_{2, j} + N_2^\star,\]
and $r_{i, j} = n_i$ for $i \in \llbracket 1, N\rrbracket \setminus\{k_1, k_2\}$; it is clear, by this definition, that $n_{k_i} \leq r_{k_i, j} \leq N_i + n_{k_i}$ for $i \in \{1, 2\}$ and $j \in \llbracket 1, \kappa \rrbracket$. Set
\[z_j = t - L(\mathfrak r_j) , \qquad j \in \llbracket 1, \kappa \rrbracket;\]
we thus have
\[z_j = t - n_{1, j} L_{k_1} - n_{2, j}  L_{k_2} - Z^\star = t - Z^\star - y_j\]
with $Z^\star = L(\mathfrak n) + N_1^\star L_{k_1} + N_2^\star L_{k_2}$. Since, by construction, $y_j \in (0, T)$ for $j \in \llbracket 1, \kappa \rrbracket$, we have $z_j \in [t - Z^\star - T, t - Z^\star]$. 

Take $\alpha \in \mathcal G(T, \mu)$. By Lemma~\ref{LemmIPE}, $\mathcal I_{k_2, \rho_{k_2}, \alpha} \cap [t - Z^\star - T, t - Z^\star]$ contains an interval $J$ of length $\ell_{k_2}$. Consider the interval $J^\prime = - J + t - Z^\star$, which is a subinterval of $[0, T]$ of length $\ell_{k_2}$. By Step~\ref{Step1}, there exists $j \in \llbracket 1, \kappa\rrbracket$ such that $y_j \in J^\prime$, and thus $z_j \in J \subset \mathcal I_{k_2, \rho_{k_2}, \alpha}$. Since $z_j = t - L(\mathfrak r_j)$, we obtain the desired result with $\mathfrak r = \mathfrak r_j$.
\end{prf}

\begin{remk}
\label{RemkRational}
The only instance in the proof of Lemma~\ref{LemmTimeIRhoAlpha} where we use the fact that ${L_{k_1}}/{L_{k_2}} \notin \mathbb Q$ is when we establish the existence of numbers $y_j$, $j = 1, \dotsc, \kappa$, of the form $y_j = n_{1, j} L_{k_1} + n_{2, j} L_{k_2}$ with $n_{1, j}, n_{2, j} \in \mathbb Z$ satisfying \eqref{YBetweenX}, which we do by using the density of the set \eqref{Setn1L1n2L2}. When ${L_{k_1}}/{L_{k_2}} \in \mathbb Q$ and we write ${L_{k_1}}/{L_{k_2}} = {p}/{q}$ for coprime $p, q \in \mathbb N^\ast$, the set given in \eqref{Setn1L1n2L2} is, by B\'ezout's Lemma,
\[\left\{L_{k_2} \frac{n_1 p + n_2 q}{q} \middlesuchthat n_1, n_2 \in \mathbb Z\right\} = \left\{k \frac{L_{k_2}}{q} \middlesuchthat k \in \mathbb Z\right\}.\]
Hence the construction of $y_j = n_{1, j} L_{k_1} + n_{2, j} L_{k_2}$ with $n_{1, j}, n_{2, j} \in \mathbb Z$ satisfying \eqref{YBetweenX} is still possible if ${L_{k_2}}/{q} < {\kappa}/{T}$, i.e., if $q > {L_{k_2} \kappa}/{T}$, and thus Lemma~\ref{LemmTimeIRhoAlpha} still holds true if ${L_{k_1}}/{L_{k_2}} = {p}/{q}$ with coprime $p, q \in \mathbb N$ and $q$ large enough.

Recalling that $\kappa = 3 \ceil{{T}/{\ell_{k_2}}}$ with $\ell_{k_2} = \min\left\{\frac{\mu (b_{k_2} - a_{k_2})}{2 T}, T\right\}$, we can even give a more explicit sufficient condition on $q$ to still have Lemma~\ref{LemmTimeIRhoAlpha}: if
\begin{equation*}
q \geq 3 L_{k_2} \left(\max\left\{\frac{2 T}{\mu (b_{k_2} - a_{k_2})}, \frac{1}{T}\right\} + \frac{1}{T}\right),
\end{equation*}
then one can easily check that $q > {L_{k_2} \kappa}/{T}$ and hence we are in the previous situation.

More explicitly, we can replace Hypothesis~\ref{HypoIrrational} by the following one.

\begin{hypo}
\label{HypoRational}
There exist $i \in \llbracket 1, N\rrbracket$ and $j \in \llbracket 1, N_d\rrbracket$ for which we have either ${L_i}/{L_j} \notin \mathbb Q$ or ${L_i}/{L_j} = {p}/{q}$ with coprime $p, q \in \mathbb N^\ast$ satisfying
\begin{equation}
q \geq 3 L_{j} \left(\max\left\{\frac{2 T}{\mu (b_{j} - a_{j})}, \frac{1}{T}\right\} + \frac{1}{T}\right).
\label{CondDenominator}
\end{equation}
\end{hypo}

Notice that condition \eqref{CondDenominator} only depends on the constants $T$, $\mu$ of the PE condition, on the length $b_j - a_j$ of the damping interval $j$ and on the length $L_j$.
\end{remk}

As a consequence of the previous lemma we deduce the following property.

\begin{lemm}
\label{LemmPEEta}
Let $T \geq \mu > 0$. There exist $\eta \in (0, 1)$ and $K \in \mathbb N$ such that, for every $k_1 \in \llbracket 1, N\rrbracket$, $k_2 \in \llbracket 1, N_d\rrbracket $ with ${L_{k_1}}/{L_{k_2}} \notin \mathbb Q$, $t \in \mathbb R$, $\mathfrak n \in \mathfrak N$, and $\alpha_{k_2} \in \mathcal G(T, \mu)$, there exists $\mathfrak r \in \mathfrak N$ with $n_j \leq r_j \leq K + n_j$, $j \in \{k_1, k_2\}$, and $r_j = n_j$ for $j \in \llbracket 1, N\rrbracket \setminus \{k_1, k_2\}$, such that
\[\varepsilon_{k_2, \mathfrak r, 0, t} \leq \eta.\]
\end{lemm}

\begin{prf}
Take $\rho_j > 0$, $j \in \llbracket 1, N_d\rrbracket$, as in Lemma~\ref{LemmIPE} and $K \in \mathbb N$ as in Lemma~\ref{LemmTimeIRhoAlpha}. Define $\eta = \max_{j \in \llbracket 1, N_d\rrbracket} e^{-\rho_j} \in (0, 1)$. Take $k_1 \in \llbracket 1, N\rrbracket$ and $k_2 \in \llbracket 1, N_d\rrbracket $ with ${L_{k_1}}/{L_{k_2}} \notin \mathbb Q$. Let $t \in \mathbb R$,  $\mathfrak n \in \mathfrak N$,  and  $\alpha_{k_2} \in \mathcal G(T, \mu)$. 
Applying Lemma~\ref{LemmTimeIRhoAlpha} at $t - L_{k_2}$, we deduce the existence of $\mathfrak r \in \mathfrak N$ such that
\[t - L_{k_2}- L(\mathfrak r)  \in \mathcal I_{k_2, \rho_{k_2}, \alpha_{k_2}}.\]
By the definition \eqref{DefIRhoAlpha} of $\mathcal I_{k_2, \rho_{k_2}, \alpha_{k_2}}$, this means that
\[\int_{t - L(\mathfrak r) - L_{k_2} + a_{k_2}}^{t - L(\mathfrak r) - L_{k_2} + b_{k_2}} \alpha_{k_2} (s) ds \geq \rho_{k_2}.\]
By \eqref{Varepsilon}, we thus obtain that $\varepsilon_{k_2, \mathfrak r, 0, t} \leq e^{-\rho_{k_2}} \leq \eta$.
\end{prf}

\begin{remk}
Notice that the hypothesis that ${L_{k_1}}/{L_{k_2}} \notin \mathbb Q$ is only used to apply Lemma~\ref{LemmTimeIRhoAlpha}, and thus, by Remark~\ref{RemkRational}, Lemma~\ref{LemmPEEta} still holds true if $k_1$ and $k_2$ are chosen as in Hypothesis~\ref{HypoRational}.
\end{remk}

\begin{figure}[ht]
\centering

\setlength{\unitlength}{1.3cm}

\begin{picture}(8.6, 8.8)(-0.6, -0.6)

\linethickness{0.02\unitlength}%
\textcolor{gray}{\multiput(0, 0)(0   , 0.404508){20}{\Line(0, 0)(7.7, 0       )}}%
\textcolor{gray}{\multiput(0, 0)(0.25, 0       ){31}{\Line(0, 0)(0  , 7.885661)}}%

\linethickness{0.03\unitlength}
\multiput(0.75, 0.404508)(0.25, 0       ){8}{\Line(0, 0)(0, 3.236068)}
\multiput(0.75, 0.404508)(0   , 0.404508){8}{\Line(0, 0)(2, 0       )}

\linethickness{0.03\unitlength}%
\textcolor{red}{\put(0   , 2.831559){\Line(0, 0)(0.25, 0)\Line(0.045, -0.08)(0.205, 0.08)\Line(0.045, 0.08)(0.205, -0.08)}}%
\textcolor{red}{\put(0.25, 0.809017){\Line(0, 0)(0.25, 0)\Line(0.045, -0.08)(0.205, 0.08)\Line(0.045, 0.08)(0.205, -0.08)}}%
\textcolor{red}{\put(0.5 , 7.281153){\Line(0, 0)(0.25, 0)\Line(0.045, -0.08)(0.205, 0.08)\Line(0.045, 0.08)(0.205, -0.08)}}%
\textcolor{red}{\put(0.75, 5.663119){\Line(0, 0)(0.25, 0)\Line(0.045, -0.08)(0.205, 0.08)\Line(0.045, 0.08)(0.205, -0.08)}}%
\textcolor{red}{\put(1.25, 4.045085){\Line(0, 0)(0.25, 0)\Line(0.045, -0.08)(0.205, 0.08)\Line(0.045, 0.08)(0.205, -0.08)}}%
\textcolor{red}{\put(2   , 2.022542){\Line(0, 0)(0.25, 0)\Line(0.045, -0.08)(0.205, 0.08)\Line(0.045, 0.08)(0.205, -0.08)}}%
\textcolor{red}{\put(2.25, 6.067627){\Line(0, 0)(0.25, 0)\Line(0.045, -0.08)(0.205, 0.08)\Line(0.045, 0.08)(0.205, -0.08)}}%
\textcolor{red}{\put(2.75, 4.449593){\Line(0, 0)(0.25, 0)\Line(0.045, -0.08)(0.205, 0.08)\Line(0.045, 0.08)(0.205, -0.08)}}%
\textcolor{red}{\put(3.5 , 0.404508){\Line(0, 0)(0.25, 0)\Line(0.045, -0.08)(0.205, 0.08)\Line(0.045, 0.08)(0.205, -0.08)}}%
\textcolor{red}{\put(3.75, 6.472136){\Line(0, 0)(0.25, 0)\Line(0.045, -0.08)(0.205, 0.08)\Line(0.045, 0.08)(0.205, -0.08)}}%
\textcolor{red}{\put(4   , 2.022542){\Line(0, 0)(0.25, 0)\Line(0.045, -0.08)(0.205, 0.08)\Line(0.045, 0.08)(0.205, -0.08)}}%
\textcolor{red}{\put(4.5 , 5.258610){\Line(0, 0)(0.25, 0)\Line(0.045, -0.08)(0.205, 0.08)\Line(0.045, 0.08)(0.205, -0.08)}}%
\textcolor{red}{\put(5   , 7.281153){\Line(0, 0)(0.25, 0)\Line(0.045, -0.08)(0.205, 0.08)\Line(0.045, 0.08)(0.205, -0.08)}}%
\textcolor{red}{\put(5.25, 3.236068){\Line(0, 0)(0.25, 0)\Line(0.045, -0.08)(0.205, 0.08)\Line(0.045, 0.08)(0.205, -0.08)}}%
\textcolor{red}{\put(5.5 , 0.809017){\Line(0, 0)(0.25, 0)\Line(0.045, -0.08)(0.205, 0.08)\Line(0.045, 0.08)(0.205, -0.08)}}%
\textcolor{red}{\put(6   , 4.854102){\Line(0, 0)(0.25, 0)\Line(0.045, -0.08)(0.205, 0.08)\Line(0.045, 0.08)(0.205, -0.08)}}%
\textcolor{red}{\put(6.5 , 6.067627){\Line(0, 0)(0.25, 0)\Line(0.045, -0.08)(0.205, 0.08)\Line(0.045, 0.08)(0.205, -0.08)}}%
\textcolor{red}{\put(7   , 2.427051){\Line(0, 0)(0.25, 0)\Line(0.045, -0.08)(0.205, 0.08)\Line(0.045, 0.08)(0.205, -0.08)}}%
\textcolor{red}{\put(7.25, 0       ){\Line(0, 0)(0.25, 0)\Line(0.045, -0.08)(0.205, 0.08)\Line(0.045, 0.08)(0.205, -0.08)}}%

\linethickness{0.02\unitlength}
\multiput(7.9, -0.25)(-0.176777, 0.176777){47}{\line(-1, 1){0.0883883}}

\put(0, 0){\circle*{0.14}}

\put(-0.1 , -0.1 ){\makebox(0, 0)[tr]{$O$}}
\put( 0.2 , -0.2 ){\makebox(0, 0)[t]{$L_1$}}
\put(-0.2 ,  0.3 ){\makebox(0, 0)[r]{$L_2$}}
\put( 7.8, -0.3    ){\makebox(0, 0)[t]{$l_t$}}
\end{picture}

\caption{Interpretation of Lemma~\ref{LemmPEEta} in the case $N = 2$ and $N_d=1$.}
\label{FigPEEta}
\end{figure}
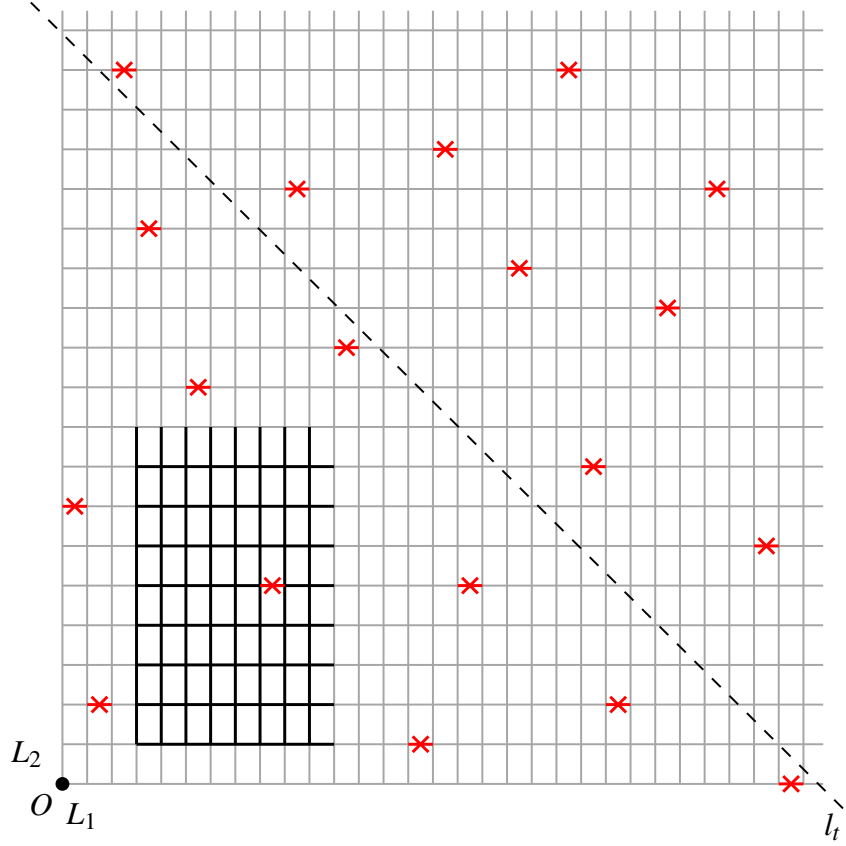

\begin{remk}
Figure~\ref{FigPEEta} helps illustrating Lemma~\ref{LemmPEEta} in the case $N = 2$. Taking $N_d=1$, the lemma states that there exists $K \in \mathbb N$ ($K=7$ in the picture) such that every rectangle of size $(K+1) L_{k_1} \times (K+1) L_{k_2}$ that we place in the grid pictured in Figure~\ref{FigPEEta} contains at least one horizontal segment (highlighted in the picture) where $\varepsilon_{k_2, \mathfrak r, 0, t} \leq \eta$. Let us remark that $K$ and $\eta$ do not depend on  $\alpha_1 \in \mathcal G(T, \mu)$: hence the position of the highlighted segments may change if we change the persistently exciting signal, but we can guarantee that on every rectangle there exists at least one such segment.
\end{remk}

We now apply Lemma~\ref{LemmPEEta} to obtain the following property, which is a preliminary step towards the exponential decay of $\vartheta^{(i)}_{j, \mathfrak n, x, t}$ in $\mathfrak N_c(\rho)$.

\begin{lemm}
\label{LemmVarthetaLeqLambda}
Suppose that Hypotheses~\ref{HypoIrrational} and \ref{HypoM} are satisfied. Let $T \geq \mu > 0$ and $\rho \in (0, 1)$. Then there exist $\lambda \in (0, 1)$ and $K \in \mathbb N^\ast$ such that, for every $\mathfrak n \in \mathfrak N_{c}(\rho)$ with $\min_{i \in \llbracket 1, N\rrbracket} n_i\geq K$, $i, j \in \llbracket 1, N\rrbracket$, $t \in \R$, and $\alpha_k \in \mathcal G(T, \mu)$, $k \in \llbracket 1, N_d\rrbracket$, we have
\[\abs{\vartheta^{(i)}_{j, \mathfrak n, L_j, t}} \leq \lambda \max_{\substack{\mathfrak p \in \mathfrak N \suchthat \abs{\mathfrak p}_{\ell^1} = K \\ p_r \leq n_r,\; \forall r \in \llbracket 1, N\rrbracket \\ s \in \llbracket 1, N\rrbracket \suchthat p_s > 0}} \abs{\vartheta^{(i)}_{s, \mathfrak n - \mathfrak p, L_s, t}}.\]
\end{lemm}

\begin{prf}
Let $\eta \in (0, 1)$ and $K_0 \in \mathbb N^\ast$ be as in Lemma~\ref{LemmPEEta}. Let $L_{\min} = \min_{i \in \llbracket 1, N\rrbracket} L_i$. We take $K = 2 K_0+1$.  

Take $\mathfrak n \in \mathfrak N_{c}(\rho)$ with $\min_{i \in \llbracket 1, N\rrbracket} n_i\geq K$ and let $k_1 \in \llbracket 1, N\rrbracket$, $k_2 \in \llbracket 1, N_d\rrbracket $ be such that ${L_{k_1}}/{L_{k_2}} \notin \mathbb Q$. Since $\mathfrak n \in \mathfrak N_c(\rho)$, one has $n_{k_i} > \rho \abs{\mathfrak n}_{\ell^1}$ for $i \in \{1, 2\}$. Take $i, j \in \llbracket 1, N\rrbracket$, $t \in\R$ and $\alpha_k \in \mathcal G(T, \mu)$ for $k \in \llbracket 1, N_d\rrbracket$. Since $\abs{\mathfrak n}_{\ell^1} \geq \min_{i \in \llbracket 1, N\rrbracket} n_i\geq K$, we can apply Proposition~\ref{PropAverage} and deduce from \eqref{VarthetaAverage} the estimate
\begin{equation}\label{EstimAverage}
\abs{\vartheta^{(i)}_{j, \mathfrak n, L_j, t}}  \leq \Theta \sum_{v \in \Phi_K(\mathfrak n)} \left[\left(\abs{m_{v_1 j}} \prod_{s=2}^K \abs{m_{v_s v_{s-1}}}\right) \left(\prod_{s=1}^K \varepsilon_{v_s, \mathfrak n - \sum_{r=1}^s \mathbf 1_{v_r}, 0, t}\right)\right],
\end{equation}
where
\[\Theta = \max_{\substack{\mathfrak p \in \mathfrak N \suchthat \abs{\mathfrak p}_{\ell^1} = K \\ p_r \leq n_r,\; \forall r \in \llbracket 1, N\rrbracket \\ s \in \llbracket 1, N\rrbracket \suchthat p_s > 0}} \abs{\vartheta^{(i)}_{s, \mathfrak n - \mathfrak p, L_s, t}}.\]

Let us now apply Lemma~\ref{LemmPEEta} to the point $\mathfrak n^\prime = \mathfrak n - K_0 \mathbf 1_{k_1} - (K_0+1) \mathbf 1_{k_2}$. Notice that $\mathfrak n^\prime \in \mathfrak N$ since $n_{k_1},n_{k_2} > K_0$. Hence there exists $\mathfrak r \in \mathfrak N$ with
\[
\begin{gathered}
n_{k_1} - K_0 \leq r_{k_1} \leq n_{k_1}, \\
n_{k_2}  -K_0-1\leq r_{k_2} \leq n_{k_2}-1 , \\
r_j = n_j \text{ for } j \in \llbracket 1, N\rrbracket \setminus\{k_1, k_2\}, \\
\end{gathered}
\]
such that $\varepsilon_{k_2, \mathfrak r, 0, t} \leq \eta$.

We next show that there exists $v_0 = (v_{0, 1}, \dotsc, v_{0, K}) \in \Phi_K(\mathfrak n)$ and $s_0 \in \llbracket 1, K\rrbracket$ such that $v_{0, s_0} = k_2$ and $\mathfrak n - \sum_{r=1}^{s_0} \mathbf 1_{v_{0, r}} = \mathfrak r$. For that purpose, take $v_0 \in \llbracket 1, N\rrbracket^K$ with $v_{0, 1} = v_{0, 2} = \dotsb = v_{0, n_{k_1} - r_{k_1}} = k_1$ and $v_{0, n_{k_1} - r_{k_1} + 1} = v_{0, n_{k_1} - r_{k_1} + 2} = \dotsb = v_{0, K} = k_2$. Such a $v_0$ is well-defined in $\Phi_K(\mathfrak n)$ since $0 \leq n_{k_1} - r_{k_1} \leq K_0 < K$. By construction, $\varphi_{k_1, K}(v_0) = n_{k_1} - r_{k_1} \leq n_{k_1}$, $\varphi_{k_2, K}(v_0) = K - (n_{k_1} - r_{k_1}) \leq K \leq n_{k_2}$ and $\varphi_{k, K}(v_0) = 0\leq n_k$ for $k \in \llbracket 1, N\rrbracket \setminus \{k_1, k_2\}$. Hence $v_0$ is in $\Phi_K(\mathfrak n)$. 
Taking $s_0 = n_{k_1} - r_{k_1} + n_{k_2} - r_{k_2} \in \llbracket 1, K\rrbracket$, we have $v_{0, s_0} = k_2$ and $\mathfrak n - \sum_{r=1}^{s_0} \mathbf 1_{v_{0, r}} = \mathfrak r$.

Let $\delta = \min_{i, j \in \llbracket 1, N\rrbracket} \abs{m_{ij}} > 0$ and $\lambda =  1 - \delta^K(1 - \eta)$. One clearly has that $\lambda$ is in $(0, 1)$, since $\eta,\delta\in(0,1)$.
Using \eqref{SumLeq1}, we get from \eqref{EstimAverage} that  
\begin{align*}
\abs{\vartheta^{(i)}_{j, \mathfrak n, L_j, t}} 
 & \leq \Theta \left[\left(\abs{m_{v_{0, 1} j}} \prod_{s=2}^K \abs{m_{v_{0, s} v_{0, s-1}}}\right) \left(\varepsilon_{v_{0, s_0}, \mathfrak n - \sum_{r=1}^{s_0} \mathbf 1_{v_{0, r}}, 0, t} \prod_{\substack{s=1 \\ s \not = s_0}}^K \varepsilon_{v_{0, s}, \mathfrak n - \sum_{r=1}^s \mathbf 1_{v_{0, r}}, 0, t}\right)\right. \\
 & \hphantom{ = \Theta[} \left. + \sum_{v \in \Phi_K(\mathfrak n) \setminus \{v_0\}} \left(\abs{m_{v_1 j}} \prod_{s=2}^K \abs{m_{v_s v_{s-1}}}\right) \left(\prod_{s=1}^K \varepsilon_{v_s, \mathfrak n - \sum_{r=1}^s \mathbf 1_{v_r}, 0, t}\right)\right] \displaybreak[0] \\
 & \leq \Theta \left[\left(\abs{m_{v_{0, 1} j}} \prod_{s=2}^K \abs{m_{v_{0, s} v_{0, s-1}}}\right) \eta +\sum_{v \in \Phi_K(\mathfrak n) \setminus \{v_0\}} \left(\abs{m_{v_1 j}} \prod_{s=2}^K \abs{m_{v_s v_{s-1}}}\right)\right] \displaybreak[0] \\
 & = \Theta \left[\left(\abs{m_{v_{0, 1} j}} \prod_{s=2}^K \abs{m_{v_{0, s} v_{0, s-1}}} \right)(\eta-1)+\sum_{v \in \Phi_K(\mathfrak n)} \left(\abs{m_{v_1 j}} \prod_{s=2}^K \abs{m_{v_s v_{s-1}}}\right) \right] \displaybreak[0] \\
 & \leq \Theta \left[\delta^K (\eta-1) + \sum_{v \in \Phi_K(\mathfrak n)} \left(\abs{m_{v_1 j}} \prod_{s=2}^K \abs{m_{v_s v_{s-1}}}\right)\right]\leq \Theta \left[1-\delta^K (1-\eta)\right] = \lambda \Theta.
\end{align*}
\end{prf}

We now  obtain the exponential decay of $\vartheta^{(i)}_{j, \mathfrak n, L_j, t}$ in the set $\mathfrak N_{c}(\sigma)$.

\begin{theo}
\label{TheoEstimMiddle}
Suppose that Hypotheses~\ref{HypoIrrational} and \ref{HypoM} are satisfied. Let $T \geq \mu > 0$ and $\sigma \in (0, 1)$. Then there exist $C\ge 1$, $\gamma > 0$ and $K \in \mathbb N^\ast$ such that, for every $\mathfrak n \in \mathfrak N_{c}(\sigma)$ with $\min_{i \in \llbracket 1, N\rrbracket} n_i \geq K$, $i, j \in \llbracket 1, N\rrbracket$, $x \in [0, L_j]$, $t \in \R$, and $\alpha_k \in \mathcal G(T, \mu)$, $k \in \llbracket 1, N_d\rrbracket$, we have
\begin{equation}
\abs{\vartheta^{(i)}_{j, \mathfrak n, x, t}} \leq C e^{-\gamma \abs{\mathfrak n}_{\ell^1}}.
\label{EstimMiddle}
\end{equation}
\end{theo}

\begin{prf}
Take $\rho = \nicefrac{\sigma}{2}$ and let $\lambda \in (0, 1)$ and $K \in \mathbb N^\ast$ be as in Lemma~\ref{LemmVarthetaLeqLambda}. For $\mathfrak n \in \mathfrak N_{c}(\rho)$ with $\min_{i \in \llbracket 1, N\rrbracket} n_i \geq K$, we set
\[q_{\max}(\mathfrak n) = \max\left\{q\in\mathbb N \mid \mathfrak n-\mathfrak r\in \mathfrak N_{c}(\rho)\mbox{ for every }\mathfrak r\in \mathfrak N \mbox{ with }|\mathfrak r|_{\ell^1}=q K\right\}.\]
From Lemma~\ref{LemmVarthetaLeqLambda}, one deduces by an immediate inductive argument that for every $q \in \llbracket 1, q_{\max}(\mathfrak n)\rrbracket$,  $i, j \in \llbracket 1, N\rrbracket$, $t \in\R$, and $\alpha_k \in \mathcal G(T, \mu)$, $k \in \llbracket 1, N_d\rrbracket$, we have
\begin{equation*}
\abs{\vartheta^{(i)}_{j, \mathfrak n, L_j, t}} \leq \lambda^q \max_{\substack{\mathfrak p \in \mathfrak N \suchthat \abs{\mathfrak p}_{\ell^1} = q K \\ p_r \leq n_r,\; \forall r \in \llbracket 1, N\rrbracket \\ s \in \llbracket 1, N\rrbracket \suchthat p_s > 0}} \abs{\vartheta^{(i)}_{s, \mathfrak n - \mathfrak p, L_s, t}}.
\end{equation*}

By Corollary~\ref{CoroLeq1}, it holds $\abs{\vartheta^{(i)}_{j, \mathfrak n, L_j, t}} \leq 1$ for every $i, j \in \llbracket 1, N\rrbracket$, $\mathfrak n \in \mathfrak N$ and $t \in\R$. Therefore,  for every $\mathfrak n\in\mathfrak N_c(\rho)$ with $\min_{i \in \llbracket 1, N\rrbracket}\mathfrak n_i\geq K$, $i, j \in \llbracket 1, N\rrbracket$,  $t \in\R$ and $\alpha_k \in \mathcal G(T, \mu)$, $k \in \llbracket 1, N_d\rrbracket$, we obtain
\begin{equation}
\abs{\vartheta^{(i)}_{j, \mathfrak n, L_j, t}} \leq \lambda^{q_{\max}(\mathfrak n)}.
\label{EstimVarthetaLambdaFinal}
\end{equation}

Notice now that, by definition of $q_{\max}$, one also has that
$$q_{\max}(\mathfrak n) + 1 \geq \frac 1 K \min\{|\mathfrak n-\mathfrak r|_{\ell^1}\mid \mathfrak r\in \mathfrak N_b(\rho)\},$$
where, according to Definition~\ref{def_Nb}, $\mathfrak N_b(\rho)=\mathfrak N\setminus \mathfrak N_c(\rho)$. One deduces at once that there exists $\xi>0$ such that, for every $\mathfrak n \in \mathfrak N_{c}(\sigma)$,
\begin{equation}
q_{\max}(\mathfrak n)+1\geq \xi \abs{\mathfrak n}_{\ell^1}.
\label{RhoEstimQ}
\end{equation}
Since $\lambda \in (0, 1)$, setting $\gamma = - \xi \log \lambda > 0$ and $C=1/\lambda$ one concludes by inserting \eqref{RhoEstimQ} into \eqref{EstimVarthetaLambdaFinal} and then using \eqref{ThetaIX-PE}.
\end{prf}


\subsection{Exponential convergence of the solutions}

To conclude the proof of Theorem~\ref{MainTheoIntro}, it suffices to combine Proposition~\ref{LemmConvCoefConvSol} with the estimates of $\vartheta^{(i)}_{j, \mathfrak n, x, t}$ for $\mathfrak n \in \mathfrak N_{b}(\sigma)$ and $\mathfrak n \in \mathfrak N_c(\rho)$ 
given by Theorems~\ref{TheoVarthetaNb} and \ref{TheoEstimMiddle}, respectively.

\begin{prf}[Proof of Theorem~\ref{MainTheoIntro}]
Let $C_1, \gamma_1 > 0$ and $\rho > 0$ be as in Theorem~\ref{TheoVarthetaNb}. Take $\sigma = \rho$ in Theorem~\ref{TheoEstimMiddle} and let $C_2, \gamma_2 > 0$ and $K \in \mathbb N^\ast$ be as in that theorem. Let
\[\gamma = \min\{\gamma_1, \gamma_2\}, \qquad C = \max\{C_1, C_2,e^{\gamma K/\rho}\}.\]
Take $i, j \in \llbracket 1, N\rrbracket$, $\mathfrak n \in \mathfrak N$, $x \in [0, L_j]$, $t \in\R$ and $\alpha_k \in \mathcal G(T, \mu)$, $k \in \llbracket 1, N_d\rrbracket$. If $\mathfrak n \in \mathfrak N_{b}(\rho)$ or $\mathfrak n \in \mathfrak N_{c}(\rho)$ with $\min_{i \in \llbracket 1, N\rrbracket} n_i\ge K$, then one concludes directly from Theorems~\ref{TheoVarthetaNb} and \ref{TheoEstimMiddle}. Finally, if $\mathfrak n \in \mathfrak N_{c}(\rho)$ with $\min_{i \in \llbracket 1, N\rrbracket} n_i< K$, note that $|\mathfrak n|_{\ell^1}\le K/\rho$. Then, by Corollary~\ref{CoroLeq1}, one has
\[\abs{\vartheta^{(i)}_{j, \mathfrak n, x, t}} \leq 1 \leq C e^{-\gamma K/\rho} \leq  C e^{-\gamma \abs{\mathfrak n}_{\ell^1}}.\]
Theorem~\ref{MainTheoIntro} now follows from Proposition~\ref{LemmConvCoefConvSol}.
\end{prf}

\begin{remk}
By Remark~\ref{RemkRational}, Hypothesis~\ref{HypoIrrational} can be replaced by Hypothesis~\ref{HypoRational} in Lem\-ma~\ref{LemmVarthetaLeqLambda} and Theorem~\ref{TheoEstimMiddle}, and so the same is also true for Theorem~\ref{MainTheoIntro}. We recall that the case $p = +\infty$ also follows from Proposition~\ref{LemmConvCoefConvSol} thanks to Remark~\ref{RemkInfty}.
\end{remk}

\appendix


\section{Well-posedness of the Cauchy problems \eqref{TranspOperatorPE}, \eqref{TranspOperatorUndamped} and \eqref{TranspOperatorAlwaysActive}}
\label{AppendWellPosed}

We are interested in this section in the proof of Theorem~\ref{TheoWellPosed}, which states the well-posedness of the Cauchy problem \eqref{TranspOperatorPE}. This is done in two steps. First, we show that the operator $A$ defined in \eqref{DefA} is the generator of a strongly continuous semigroup $\{e^{tA}\}_{t \geq 0}$, thus establishing the well-posedness of the undamped system. We then consider the operator $B(t) = \sum_{i=1}^{N_d} \alpha_i(t) B_i$ as a bounded time-dependent perturbation of $A$ in order to conclude the well-posedness of \eqref{TranspOperatorPE}.


\subsection{Preliminaries}

Let $\mathsf X$ be a Banach space and let $A: D(A) \subset \mathsf X \to \mathsf X$ be an operator in $\mathsf X$. The definitions of strong and weak solutions for the Cauchy problem associated with $A$ can be found for instance in \cite{Pazy1983Semigroups}. Recall that if $A$ is densely defined with a non-empty resolvent set, then the Cauchy problem associated with $A$ has a unique strong solution for each initial condition $z_0 \in D(A)$ if and only if $A$ is the generator of a strongly continuous semigroup $\{e^{tA}\}_{t \geq 0}$. In this case, the solution is given by $z(t) = e^{tA} z_0$ (see, for instance, \cite[Chapter 4, Theorem 1.3]{Pazy1983Semigroups}). Then, $t \mapsto e^{tA} z_0$ is a well-defined continuous function for every $z_0 \in \mathsf X$ and it is the unique weak solution of  the Cauchy problem associated with $A$.

\begin{defi}
A family of operators $\{T(t, s)\}_{t \geq s \geq 0} \subset \mathcal L(\mathsf X)$ is an \emph{evolution family} on $\mathsf X$ if
\begin{enumerate}[label={\bf\roman*.}, ref={\roman*}]
\item $T(s, s) = \id_{\mathsf X}$ for every $s \geq 0$,
\item $T(t, s) = T(t, \tau) T(\tau, s)$ for every $t \geq \tau \geq s \geq 0$,
\item for every $z \in \mathsf X$, $(t, s) \mapsto T(t, s) z$ is continuous for every $t \geq s \geq 0$.
\end{enumerate}
An evolution family is \emph{exponentially bounded} if it satisfies further the following property.
\begin{enumerate}[label={\bf\roman*.}, ref={\roman*}, resume]
\item There exist $M \geq 1$ and $\omega \in \mathbb R$ such that $\norm{T(t, s)}_{\mathcal L(\mathsf X)} \leq M e^{\omega (t - s)}$ for every $t \geq s \geq 0$.
\end{enumerate}
\end{defi}
For references on evolution families see, for instance, \cite{Chicone1999Evolution,Phillips1953Perturbation,Kato1973Linear}. We are interested here in a family of the form $A(t) = A + B(t)$ where $A$ is the generator of a strongly continuous semigroup and $B \in L^\infty(\mathbb R_+, \mathcal L(\mathsf X))$. We shall use here the following notions of solution.

\begin{defi}
\label{DefSolCauchyNonAuto}
Consider the problem
\begin{equation}
\left\{
\begin{aligned}
\dot z(t) & = (A + B(t)) z(t), & \quad & t \geq s \geq 0, \\
z(s) & = z_0, & &
\end{aligned}
\right.
\label{CauchyNonAutoAB}
\end{equation}
where $A$ is the generator of a strongly continuous semigroup $\{e^{tA}\}_{t \geq 0}$ and $B \in L^\infty(\mathbb R_+, \mathcal L(\mathsf X))$.

\begin{enumerate}[label={\bf \roman*.}, ref={\roman*}]
\item We say that $z: \left[s, +\infty\right) \to \mathsf X$ is a \emph{regular solution} of \eqref{CauchyNonAutoAB} if $z$ is continuous, $z(s) = z_0$, $z(t) \in D(A)$ for every $t \geq s$, $z$ is differentiable for almost every $t \geq s$, $\dot z \in L^\infty_{\text{\rm loc}}(\left[s, +\infty\right), \mathsf X)$ and $\dot z(t) = (A + B(t)) z(t)$ for almost every $t \geq s$.

\item We say that $z: \left[s, +\infty\right) \to \mathsf X$ is a \emph{mild solution} of \eqref{CauchyNonAutoAB} if $z$ is continuous and, for every $t \geq s$, we have
\[
z(t) = e^{(t-s)A} z_0 + \int_s^t e^{(t-\tau) A} B(\tau) z(\tau) d\tau.
\]
\end{enumerate}
\end{defi}
Here, the integrals of $\mathsf X$-valued functions should be understood as Bochner integrals; see, for instance, \cite{Yosida1980Functional}. In the following proposition, we summarize the main facts needed for the present paper.

\begin{prop}
\label{PropRegularMild}
Let $A$ be the generator of a strongly continuous semigroup $\{e^{tA}\}_{t \geq 0}$ and $B \in L^\infty(\mathbb R_+, \mathcal L(\mathsf X))$. Then, the following holds true: 

\begin{enumerate}[label={\bf\roman*.}, ref={\roman*}]
\item\label{RegularMild} every regular solution of \eqref{CauchyNonAutoAB} is also a mild solution;

\item\label{IntegralEquation} there exists a unique family $\{T(t, s)\}_{t \geq s \geq 0}$ of bounded operators in $\mathsf X$ such that $(t, s) \mapsto T(t, s) z$ is continuous for every $z \in \mathsf X$ and
\begin{equation}
T(t, s) z = e^{(t-s)A} z + \int_s^t e^{(t-\tau) A} B(\tau) T(\tau, s) z d\tau, \qquad \forall z \in \mathsf X.
\label{IntegralEquationT}
\end{equation}
Furthermore, this family is an exponentially bounded evolution family;

\item\label{MildSol} for every $z_0 \in \mathsf X$, \eqref{CauchyNonAutoAB} admits a unique mild solution $z$, given by $z(t) = T(t, s) z_0$.
\end{enumerate}
\end{prop}


\subsection{Proof of Theorem~\ref{TheoWellPosedUndamped}}

Since $\sum_{i=1}^{N_d} B_i$ is a bounded operator, it suffices to show Theorem~\ref{TheoWellPosedUndamped} for $A$.

\begin{prop}
\label{PropFermeDensementDefini}
Let $p \in \left[1, +\infty\right)$. The operator $A$ is closed and densely defined. Moreover, $D(A)$ endowed with the graph norm is a Banach space compactly embedded in $\mathsf X_p$.
\end{prop}

\begin{prf}
The proposition follows straightforwardly by the remark that the graph norm on $D(A)$ coincides with the usual norm in $\prod_{i=1}^N W^{1, p}(0, L_i)$, that is,
\[\norm{z}_{D(A)}^p = \sum_{i=1}^N \left(\norm{u_i}_{L^p(0, L_i)}^p + \norm{u_i^\prime}_{L^p(0, L_i)}^p\right)\]
for $z = (u_1, \dotsc, u_N) \in D(A)$.
\end{prf}

\begin{prop}
\label{PropResolventNonEmpty}
Let $p \in \left[1, +\infty\right)$. The resolvent set $\rho(A)$ of $A$ is nonempty.
\end{prop}

\begin{prf}
Since $A$ is closed, we have $\lambda \in \rho(A)$ if and only if $\lambda - A$ is a bijection from $D(A)$ to $\mathsf X_p$. 
A direct computation based on explicit formulas yields that $\lambda - A$ is a bijection as soon as $\lambda \in \mathbb R$ with $\lambda > \frac{\log\abs{M}_{\ell^2}}{L_{\min}}$, where $L_{\min} = \min_{i \in \llbracket 1, N \rrbracket} L_i$.
\end{prf}

We now turn to a result of existence of solutions of \eqref{TranspOperatorUndamped} when $z_0 \in D(A)$.

\begin{theo}
\label{TheoExistUndamped}
For every $z_0 = (u_{1, 0}, \dotsc, u_{N, 0}) \in D(A)$, there exists a unique strong solution $z = (u_1, \dotsc, u_N) \in \mathcal C^0(\mathbb R_+, D(A)) \cap \mathcal C^1(\mathbb R_+, \mathsf X_p)$ of \eqref{TranspOperatorUndamped}.
\end{theo}

\begin{prf}
Let $T_0 > 0$ be such that $T_0 <L_{\min}$. Note that it suffices to show the theorem for solutions in $\mathcal C^0([0, T_0], D(A)) \cap \mathcal C^1([0, T_0], \mathsf X_p)$, since $T_0$ does not depend on $z_0 \in D(A)$. Let $z_0 = (u_{1, 0}, \dotsc, u_{N, 0}) \in D(A)$. It follows easily from the transport equation and the transmission condition \eqref{TransmissionMatrix} that a solution $t\mapsto z(t)=(u_1(t),\dotsc,u_N(t))$ of \eqref{TranspOperatorUndamped} necessarily satisfies
\begin{equation}
u_i(t, x) = 
\begin{dcases*}
\sum_{j=1}^N m_{ij} u_{j, 0}(L_j - t + x) & if $0 \leq x \leq t$, \\
u_{i, 0}(x - t) & if $x > t$.
\end{dcases*}
\label{EqSol}
\end{equation}
Conversely, if $z=(u_1,\dotsc,u_N)$ is given by \eqref{EqSol}, then it solves \eqref{TranspOperatorUndamped} and has $z_0$ as initial condition. Moreover, one checks by direct computations that $z$ fulfills the required regularity properties. 
\end{prf}

\begin{prf}[Proof of Theorem~\ref{TheoWellPosedUndamped}]
From Propositions~\ref{PropFermeDensementDefini} and \ref{PropResolventNonEmpty} and Theorem~\ref{TheoExistUndamped}, we obtain that $A$ generates a strongly continuous semigroup $\{e^{tA}\}_{t \geq 0}$ (see, for instance, \cite[Chapter 4, Theorem 1.3]{Pazy1983Semigroups}). Since $\sum_{i=1}^{N_d} B_i \in \mathcal L(\mathsf X_p)$, $A + \sum_{i=1}^{N_d} B_i$  also generates a strongly continuous semigroup (see \cite[Chapter 3, Theorem 1.1]{Pazy1983Semigroups}).
\end{prf}

\begin{remk}
In the particular case $p = 2$ and $\abs{M}_{\ell^2} \leq 1$, one may conclude more easily that $A$ is the generator of a strongly continuous semigroup $\{e^{tA}\}_{t \geq 0}$, without having to construct the explicit formula \eqref{EqSol} for the solution as we did in Theorem~\ref{TheoExistUndamped}. Indeed, a straightforward computation shows that, for any $M \in \mathcal M_N(\mathbb R)$, the adjoint operator $A^\ast$ of $A$ is given by
\[
\begin{gathered}
D(A^\ast) = \left\{(u_1, \dotsc, u_N) \in \prod_{i=1}^N H^1(0, L_i) \middlesuchthat u_i(L_i) = \sum_{j=1}^N m_{ji} u_j(0) \right\}, \\
A^\ast(u_1, \dotsc, u_N) = \left(\frac{d u_1}{dx}, \dotsc, \frac{d u_N}{dx}\right).
\end{gathered}
\]
Also, for any $z = (u_1, \dotsc, u_N) \in D(A)$, we have
\[\scalprod{z}{Az} = -\sum_{i=1}^N \int_0^{L_i} u_i u_i^\prime = \frac{1}{2} \sum_{i=1}^N (u_i(0)^2 - u_i(L_i)^2) \leq \frac{\abs{M}_{\ell^2}^2 - 1}{2} \sum_{i=1}^N u_i(L_i)^2,\]
since, by \eqref{TransmissionMatrix}, we have $\sum_{i=1}^N u_i(0)^2 \leq \abs{M}_{\ell^2}^2 \sum_{i=1}^N u_i(L_i)^2$. Thus, if $\abs{M}_{\ell^2} \leq 1$, we have $\scalprod{z}{Az} \leq 0$ for every $z \in D(A)$, so that $A$ is dissipative. A similar computation holds for $A^\ast$, with $M$ replaced by its transpose $M^\transp$, showing that $A^\ast$ is also dissipative. Hence $A$ generates a strongly continuous semigroup of contractions $\{e^{tA}\}_{t \geq 0}$ (see, for instance, \cite[Chapter 1, Theorem 4.4]{Pazy1983Semigroups}).
\end{remk}


\subsection{Proof of Theorem~\ref{TheoWellPosed}}

Thanks to Theorem~\ref{TheoWellPosedUndamped} and Items \ref{IntegralEquation} and \ref{MildSol} of Proposition~\ref{PropRegularMild}, there exists a unique evolution family $\{T(t, s)\}_{t \geq s \geq 0}$ associated with \eqref{TranspOperatorPE} such that, for every $z_0 \in \mathsf X_p$, $t \mapsto T(t, s) z_0$ is the unique mild solution of \eqref{TranspOperatorPE} with initial condition $z(s) = z_0$. In order to complete the proof of Theorem~\ref{TheoWellPosed}, it suffices to show that this solution is actually regular when $z_0 \in D(A)$. To do so, we study the explicit formula for the solutions of \eqref{TranspOperatorPE} for small time, as we did with the undamped system \eqref{TranspOperatorUndamped} in Theorem~\ref{TheoExistUndamped}.

\begin{prf}[Proof of Theorem~\ref{TheoWellPosed}]
Let $T_0 > 0$ be such that $T_0 <L_{\min}$. As in Theorem~\ref{TheoExistUndamped}, it suffices to show the theorem for solutions in $\mathcal C^0 (\left[s, s + T_0\right], \mathsf X_p)$. Since the class $L^\infty(\R,[0,1])$ is invariant by time-translation, we can also suppose without loss of generality that $s = 0$. In order to simplify the notations, we define the function $\varphi_i : \mathbb R_+ \times \mathbb R \to \mathbb R_+^\ast$ for $i \in \llbracket 1, N\rrbracket$ by
\[\varphi_i(t, x) = e^{-\int_0^t \alpha_i(s) \chi_i(x - t + s) ds},\]
where we extend the function $\chi_i$ to $\mathbb R$ by $0$ outside its interval of definition $[0, L_i]$. (In particular, $\varphi_i \equiv 1$ for $i \in \llbracket N_d + 1, N\rrbracket$.) We have that both $\varphi_i$ and ${1}/{\varphi_i}$ belong to $ L^\infty(\mathbb R_+ \times \mathbb R) \cap \mathcal C^0(\mathbb R_+ \times \mathbb R)$ and $ W^{1, \infty}(\mathbb R_+ \times \mathbb R)$.

Let $z_0 = (u_{1, 0}, \dotsc, u_{N, 0}) \in D(A)$. We claim that, for $0 \leq t \leq T_0$, the function $t \mapsto z(t) = (u_1(t), \dotsc, u_N(t))$ given by
\[
u_i(t, x) = 
\begin{dcases*}
\frac{\varphi_i(t, x)}{\varphi_i(t - x, 0)} \sum_{j=1}^{N} m_{ij} \varphi_j(t - x, L_j) u_{j, 0}(x + L_j - t), & if $0 \leq x \leq t$, \\
\varphi_i(t, x) u_{i, 0}(x - t), & if $x > t$,
\end{dcases*}
\]
is in $\mathcal C^0([0, T_0], \mathsf X_p)$, $z(0) = z_0$, $z(t) \in D(A)$ for every $t \in [0, T_0]$, $z$ is differentiable for almost every $t \in [0, T_0]$, $\dot z \in L^\infty([0, T_0], \mathsf X_p)$ and $\dot z(t) = A z(t) + \sum_{i=1}^{N_d} \alpha_i(t) B_i z(t)$ for almost every $t \in [0, T_0]$. Indeed, it is clear that $z$ is well-defined, $z(0) = z_0$, and $z(t) \in \mathsf X_p$ for every $t \in [0, T_0]$. It is also clear, thanks to the regularity properties of $\varphi_i$, that, for every $t \in [0, T_0]$, $u_i(t) \in W^{1, p}(0, t)$ and $u_i(t) \in W^{1, p}(t, T_0)$, and, since $x \mapsto u_i(t, x)$ is continuous at $x = t$, we conclude that $u_i(t) \in W^{1, p}(0, L_i)$. Furthermore,
\[u_i(t, 0) = \sum_{j=1}^N m_{ij} \varphi_j(t, L_j) u_{j, 0}(L_j - t) = \sum_{j=1}^N m_{ij} u_j(t, L_j),\]
and thus $z(t) \in D(A)$ for every $t \in [0, T_0]$. By the same argument, we also obtain that $u_i(\cdot, x) \in W^{1, p}(0, T_0)$ for every $x \in [0, L_i]$. Computing $u_i(t + h, x) - u_i(t, x)$ for $t, t + h \in [0, T_0]$ also shows, by a straightforward estimate, that $\norm{u_i(t + h) - u_i(t)}_{L^p(0, L_i)} \to 0$ as $h \to 0$, and thus $z \in \mathcal C^0([0, T_0], \mathsf X_p)$.

Since $u_i(\cdot, x) \in W^{1, p}(0, T_0)$ for every $x \in [0, L_i]$, one can also compute $\partial_t u_i(t, x)$, and it is easy to verify that $\partial_t u_i \in L^\infty([0, T_0], L^p(0, L_i))$. Hence $z$ is differentiable almost everywhere, with $\dot z = (\partial_t u_1, \dotsc, \partial_t u_N) \in L^\infty([0, T_0], \mathsf X_p)$.

Notice now that $\dot z(t) - A z(t) = (\partial_t u_1(t) + \partial_x u_1(t), \dotsc, \partial_t u_N(t) + \partial_x u_N(t))$ is given by
\[
\partial_t u_i(t, x) + \partial_x u_i(t, x) = 
\left\{
\begin{aligned}
& \frac{\partial_t \varphi_i(t, x) + \partial_x \varphi_i(t, x)}{\varphi_i(t - x, 0)} \sum_{j=1}^{N} m_{ij} \varphi_j(t - x, L_j) u_{j, 0}(x + L_j - t), \\
& \begin{aligned}
& & \qquad\qquad\qquad & \text{ if } 0 \leq x \leq t, \\
& [\partial_t \varphi_i(t, x) + \partial_x \varphi_i(t, x)] u_{i, 0}(x - t), & & \text{ if } x > t,
\end{aligned}
\end{aligned}
\right.
\]
and one can compute that $\partial_t \varphi_i(t, x) + \partial_x \varphi_i(t, x) = - \alpha_i(t) \chi_i(x) \varphi_i(t, x)$ almost everywhere, so that
\[
\partial_t u_i(t, x) + \partial_x u_i(t, x) = -\alpha_i(t) \chi_i(x) u_i(t, x).
\]
Thus $\dot z(t) - A z(t) = - \sum_{i=1}^{N_d} \alpha_i(t) B_i z(t)$ for almost every $t \in [0, T_0]$, which concludes the proof of existence. Uniqueness results from the fact that  every regular solution is in particular a mild solution, which is unique, according to Proposition~\ref{PropRegularMild}.
\end{prf}


\section{Asymptotic behavior of \eqref{TwoCircles}}
\label{AppendAsymptotic}

We consider here System~\eqref{TwoCircles} from Example~\ref{ExplTwoCircles}. Let $\mathsf X_2$ be the Hilbert space $\mathsf X_2 = L^2(0, L_1) \times L^2(0, L_2)$. The goal of this section is to prove Theorem~\ref{TheoAsympTwoCircles} concerning the asymptotic behavior of \eqref{TwoCircles} when $\chi \equiv 0$, and also to show the existence of periodic solutions to \eqref{TwoCircles} with a persistent damping when ${L_1}/{L_2} \in \mathbb Q$ and $b - a$ is small enough. The proof of Theorem~\ref{TheoAsympTwoCircles}.\ref{CaseIrrational} being based on LaSalle Principle, we recall its formulation in a Banach space in Section~\ref{AppendLyapunov}.


\subsection{LaSalle Principle in a Banach space}
\label{AppendLyapunov}

In this section, $\mathsf X$ denotes a Banach space and $A: D(A) \subset \mathsf X \to \mathsf X$ is a linear operator in $\mathsf X$ that generates a strongly continuous semigroup $\{e^{tA}\}_{t \geq 0}$.

\begin{defi}
\label{DefiLaSalle}
\begin{enumerate}[label={\bf\roman*.}, ref={\roman*}]
\item\label{DefiOmegaLimit} For $z_0 \in \mathsf X$, the {\em $\omega$-limit set $\omega(z_0)$} is the set of $z \in \mathsf X$ such that there exists a nondecreasing sequence $(t_n)_{n \in \mathbb N}$ in $\mathbb R_+$ with $t_n \to +\infty$ as $n \to \infty$ such that $e^{t_n A} z_0 \to z$ in $\mathsf X$ as $n \to \infty$. A set $M \subset \mathsf X$ is {\em positively invariant} under $\{e^{tA}\}_{t \geq 0}$ if, for every $z_0 \in M$ and $t \geq 0$, we have $e^{tA} z_0 \in M$. For $E \subset \mathsf X$, the {\em maximal positively invariant subset} $M$ of $E$ is the union of all positively invariant sets contained in $E$.

\item A {\em Lyapunov function} for $\{e^{tA}\}_{t \geq 0}$ is a continuous function $V: \mathsf X \to \mathbb R_+$ such that
\[\dot V(z) = \limsup_{t \to 0+}\frac{V(e^{tA}z) - V(z)}{t} \leq 0, \qquad \forall z \in \mathsf X.\]
\end{enumerate}
\end{defi}

The following results can be found in \cite{Hale1969Dynamical, Henry1981Geometric, Slemrod1988LaSalle}.

\begin{theo}
\label{TheoOmegaLimLaSalle}
\begin{enumerate}[label={\bf\roman*.}, ref={\roman*}]
\item Suppose that $\{e^{tA}z_0 \suchthat t \geq 0\}$ is precompact in $\mathsf X$. Then $\omega(z_0)$ is a nonempty, compact, connected, positively invariant set. \label{TheoOmegaLim}
\item Let $V$ be a Lyapunov function on $\mathsf X$, define $E = \{z \in \mathsf X \suchthat \dot V(z) = 0\}$ and let $M$ be the maximal positively invariant subset of $E$. If $\{e^{tA} z_0 \suchthat t \geq 0\}$ is precompact in $\mathsf X$, then $\omega(z_0) \subset M$. \label{TheoLaSalle}
\end{enumerate}
\end{theo}


\subsection{Asymptotic behavior when ${L_1}/{L_2} \notin \mathbb Q$}
\label{SecIrrational}

Let us turn to the proof of Theorem~\ref{TheoAsympTwoCircles}.\ref{CaseIrrational}. We consider the undamped system
\begin{equation}
\left\{
\begin{aligned}
 & \partial_t u_i(t, x) + \partial_x u_i(t, x) = 0, & \quad & t \in \mathbb R_+,\; x \in [0, L_i],\; i \in \{1, 2\}, \\
 & u_1(t, 0) = u_2(t, 0) = \frac{u_1(t, L_1) + u_2(t, L_2)}{2}, & & t \in \mathbb R_+, \\
 & u_i(0, x) = u_{i, 0}(x), & & x \in [0, L_i],\; i \in \{1, 2\}.
\end{aligned}
\right.
\label{AppTwoCirclesUndamped}
\end{equation}

Let $V: \mathsf X_2 \to \mathbb R$ be the function $V(z) = \norm{z}_{\mathsf X_2}^2$ and $A$ be the operator given in Definition~\ref{DefABi} in the case $p = 2$, $N = 2$ and $m_{ij} = \nicefrac{1}{2}$ for $i, j \in \{1, 2\}$, which is the operator associated with System~\eqref{AppTwoCirclesUndamped}. By Theorem~\ref{TheoWellPosedUndamped}, $A$ is the generator of a strongly continuous semigroup $\{e^{tA}\}_{t \geq 0}$.

\begin{lemm}
\label{LemmLyapunov}
The function $V$ is a Lyapunov function for $\{e^{tA}\}_{t \geq 0}$. If $z = (u_1, u_2) \in D(A)$, we have
\[\dot V(z) = 2\scalprod{z}{Az} = -\frac{(u_1(L_1) - u_2(L_2))^2}{2}.\]
\end{lemm}

\begin{prf}
Notice first that, for $z = (u_1, u_2) \in D(A)$, we have $\scalprod{z}{Az} = - {(u_1(L_1) - u_2(L_2))^2}/{4} \leq 0$. Take $z \in D(A)$, so that $t \mapsto e^{tA} z$ is continuously differentiable in $\mathbb R_+$; thus $t \mapsto V(e^{tA} z)$ is continuously differentiable in $\mathbb R_+$ with $\frac{d}{dt} V(e^{tA} z) = 2 \scalprod{e^{tA} z}{A e^{tA} z}$. Hence, for every $z \in D(A)$,  $\dot V(z) = 2 \scalprod{z}{A z} \leq 0$. This also shows that $\norm{e^{tA} z}_{\mathsf X_2} \leq \norm{z}_{\mathsf X_2}$ for every $z \in D(A)$ and $t \geq 0$, and, by the density of $D(A)$ in $\mathsf X_2$, we obtain that $\norm{e^{tA} z}_{\mathsf X_2} \leq \norm{z}_{\mathsf X_2}$ for every $z \in \mathsf X_2$ and $t \geq 0$, i.e., $\{e^{tA}\}_{t \geq 0}$ is a contraction semigroup. Thus $\dot V(z) \leq 0$ for every $z \in \mathsf X_2$, and so $V$ is a Lyapunov function for $\{e^{tA}\}_{t \geq 0}$.
\end{prf}

It is then immediate to prove the following.

\begin{lemm}
\label{LemmOrbitePrecompacte}
For every $z \in D(A)$, $\{e^{tA} z \suchthat t \geq 0\}$ is precompact in $\mathsf X_2$ and $\omega(z) \subset D(A)$.
\end{lemm}

Set $E = \{z \in \mathsf X_2 \suchthat \dot V(z) = 0\}$ and let $M$ be the maximal positively invariant subset of $E$.

\begin{lemm}
\label{LemmDAcMConstant}
Suppose ${L_1}/{L_2} \notin \mathbb Q$. Then
\[D(A) \cap M = \{(\lambda, \lambda) \in L^2(0, L_1) \times L^2(0, L_2) \suchthat \lambda \in \mathbb R\},\]
i.e., $D(A) \cap M$ is the set of constant functions on $L^2(0, L_1) \times L^2(0, L_2)$.
\end{lemm}

\begin{prf}
Take $z_0 = (u_{1, 0}, u_{2, 0}) \in D(A) \cap M$. By the positive invariance of $M$, $e^{tA} z_0 \in D(A) \cap M$ for every $t \geq 0$.

Let us denote $z(t) = (u_1(t), u_2(t)) = e^{tA} z_0$, which is a strong solution of \eqref{AppTwoCirclesUndamped} with initial condition $z_0$. Since $z(t) \in M$, we have $\dot V(z(t)) = 0$ for every $t \geq 0$, which means, by Lemma~\ref{LemmLyapunov}, that $u_1(t, L_1) = u_2(t, L_2)$ for every $t \geq 0$. Then we have that
\[u_1(t, 0) = u_2(t, 0) = \frac{u_1(t, L_1) + u_2(t, L_2)}{2} = u_1(t, L_1) = u_2(t, L_2), \qquad \forall t \geq 0.\]
Without any loss of generality, we can suppose that $L_1 \leq L_2$. For $t \geq L_1$ and $x \in [0, L_1]$, we have that
\[u_1(t, x) = u_1(t-x, 0) = u_2(t-x, 0) = u_2(t, x),\]
and so $u_1(t, x) = u_2(t, x)$. Now, for $t \geq L_1$ and $x \in [0, L_1]$, we have that
\[u_1(t + L_1, x) = u_1(t + L_1 - x, 0) = u_1(t + L_1 - x, L_1) = u_1(t - x, 0) = u_1(t, x)\]
and thus $\left[L_1, +\infty\right) \ni t \mapsto u_1(t, x)$ is a $L_1$-periodic function for every $x \in [0, L_1]$. Similarly, for $t \geq L_2$ and $x \in [0, L_2]$, we have that
\[u_2(t + L_2, x) = u_2(t + L_2 - x, 0) = u_2(t + L_2 - x, L_2) = u_2(t - x, 0) = u_2(t, x)\]
and thus $\left[L_2, +\infty\right) \ni t \mapsto u_2(t, x)$ is a $L_2$-periodic function for every $x \in [0, L_2]$. Since $u_1(t, x) = u_2(t, x)$ for $t \geq L_1$, $x \in [0, L_1]$, we obtain that $\left[L_2, +\infty\right) \ni t \mapsto u_1(t, x)$ is both $L_1$-periodic and $L_2$-periodic for every $x \in [0, L_1]$, and the fact that ${L_1}/{L_2} \notin \mathbb Q$ thus implies that $\left[L_2, +\infty\right) \ni t \mapsto u_1(t, x) = u_2(t, x) \eqqcolon \lambda(x)$ is constant for every $x \in [0, L_1]$. Clearly, $\lambda(x)$ does not depend on $x$ since
\[
\lambda(x) = u_1(t, x) = u_1(t-x, 0) = \lambda(0), \qquad \forall t \geq L_2 + L_1,\; \forall x \in [0, L_1],
\]
and so we shall denote this constant value simply by $\lambda$. We thus have that
\[u_1(t, x) = u_2(t, x) = \lambda, \qquad \forall t \geq L_2,\; \forall x \in [0, L_1].\]
We deduce at once that  $u_{1, 0}$ and $u_{2, 0}$ are both equal to the constant function $\lambda$ and hence
\[D(A) \cap M \subset \{(\lambda, \lambda) \in L^2(0, L_1) \times L^2(0, L_2) \suchthat \lambda \in \mathbb R\}.\]
The converse inclusion is trivial and this concludes the proof.
\end{prf}

Suppose now that $z_0 \in D(A)$. By Lemma~\ref{LemmOrbitePrecompacte} and Theorem~\ref{TheoOmegaLimLaSalle}.\ref{TheoLaSalle}, we have that $\omega(z_0) \subset D(A) \cap M$ and thus, by Lemma~\ref{LemmDAcMConstant}, if ${L_1}/{L_2} \notin \mathbb Q$, we get that every function in $\omega(z_0)$ is constant. We now wish to show that $\omega(z_0)$ contains only one point in $\mathsf X_2$, which will imply that $e^{tA} z_0$ converges to this function as $t \to +\infty$. To do so, we study a conservation law for \eqref{AppTwoCirclesUndamped}.

We define $U: \mathsf X_2 \to \mathbb R$ by
\[U(u_1, u_2) = \frac{1}{L_1 + L_2} \left(\int_0^{L_1} u_1(x) dx + \int_0^{L_2} u_2(x) dx\right).\]
Notice that $U$ is well defined and continuous in $\mathsf X_2$ since we have the continuous embedding $\mathsf X_2 \hookrightarrow L^1(0, L_1) \times L^1(0, L_2)$.

\begin{lemm}
\label{LemmUConst}
For every $z \in \mathsf X_2$ and $t \geq 0$, we have $U(e^{tA}z) = U(z)$.
\end{lemm}

\begin{prf}
By the density of $D(A)$ in $\mathsf X_2$ and by the continuity of $U$, it suffices to show this for $z \in D(A)$. In this case, the function $t \mapsto U(e^{tA}z)$ is differentiable in $\mathbb R_+$, and, noting $e^{tA}z = (u_1(t), u_2(t))$, we have by a trivial computation that 
$\frac{d}{dt} U(e^{tA}z) =0$.
\end{prf}

Define the operator $L$ on $\mathsf X_2$ by $L z = (U(z), U(z))$ and notice that $L \in \mathcal L(\mathsf X_2)$. The main result of this section is the following, which proves Theorem~\ref{TheoAsympTwoCircles}.\ref{CaseIrrational} and gives the explicit value of the constant $\lambda$.

\begin{theo}
\label{TheoConvConst}
Suppose ${L_1}/{L_2} \notin \mathbb Q$. Then, for every $z_0 \in \mathsf X_2$, $\lim_{t \to +\infty} e^{tA} z_0 = L z_0$.
\end{theo}

\begin{prf}
Since $L$ is a bounded operator and the semigroup $\{e^{tA}\}_{t \geq 0}$ is uniformly bounded, it suffices by density to show this result for $z_0 \in D(A)$. By Lemma~\ref{LemmOrbitePrecompacte} and Theorem~\ref{TheoOmegaLimLaSalle}.\ref{TheoLaSalle}, we have $\omega(z_0) \subset D(A) \cap M$ and thus, by Lemma~\ref{LemmDAcMConstant}, every function in $\omega(z_0)$ is constant. Let $z = (\lambda, \lambda) \in \omega(z_0)$ with $\lambda \in \mathbb R$ and take $(t_n)_{n \in \mathbb N}$ a nondecreasing sequence in $\mathbb R_+$ with $t_n \to +\infty$ as $n \to \infty$ such that $e^{t_n A}z_0 \to z$ in $\mathsf X_2$ as $n \to \infty$. By the continuity of $U$ and by Lemma~\ref{LemmUConst}, we obtain that
\[\lambda = U(z) = \lim_{n \to \infty}U(e^{t_n A}z_0) = U(z_0)\]
and thus $z = L z_0$. Hence $\omega(z_0) = \{L z_0\}$ and, by definition of $\omega(z_0)$, this shows that $e^{tA} z_0 \to L z_0$ as $t \to +\infty$, which gives the desired result.
\end{prf}


\subsection{Periodic solutions for the undamped system}
\label{SecPeriodic}

We now turn to a constructive proof of Theorem~\ref{TheoAsympTwoCircles}.\ref{CaseRational}.

\begin{prf}[Proof of Theorem~\ref{TheoAsympTwoCircles}.\ref{CaseRational}]
Take $p, q \in \mathbb N^\ast$ coprime such that ${L_1}/{L_2} = {p}/{q}$. Let $\ell = {L_1}/{p} = {L_2}/{q}$. Take $\varphi \in \mathcal C^\infty_c(\mathbb R)$ with support included in $(0, \ell)$. For $x \in [0, L_1]$, we define $u_{1, 0}$ by
\begin{equation}
u_{1, 0}(x) = \sum_{k = -\infty}^{+\infty} \varphi(x - k \ell).
\label{u10Per}
\end{equation}
Notice that, for each $x \in \mathbb R$, there exists at most one $k \in \mathbb Z$ such that $\varphi(x - k \ell) \not = 0$. In particular, $u_{1, 0} \in \mathcal C^\infty([0, L_1])$. Similarly we define $u_{2, 0}\in \mathcal C^\infty([0, L_2])$ by the same expression
\begin{equation}
u_{2, 0}(x) = \sum_{k = -\infty}^{+\infty} \varphi(x - k \ell),\quad x \in [0, L_2].
\label{u20Per}
\end{equation}
Define $u_j(t, x) =u_{j,0}(x-t)$ for $j=1, 2$. Since $L_1 = p \ell$, $L_2 = q \ell$, we have 
\[u_1(t, L_1) = u_2(t, L_2) = u_1(t, 0) = u_2(t, 0).\]
Thus, $(u_1, u_2)$ is the unique solution of \eqref{AppTwoCirclesUndamped} with initial data $z_0 = (u_{1, 0}, u_{2, 0})$. It is periodic in time, and non-constant if $\varphi$ is chosen to be non-constant.
\end{prf}


\subsection{Periodic solutions for the persistently damped system}
\label{SecPEPeriodic}

Sections~\ref{SecIrrational} and \ref{SecPeriodic} present the asymptotic behavior of \eqref{AppTwoCirclesUndamped} in the cases ${L_1}/{L_2} \notin \mathbb Q$ and ${L_1}/{L_2} \in \mathbb Q$, showing that all solutions converge to a constant in the first case and that periodic solutions exist in the second one. When considering System~\eqref{TwoCircles} with a persistent damping, the fact that all its solutions converge exponentially to the origin when ${L_1}/{L_2} \notin \mathbb Q$ is a consequence of our main result, Theorem~\ref{MainTheoIntro}. However, if ${L_1}/{L_2} \in \mathbb Q$ and the damping interval $[a, b]$ is small enough, one may obtain periodic solutions, thus showing that Theorem~\ref{MainTheoIntro} cannot hold in general for ${L_1}/{L_2} \in \mathbb Q$ and any length of damping interval.

\begin{theo}
\label{TheoPEPeriodic}
Suppose that ${L_1}/{L_2} \in \mathbb Q$. Then there exists $\ell_0 > 0$ such that, if $b - a \leq \ell_0$, there exists $\alpha \in \mathcal G(4\ell_0, \ell_0)$ for which \eqref{TwoCircles} admits a non-zero periodic solution.
\end{theo}

\begin{prf}
We consider here the construction of a periodic solution for \eqref{AppTwoCirclesUndamped} done in the proof of Theorem~\ref{TheoAsympTwoCircles}.\ref{CaseRational}. Take $p, q \in \mathbb N^\ast$ coprime such that ${L_1}/{L_2} = {p}/{q}$ and note $\ell = {L_1}/{p} = {L_2}/{q}$. Take $\varphi \in \mathcal C^\infty_c(\mathbb R)$ not identically zero with support included in $(0, \ell/2)$. By proceeding as in  Theorem~\ref{TheoAsympTwoCircles}.\ref{CaseRational} we get a periodic non-zero solution $(u_1,u_2)$ of \eqref{AppTwoCirclesUndamped}
given by
\begin{equation}
u_1(t, x) = \sum_{k = -\infty}^{+\infty} \varphi(x - t - k \ell), \qquad u_2(t, x) = \sum_{k = -\infty}^{+\infty} \varphi(x - t - k \ell).
\label{PeriodicSolution}
\end{equation}

Take $\ell_0 = \nicefrac{\ell}{4}$ and suppose that $b - a \leq \ell_0$. We construct a periodic signal $\alpha: \mathbb R \to \{0, 1\}$ defined by
\[
\alpha(t) = 
\begin{dcases}
0, & \text{if } t \in \bigcup_{n \in \mathbb Z} \left[a - \left(n + \nicefrac{1}{2}\right)\ell, b - n \ell\right], \\
1, & \text{otherwise.}
\end{dcases}
\]
This defines a periodic signal $\alpha$ with period $T = \ell = 4\ell_0$
which belongs to $\mathcal G(4\ell_0,\ell_0)$. One then easily checks that $\alpha(t) \chi(x) u_2(t, x) = 0$ for every $(t, x) \in \mathbb R_+ \times [0, L_2]$, and so \eqref{PeriodicSolution} satisfies \eqref{TwoCircles}.
\end{prf}


\section{Proof of Theorem~\ref{TheoSolExplicite}}
\label{AppendExplicit}

From now on, we use the convention $\beta^{(i)}_{j, \mathfrak n} = 0$ if $\mathfrak n = (n_1, \dotsc, n_N) \in \mathbb Z^N$ is such that $n_k < 0$ for some index $k \in \llbracket 1, N\rrbracket$, so that \eqref{BetaRecurrence} can be written as
\[
\beta^{(i)}_{j, \mathfrak n} = \sum_{k=1}^N m_{kj} \beta^{(i)}_{k, \mathfrak n - \mathbf 1_k}.
\]
One then gets by a trivial induction the following result. 

\begin{lemm}
\label{LemmBetaSecond}
For every $i, j \in \llbracket 1, N \rrbracket$ and $\mathfrak n = (n_1, \dotsc, n_N) \in \mathfrak N \setminus \{\mathbf 0\}$, we have
\begin{equation}
\beta^{(i)}_{j, \mathfrak n} = \sum_{k=1}^N m_{ik} \beta^{(k)}_{j, \mathfrak n - \mathbf 1_k}.
\label{BetaSecondRecurrence}
\end{equation}
\end{lemm}

We can now turn to the proof of Theorem~\ref{TheoSolExplicite}.

\begin{prf}[Proof of Theorem~\ref{TheoSolExplicite}]
Let $z_0 = (u_{1, 0}, \dotsc, u_{N, 0}) \in D(A)$ and let $z = (u_1, \dotsc, u_N)$ be defined by \eqref{ExplicitSolution}, with $u_i(t, 0)$ given by \eqref{ExplicitPhi}. Notice that $u_i(\cdot, 0)$ is defined everywhere on $\mathbb R_+$ and is measurable, so that $u_i$ is defined everywhere on $\mathbb R_+ \times [0, L_i]$ and is measurable. Note also that $u_i$ is well defined, since
\[
u_i(0, 0) = \sum_{j=1}^N \beta^{(i)}_{j, \mathbf 0} u_{j, 0}(L_j) = u_{i, 0}(0)
\]
thanks to the fact that $\beta^{(i)}_{j, \mathbf 0} = m_{ij}$ and $z_0 \in D(A)$.

Let $T_0 > 0$ be as in the proof of Theorem~\ref{TheoExistUndamped}. The unique solution of \eqref{Undamped} with initial condition $z_0$ is then given by \eqref{EqSol} for $0 \leq t \leq T_0$, and, in order to prove the theorem, it suffices to show that, for every $t,\tau$ with $0\leq \tau\leq t\leq  \tau + T_0$, we have
\begin{equation}
u_i(t, x) = 
\begin{dcases*}
\sum_{j=1}^N m_{ij} u_j(\tau, L_j - (t - \tau) + x), & if $0 \leq x \leq t - \tau$, \\
u_i(\tau, x - t + \tau), & if $x > t - \tau.$
\end{dcases*}
\label{SolIntervalT0}
\end{equation}
Indeed, if this is proved, we apply it to $\tau = 0$ to obtain that $z$ is the solution of \eqref{Undamped} with initial condition $z_0$ for $0 \leq t \leq T_0$, and also that $z(t) \in D(A)$ for every $t \in [0, T_0]$, and a simple induction allows us to conclude.

We prove \eqref{SolIntervalT0} by considering three cases.

{\bf Case 1.} $0 \leq t - \tau < x$ and $t \leq x$.
By \eqref{ExplicitSolution}, we have $u_i(t, x) = u_{i, 0}(x - t) = u_i(\tau, x - t + \tau)$.

{\bf Case 2.}  $0 \leq t - \tau < x$ and $t > x$. Since $x - t + \tau < \tau$, we have in this case $u_i(\tau, x - t + \tau) = u_i(\tau - x + t - \tau, 0) = u_i(t - x, 0) = u_i(t, x)$.

{\bf Case 3.} $t - \tau \geq x$.
We notice first that it suffices to consider the case $x = 0$. Indeed, if $t - \tau \geq x$, then $t \geq x$, and it follows clearly by the definition \eqref{ExplicitSolution} of $u_i$ that $u_i(t, x) = u_i(t - x, 0)$ for $t \geq x \geq 0$. On the other hand, let us denote
\[
v_{i, \tau}(t, x) = 
\begin{dcases*}
\sum_{j=1}^N m_{ij} u_j(\tau, L_j - (t - \tau) + x) & if $0 \leq x \leq t - \tau$, \\
u_i(\tau, x - t + \tau) & if $x > t - \tau$.
\end{dcases*}
\]
It is also clear that $v_{i, \tau}(t, x) = v_{i, \tau}(t - x, 0)$, and thus it suffices to show that $u_i(t, 0) = v_{i, \tau}(t, 0)$ for every $t \in [\tau, \tau + T_0]$ in order to conclude that $u_i(t, x) = u_i(t-x, 0) = v_{i, \tau}(t - x, 0) = v_{i, \tau}(t, x)$ for every $t \in [\tau, \tau + T_0]$ and $x \in [0, L_i]$ with $t - \tau \geq x$.

For $t \in [\tau, \tau + T_0]$, we have $v_{i, \tau}(t, 0) = \sum_{j=1}^N m_{ij} u_j(\tau, L_j - t + \tau)$.
Furthermore
\[
u_j(\tau, L_j - t + \tau) = 
\begin{dcases*}
u_{j, 0}(L_j - t), & if $t \leq L_j$, \\
u_j(t - L_j, 0), & if $t \geq L_j$,
\end{dcases*}
\]
and so
\begin{align*}
v_{i, \tau}(t, 0) & = \sum_{j=1}^N m_{ij} u_j(\tau, L_j - t + \tau) \displaybreak[0]
  = \sum_{\substack{j=1 \\ L_j \leq t}}^N m_{ij} u_j(t - L_j, 0) + \sum_{\substack{j=1 \\ L_j > t}}^N m_{ij} u_{j, 0}(L_j - t) \displaybreak[0]\\%
 & = \sum_{\substack{j=1 \\ L_j \leq t}}^N m_{ij} \sum_{k=1}^N \sum_{\substack{\mathfrak n \in \mathfrak N_k \\ L(\mathfrak n) \leq t - L_j}} \beta^{(j)}_{k, \mathfrak n + \floor{\frac{t - L_j - L(\mathfrak n)}{L_k}} \mathbf 1_k} u_{k, 0}(L_k - \{t - L_j - L(\mathfrak n)\}_{L_k}) \\
 & \hphantom{=} {} + \sum_{\substack{j=1 \\ L_j > t}}^N m_{ij} u_{j, 0}(L_j - t) \displaybreak[0]\\%
 & = \sum_{\substack{k=1 \\ L_k \leq t}}^N \sum_{\substack{j=1 \\ L_j \leq t}}^N \sum_{\substack{\mathfrak n \in \mathfrak N_k \\ L(\mathfrak n) \leq t - L_j}} m_{ij} \beta^{(j)}_{k, \mathfrak n + \floor{\frac{t - L_j - L(\mathfrak n)}{L_k}} \mathbf 1_k} u_{k, 0}(L_k - \{t - L_j - L(\mathfrak n)\}_{L_k}) \\
 &+ \sum_{\substack{k=1 \\ L_k > t}}^N \sum_{\substack{j=1 \\ L_j \leq t}}^N \sum_{\substack{\mathfrak n \in \mathfrak N_k \\ L(\mathfrak n) \leq t - L_j}} m_{ij} \beta^{(j)}_{k, \mathfrak n} u_{k, 0}(L_k - (t - L_j - L(\mathfrak n))) 
 + \sum_{\substack{j=1 \\ L_j > t}}^N \beta^{(i)}_{j, \mathbf 0} u_{j, 0}(L_j - t).
\end{align*}
Set
\[A_1(t) = \{(k, j, \mathfrak n) \in \llbracket 1, N \rrbracket \times \llbracket 1, N \rrbracket \times \mathfrak N \suchthat L_k \leq t, \; L_j \leq t, \; \mathfrak n \in \mathfrak N_k, \; L(\mathfrak n) \leq t - L_j\},\]
\[A_2(t) = \{(k, j, \mathfrak n) \in \llbracket 1, N \rrbracket \times \llbracket 1, N \rrbracket \times \mathfrak N \suchthat L_k > t, \; L_j \leq t, \; \mathfrak n \in \mathfrak N_k, \; L(\mathfrak n) \leq t - L_j\},\]
so that we can write
\begin{equation}
\begin{aligned}
v_{i, \tau}(t, 0) & = \sum_{(k, j, \mathfrak n) \in A_1(t)} m_{ij} \beta^{(j)}_{k, \mathfrak n + \floor{\frac{t - L_j - L(\mathfrak n)}{L_k}} \mathbf 1_k} u_{k, 0}(L_k - \{t - L_j - L(\mathfrak n)\}_{L_k}) \\
 & \hphantom{=} {} + \sum_{(k, j, \mathfrak n) \in A_2(t)} m_{ij} \beta^{(j)}_{k, \mathfrak n} u_{k, 0}(L_k - (t - L_j - L(\mathfrak n))) 
  + \sum_{\substack{j=1 \\ L_j > t}}^N \beta^{(i)}_{j, \mathbf 0} u_{j, 0}(L_j - t).
\end{aligned}
\label{viInFunctionofAi}
\end{equation}

Set
\begin{multline*}
B_1(t) = \left\{(k, j, \mathfrak m) \in \llbracket 1, N \rrbracket \times \llbracket 1, N \rrbracket \times \mathfrak N \middlesuchthat L_k \leq t, \; \mathfrak m \in \mathfrak N_k,\right. \\ \left. L(\mathfrak m) \leq t, \; \mathfrak m + \floor{\frac{t - L(\mathfrak m)}{L_k}} \mathbf 1_k = (r_1, \dotsc, r_N) \text{ with } r_j \geq 1\right\},
\end{multline*}
\begin{multline*}
B_2(t) = \{(k, j, \mathfrak m) \in \llbracket 1, N \rrbracket \times \llbracket 1, N \rrbracket \times \mathfrak N \suchthat L_k > t, \\ \mathfrak m = (m_1, \dotsc, m_N) \in \mathfrak N_k \setminus\{\mathbf 0\}, \; m_j \geq 1, \; L(\mathfrak m) \leq t\},
\end{multline*}
and define the functions $\varphi_\lambda: \llbracket 1, N \rrbracket \times \llbracket 1, N \rrbracket \times \mathfrak N \to \llbracket 1, N \rrbracket \times \llbracket 1, N \rrbracket \times \mathfrak N$, $\lambda = 1, 2$, by
\[
\varphi_1(k, j, \mathfrak n) = 
\begin{dcases*}
(k, j, \mathfrak n + \mathbf 1_j) & if $k \not = j$, \\
(k, j, \mathfrak n) & if $k = j$,
\end{dcases*}
\qquad
\varphi_2(k, j, \mathfrak n) = (k, j, \mathfrak n + \mathbf 1_j).
\]
We claim that $\varphi_\lambda$ is a bijection from $A_\lambda(t)$ to $B_\lambda(t)$, $\lambda = 1, 2$. Indeed, it is easy to verify that the image of $A_\lambda(t)$ by $\varphi_\lambda$ is included in $B_\lambda(t)$ and that $\varphi_\lambda: A_\lambda(t) \to B_\lambda(t)$ is injective for $\lambda = 1, 2$. Let us check that these functions are surjective.

If $(k, j, \mathfrak m) \in B_1(t)$, we note $(r_1, \dotsc, r_N) = \mathfrak m + \floor{{(t - L(\mathfrak m))}/{L_k}} \mathbf 1_k$ and we set $\mathfrak n = \mathfrak m$ if $k = j$ and $\mathfrak n = \mathfrak m - \mathbf 1_j$ if $k \not = j$. Notice first that, if $k \not = j$, then $m_j = r_j \geq 1$, so that $\mathfrak m - \mathbf 1_j \in \mathfrak N$, and thus, in both cases $k = j$ and $k \not = j$, we have $(k, j, \mathfrak n) \in \llbracket 1, N \rrbracket \times \llbracket 1, N \rrbracket \times \mathfrak N$, and clearly $\varphi_1(k, j, \mathfrak n) = (k, j, \mathfrak m)$, so that, in order to conclude that $\varphi_1: A_1(t) \to B_1(t)$ is surjective, it suffices to show that $(k, j, \mathfrak n) \in A_1(t)$. We clearly have $L_k \leq t$ and $\mathfrak n \in \mathfrak N_k$. If $j = k$, we have $L_j = L_k \leq t$ and, since $r_j = r_k = \floor{{(t - L(\mathfrak m))}/{L_k}}$ and $r_j \geq 1$, we have ${(t - L(\mathfrak m))}/{L_k} \geq 1$, i.e., $L(\mathfrak n) = L(\mathfrak m) \leq t - L_k = t - L_j$, so that $(k, j, \mathfrak n) \in A_1(t)$. If $j \not = k$, we have $m_j = r_j \geq 1$, so that $L_j \leq m_j L_j \leq L(\mathfrak m) \leq t$; also, $L(\mathfrak n) = L(\mathfrak m) - L_j \leq t - L_j$, so that $(k, j, \mathfrak n) \in A_1(t)$. Hence $(k, j, \mathfrak n) \in A_1(t)$ in both cases, and thus $\varphi_1: A_1(t) \to B_1(t)$ is surjective.

If $(k, j, \mathfrak m) \in B_2(t)$, we set $\mathfrak n = \mathfrak m - \mathbf 1_j$, so that $\mathfrak n \in \mathfrak N$ and $\varphi_2(k, j, \mathfrak n) = (k, j, \mathfrak m)$. Now, it is clear that $L_k > t$ and $\mathfrak n \in \mathfrak N_k$, and we have $L_j \leq m_j L_j \leq L(\mathfrak m) \leq t$ and $L(\mathfrak n) = L(\mathfrak m) - L_j \leq t - L_j$. Hence $(k, j, \mathfrak n) \in A_2(t)$, and thus $\varphi_2: A_2(t) \to B_2(t)$ is surjective.

Thanks to the bijections $\varphi_\lambda: A_\lambda(t) \to B_\lambda(t)$, $\lambda = 1, 2$, we can rewrite \eqref{viInFunctionofAi} as
\begin{equation*}
\begin{aligned}
v_{i, \tau}(t, 0) & = \sum_{(k, j, \mathfrak m) \in B_1(t)} m_{ij} \beta^{(j)}_{k, \mathfrak m - \mathbf 1_j + \floor{\frac{t - L(\mathfrak m)}{L_k}} \mathbf 1_k} u_{k, 0}(L_k - \{t - L(\mathfrak m)\}_{L_k}) \\
 & \hphantom{=} {} + \sum_{(k, j, \mathfrak m) \in B_2(t)} m_{ij} \beta^{(j)}_{k, \mathfrak m - \mathbf 1_j} u_{k, 0}(L_k - (t - L(\mathfrak m))) 
+ \sum_{\substack{j=1 \\ L_j > t}}^N \beta^{(i)}_{j, \mathbf 0} u_{j, 0}(L_j - t),
\end{aligned}
\end{equation*}
and so, by applying Lemma~\ref{LemmBetaSecond}, we obtain
\begin{align*}
v_{i, \tau}(t, 0) & = \sum_{\substack{k = 1 \\ L_k \leq t}}^N \sum_{\substack{\mathfrak m \in \mathfrak N_k \\ L(\mathfrak m) \leq t}} \sum_{\substack{j = 1 \\ \left(\mathfrak m + \floor{\frac{t - L(\mathfrak m)}{L_k}} \mathbf 1_k\right)_j \geq 1}}^N m_{ij} \beta^{(j)}_{k, \mathfrak m + \floor{\frac{t - L(\mathfrak m)}{L_k}} \mathbf 1_k - \mathbf 1_j} u_{k, 0}(L_k - \{t - L(\mathfrak m)\}_{L_k}) \\
 & \hphantom{=} {} + \sum_{\substack{k = 1 \\ L_k > t}}^N \sum_{\substack{\mathfrak m \in \mathfrak N_k \setminus \{\mathbf 0\} \\ L(\mathfrak m) \leq t}} \sum_{\substack{j=1 \\ m_j \geq 1}}^N m_{ij} \beta^{(j)}_{k, \mathfrak m - \mathbf 1_j} u_{k, 0}(L_k - (t - L(\mathfrak m))) 
+ \sum_{\substack{k=1 \\ L_k > t}}^N \beta^{(i)}_{k, \mathbf 0} u_{k, 0}(L_k - t) \displaybreak[0]\\%
 & = \sum_{\substack{k = 1 \\ L_k \leq t}}^N \sum_{\substack{\mathfrak m \in \mathfrak N_k \\ L(\mathfrak m) \leq t}} \beta^{(i)}_{k, \mathfrak m + \floor{\frac{t - L(\mathfrak m)}{L_k}} \mathbf 1_k} u_{k, 0}(L_k - \{t - L(\mathfrak m)\}_{L_k}) \\
 & \hphantom{=} {} + \sum_{\substack{k = 1 \\ L_k > t}}^N \sum_{\substack{\mathfrak m \in \mathfrak N_k \setminus \{\mathbf 0\} \\ L(\mathfrak m) \leq t}} \beta^{(i)}_{k, \mathfrak m} u_{k, 0}(L_k - (t - L(\mathfrak m))) + \sum_{\substack{k=1 \\ L_k > t}}^N \beta^{(i)}_{k, \mathbf 0} u_{k, 0}(L_k - t) \displaybreak[0]\\%
 & = \sum_{k = 1}^N \sum_{\substack{\mathfrak m \in \mathfrak N_k \\ L(\mathfrak m) \leq t}} \beta^{(i)}_{k, \mathfrak m + \floor{\frac{t - L(\mathfrak m)}{L_k}} \mathbf 1_k} u_{k, 0}(L_k - \{t - L(\mathfrak m)\}_{L_k}) = u_i(t, 0).
\end{align*}
Hence $v_{i, \tau}(t, 0) = u_i(t, 0)$ for every $t \in [\tau, \tau + T_0]$, which thus concludes the proof of \eqref{SolIntervalT0}.
\end{prf}


\section{A combinatorial estimate}

In order to estimate the right-hand side of \eqref{BinomialBound}, one needs the following lemma.

\begin{lemm}
\label{LemmCoeffBinom}
Let $\nu \in (0, 1)$. There exist $\rho \in \left(0, \nicefrac{1}{2}\right)$, $C, \gamma > 0$ such that, for every $n \in \mathbb N$ and $k \in \llbracket 0, \rho n \rrbracket$, we have
\begin{equation}
\binom{n}{k} \nu^{n} \leq C e^{-\gamma n}.
\label{EstimCoeffBinom}
\end{equation}
\end{lemm}

\begin{prf}
For $n \in \mathbb N$, consider the function $f_n(k) = \binom{n}{k} \nu^{n}$
defined for $k \in \llbracket 0, n\rrbracket$. Since $k \mapsto \binom{n}{k}$ is increasing for $k \in \left\llbracket 0, \nicefrac{n}{2}\right\rrbracket$, if $k \in \llbracket 0, \rho n \rrbracket$ for a certain $\rho \in \left(0, \nicefrac{1}{2}\right)$, then
\begin{equation}
f_n(k) \leq f_n(\floor{\rho n}).
\label{fnkleqfnrho}
\end{equation}

Let us estimate $f_n(\floor{\rho n})$ for $n$ large. As $n \to +\infty$, using Stirling's approximation $\log n! = n \log n - n + O(\log n)$, we get
\[
\begin{aligned}
\log f_n(\floor{\rho n}) & = \log\binom{n}{\floor{\rho n}} + n \log \nu \\
 & = n \log n - \rho n \log (\rho n) - (1 - \rho) n \log [(1-\rho) n] + n \log \nu + O(\log n) \\
 & = n \left[g(\rho) + O\left(\frac{\log n}{n}\right)\right],
\end{aligned}
\]
where the function $g: (0, 1) \to \mathbb R$ is defined by
\[g(\rho) = \rho \log\left(\frac{1}{\rho}\right) + (1-\rho) \log\left(\frac{1}{1-\rho}\right) + \log \nu.\]
It is a continuous function of $\rho \in (0, 1)$ and $\lim_{\rho \to 0} g(\rho) = \log \nu < 0$, hence there exists $\rho \in \left(0, \nicefrac{1}{2}\right)$ depending only on $\nu$ such that $g(\rho) \leq \frac{1}{2} \log\nu < 0$. For this value of $\rho$, we have
\[
\log f_n(\floor{\rho n}) \leq \frac{n}{2} \log\nu + O(\log n) \leq \frac{n}{4} \log\nu + O(1).
\]
The result follows by combining the above with \eqref{fnkleqfnrho}.
\end{prf}

\bibliographystyle{abbrv}
\bibliography{Bib}

\end{document}